\documentclass[12pt]{amsart}

\usepackage{ae,lmodern} 
\usepackage[utf8]{inputenc} 
\usepackage[T1]{fontenc} 
\usepackage[english]{babel}

\usepackage{amsthm}
\usepackage{geometry}
\usepackage{dsfont}
\usepackage{MnSymbol}
\usepackage{comment}
\usepackage[hidelinks]{hyperref}
\usepackage{stmaryrd}

\usepackage[
   backend=biber,        
   sorting=nyt,          
   citestyle=alphabetic,
   bibstyle=alphabetic,  
   isbn=false, 
   url=false, 
   doi=false, 
   eprint=false,
]{biblatex}
\AtEveryBibitem{%
  \clearfield{note}%
}

\usepackage{indentfirst}

\let\oldtocsection=\tocsection
\let\oldtocsubsection=\tocsubsection
\let\oldtocsubsubsection=\tocsubsubsection

\renewcommand{\tocsection}[2]{\hspace{0em}\oldtocsection{#1}{#2}}
\renewcommand{\tocsubsection}[2]{\hspace{1em}\oldtocsubsection{#1}{#2}}
\renewcommand{\tocsubsubsection}[2]{\hspace{2em}\oldtocsubsubsection{#1}{#2}}

\newtheorem{Theorem}{Theorem}[section]
\newtheorem{Corollary}[Theorem]{Corollary}
\newtheorem{Lemma}[Theorem]{Lemma}
\newtheorem{Proposition}[Theorem]{Proposition}
\newtheorem*{Question*}{Question}

\theoremstyle{remark}
\newtheorem{Remark}[Theorem]{Remark}

\theoremstyle{definition}
\newtheorem{Definition}{Definition}[section]

\newcommand{\p}[2][]{%
\mathsf{p}_{#1}({#2})}
\newcommand{\phat}[2][]{%
\hat{\mathsf{p}}_{#1}({#2})}
\newcommand{\s}[2][]{%
\mathsf{s}_{#1}({#2})}
\newcommand{\sm}[2][]{%
\mathsf{sw}_{#1}({#2})}
\newcommand{\cycl}[2][]{%
\overset{\curvearrowleft #1}{#2}}
\newcommand{\N}{\mathbb{N}}
\newcommand{\Z}{\mathbb{Z}}
\newcommand{\Q}{\mathbb{Q}}
\newcommand{\F}{\mathbb{F}}
\newcommand{\R}{\mathbb{R}}
\newcommand{\C}{\mathbb{C}}

\newcommand{\HH}{\mathbb{H}}
\newcommand{\eps}{\varepsilon}

\newcommand{\lf}{\lfloor}
\newcommand{\rf}{\rfloor}

\newcommand{\la}{\langle}
\newcommand{\ra}{\rangle}

\newcommand{\slope}{\mathrm{Slope}}

\title[Equivalence of primitive-stable and Bowditch actions]{Equivalence of primitive-stable and Bowditch actions of the free group of rank two on Gromov-hyperbolic spaces}
\author{Suzanne Schlich}

\begin{document}

\renewcommand{\labelitemi}{$\bullet$}

\newgeometry{hmargin=2cm,vmargin=2cm}
\maketitle

\begin{abstract}
We prove that the set of Bowditch representations (introduced by Bowditch in 1998, then generalized by Tan, Wong and Zhang in 2008) and the set of primitive-stable representations (introduced by Minsky in 2013) of the free group of rank two in the isometry group of a Gromov-hyperbolic space are equal. The case of $\mathrm{PSL}(2,\mathbb{C})$-representations has already been proved by Lee and Xu and independently Series. Our proof in this context is independent.
\end{abstract}

\tableofcontents
\restoregeometry

\section{Introduction}

When $\Gamma$ is a finitely generated group and $G$ a topological group, we can consider the space $\mathrm{Hom}(\Gamma,G)$ of representations from $\Gamma$ to $G$ endowed with the compact-open topology. The group $G$ acts on $\mathrm{Hom}(\Gamma,G)$ by conjugation and gives rise to a quotient $\mathrm{Hom}(\Gamma,G)/G$. Since this quotient might not be Hausdorff, we can instead consider its Hausdorffization, that is the largest Hausdorff quotient of $\mathrm{Hom}(\Gamma,G)/G$, which is called the character variety and denoted $\chi(\Gamma,G)=\mathrm{Hom}(\Gamma,G)//G$. Note that in the usual case where $G$ is an algebraic Lie group, the character variety is endowed with the structure of an algebraic variety, which justifies the terminology. Let $\mathrm{Out}(\Gamma)=\mathrm{Aut}(\Gamma)/\mathrm{Inn}(\Gamma)$ be the outer automorphism group of $\Gamma$, that is the quotient of the automorphism group of $\Gamma$ by the inner automorphisms. In this context, there is a natural action of $\mathrm{Out}(\Gamma)$ on the character variety by precomposition as follows :
\begin{equation*}
\begin{array}{rll}
    \mathrm{Out}(\Gamma) \times \chi(\Gamma,G) & \longrightarrow & \chi(\Gamma,G) \\
     ([\phi],[\rho]) & \longmapsto & [\phi].[\rho]=[\rho \circ \phi^{-1}]
\end{array}
\end{equation*}
This general setting has a natural geometric interpretation when $\Gamma=\pi_1(S)$ is the fundamental group of~$S$, a closed oriented surface of negative Euler characteristic and $G=\mathrm{PSL}(2,\R)$. The Teichmüller space of $S$, $\mathcal{T}(S)$, which is the space of all hyperbolic structures on $S$, is known to identify as a connected component of the character variety $\chi(\pi_1(S),\mathrm{PSL}(2,\R))$ (see Goldman \cite{goldman_topological_1988}, Farb-Margalit \cite{farb_primer_2012}). Furthermore, the mapping class group $\mathrm{MCG}(S)$ of $S$, which is the group of isotopy classes of orientation-preserving homeomorphisms of $S$, may be identified with an index two subgroup of $\mathrm{Out}(\pi_1(S))$. The Fricke's Theorem then asserts that the mapping class group acts on the Teichmüller space properly discontinuously. A bit more generally, for $\Gamma$ a hyperbolic group and $G=\mathrm{SO}_0(n,1)$ the isometry group of the hyperbolic $n$-space $\HH^n$, we can consider the space $\mathrm{CC}(\Gamma,G)$ of convex-cocompact representations, defined as the representations which preserve a non-empty convex subset of $\HH^n$ on which the action is properly discontinuous and cocompact. This is equivalent to asking that the orbit map of the representation is a quasi-isometric embedding. As for the Teichmüller space, it can be shown that $\mathrm{CC}(\Gamma,G)$ is open and that the outer automorphism group $\mathrm{Out}(\Gamma)$ acts properly discontinuously on it. Thus, a natural question is whether this domain of discontinuity is maximal or not in the character variety. 

\subsection{Bowditch representations}
In this paper, we will study the case where $\Gamma=\F_2$ is the free group of rank two and $G=\mathrm{Isom}(X)$, where $X$ is a $\delta$-hyperbolic, geodesic and visibility space. A \emph{geodesic} space means that every two points in $X$ can be connected by a geodesic segment, and \emph{visibility} space that every two distinct points in $\partial X$, the boundary of $X$, are the endpoints of some bi-infinite geodesic in $X$. \\

We fix $\{a,b\}$ a free generating set for $\F_2$. We denote by $|\gamma|$ the word length of an element $\gamma$ in $\F_2$ relatively to the set of two generators $\{a,b\}$ and $\Vert \gamma \Vert $ its cyclically reduced word length. Recall that a different choice of the set of generators would give a bi-Lipschitz equivalent (cyclically reduced) word length. When given an isometry $A$ of the metric space $X$, we can consider its \emph{displacement length}, that is the non-negative real $l(A):=\underset{x \in X} \inf d(Ax,x)$. For a representation $\rho : \F_2 \to \mathrm{Isom}(X)$, we will denote by $l_\rho(\gamma):=l(\rho(\gamma))$ the displacement length of $\rho(\gamma)$, for $\gamma \in \F_2$. We can always compare, using the triangle inequality, the displacement length and the cyclically reduced word length in the following way :  
\begin{equation} \label{majo-up}
    \forall \gamma \in \F_2, \hspace{0,5cm} l_\rho(\gamma)\leq C' \Vert \gamma \Vert 
\end{equation}
Here, $C'$ is a constant that can be chosen to be the maximum $C':=\max \{d(\rho(a)o,o),d(\rho(b)o,o)\}$, where $\{a,b\}$ is a free generating set for $\F_2$ and $o$ any basepoint of $X$. \\ 
On the other hand, the converse inequality, asking whether or not there exists two constants $C$ and $D$, such that for any element $\gamma \in \F_2$, we have : $\displaystyle \frac{1}{C} \Vert \gamma \Vert -D \leq l_\rho(\gamma)$, is of course not always satisfied. This condition precisely gives rise to the notion of convex-cocompact representations already mentioned above. \\

In the free group $\F_2$, we say that an element $\gamma \in \F_2$ is \emph{primitive} if it is part of a free basis of $\F_2$. We denote by $\mathcal{P}(\F_2)$ the set of  primitive elements in $\F_2$. \\
Following the work of Bowditch (\cite{bowditch_markoff_1998}), a broader class of representations can be obtained when considering this last inequality only for primitive elements in $\F_2$, and this leads to the notion of what we call a  Bowditch representation. 
\begin{Definition}[Bowditch representations] \label{def-bowditch}
Let $\rho : \mathbb{F}_2 \to \mathrm{Isom}(X)$ be a representation and $C \geq 1, D \geq 0$ two constants. We say that $\rho$ is a \emph{Bowditch representation of constants $(C,D)$} if :
\begin{equation}
    \forall \gamma \in \mathcal{P}(\mathbb{F}_2), \hspace{0,5cm} \frac{1}{C}\Vert \gamma \Vert -D \leq l_\rho(\gamma) 
\end{equation}
We say that $\rho$ is a \emph{Bowditch representation} if there exist two constants $C \geq 1$ and $D \geq 0$ such that $\rho$ is a Bowditch representation of constants $(C,D)$.
\end{Definition}
Denote by $\mathcal{BQ}(\F_2,X)$ the set of Bowditch representations from $\F_2$ to $\mathrm{Isom}(X)$. \\ 

The original definition by Bowditch (\cite{bowditch_markoff_1998}) was given for representations of $\F_2$ into $\mathrm{PSL}(2,\C)$. He defined these representations using the traces in $\mathrm{PSL}(2,\C)$ of the image of primitive elements, as follows :
\begin{enumerate}
    \item \label{trace-commutateur} $\mathrm{Tr}(\rho([a,b]))=-2$, where $[a,b]$ denotes the commutator of two generators $a$ and $b$ of $\F_2$.
    \item \label{all-hyperbolic}  For all $\gamma \in \mathcal{P}(\F_2), \, \, \, \mathrm{Tr}(\rho(\gamma)) \notin [-2,2]$ 
    \item  \label{set-finite} The set $\{\gamma \in \mathcal{P}(\F_2) \, \, : \, \, |\mathrm{Tr}(\rho(\gamma))| \leq 2 \}$ is finite. 
\end{enumerate}

Bowditch defines $\mathcal{BQ}$ to be the space of representations of $\F_2$ into $\mathrm{PSL}(2,\C)$ (modulo conjugation) satisfying the three previous conditions.  Bowditch shows (\cite{bowditch_markoff_1998}, Theorem 3.16), that $\mathcal{BQ}$ is open in the character variety $\chi(\F_2,\mathrm{PSL}(2,\C))$ and that the outer automorphism group of $\F_2$, $\mathrm{Out}(\F_2)$, acts properly discontinuously on $\mathcal{BQ}$. Hence $\mathcal{BQ}$ produces an open domain of discontinuity for the action of $\mathrm{Out}(\F_2)$ on $\chi(\F_2,\mathrm{PSL}(2,\C))$. In addition, Bowditch shows that (\cite{bowditch_markoff_1998}, Theorem 2) :
\begin{equation} \label{ineq-BQ}
 \text{for all } \rho \in \mathcal{BQ}, \quad \text{ there exists a constant } C >0 \text{ such that for all } \gamma \in \mathcal{P}(\F_2), \quad \frac{1}{C}\Vert\gamma\Vert  \leq l_\rho(\gamma). \end{equation} 
 Note that in this inequality, $l_\rho(\gamma)$ makes sense because $\rho(\gamma)$ is an isometry of the usual hyperbolic space of dimension 3 (recall that $\mathrm{PSL(2,\C)}=\mathrm{Isom}^+(\HH^3)$).  \\
 The work of Bowditch was later on generalized by Tan, Wong and Zhang in \cite{tan_generalized_2008} for representations satisfying only conditions \ref{all-hyperbolic}. and \ref{set-finite}. (that is when $\mathrm{Tr}(\rho([a,b]))=\tau$, for any $\tau \in \C$). They showed, as in the case of Bowditch where \ref{trace-commutateur}. also holds ($\mathrm{Tr}(\rho([a,b]))=-2$), that these representations form an open subspace of the character variety on which the outer automorphism group acts properly discontinuously. They also showed that the inequality \eqref{ineq-BQ} holds. It is easy to check that the converse is also true : when a representation satisfy \eqref{ineq-BQ}, it automatically satisfy conditions \ref{all-hyperbolic}. and \ref{set-finite}. Thus, the  inequality \eqref{ineq-BQ}, which does not make use anymore of the traces of elements in~$\mathrm{PSL}(2,\C)$, can thus be generalised to $\mathrm{Isom}(X)$ as done in definition \ref{def-bowditch}. Note that the additive constant in definition~\ref{def-bowditch} plays no major role. 

\subsection{Primitive-stable representations}

Consider $\mathcal{C}$ the Cayley graph of the free group of rank two $\F_2$ with respect to the free generating set $\{a,b\}$ chosen in previous section. This graph comes equipped with the word metric, that we denote again $d$ (in context, there should be no ambiguity with the metric $d$ of the metric space $X$) and which satisfies : $d(u,v)=|u^{-1}v|$, where $|\cdot|$ is the word length. We will sometimes refer to the vertices of $\mathcal{C}$ as the integer points of $\mathcal{C}$. For $\gamma \in \F_2$, we denote by $L_\gamma$ the axis of $\gamma$ in the Cayley graph $\mathcal{C}$. Note that when $\gamma$ is cyclically reduced, $L_\gamma=\underset{n \in \Z}{\bigcup} [\gamma^n,\gamma^{n+1}]$. We will refer to the geodesics $L_\gamma$ with $\gamma \in \mathcal{P}(\F_2)$ as \emph{primitive-leaves}. \\

Fix $o$ a basepoint in $X$. For every representation $\rho$ of $\F_2$ in $\mathrm{Isom}(X)$ we define the orbit map $\tau_\rho$ of $\rho$ to be the unique $\rho$-equivariant map from the Cayley graph $\mathcal{C}$ of $\F_2$ into $X$ such that $\tau_\rho(1)=o$ and each edge of $\mathcal{C}$ is mapped to a geodesic segment in $X$. Thus we have that $\tau_\rho$ is continuous and that for any vertex~$g \in \F_2, \, \tau_\rho(g)=\rho(g)o$. Moreover, $\tau_\rho$ is Lipschitz, with Lipschitz constant $C'$, where $C'$ can be chosen to be the maximum $C':=\max \{d(\rho(a)o,o),d(\rho(b)o,o)\}$ (with $\{a,b\}$ the free generating set for $\F_2$ used to define the Cayley graph $\mathcal{C}$). This is, as for the inequality \eqref{majo-up}, a consequence of the triangle inequality. \\

Before stating the definition of primitive-stable representations, we recall the definition of a quasi-isometric embedding : 

\begin{Definition} \label{quasi-iso-embedding}
Let $(\mathcal{X},d_{\mathcal{X}})$ and $(\mathcal{Y},d_{\mathcal{Y}})$ be two metric spaces and $C > 0, D\geq 0$ be two constants. \\ 
    We say that a map $f : \mathcal{X} \to \mathcal{Y}$ is a $(C,D)$-\emph{quasi-isometric embedding} if for all points $x$ and $x'$ in~$\mathcal{X}$, we have : 
    \begin{equation} \label{quasi-iso-emb}
        \frac{1}{C}d_{\mathcal{X}}(x,x') -D \leq d_{\mathcal{Y}}(f(x),f(x')) \leq Cd_{\mathcal{X}}(x,x') +D 
    \end{equation}
\end{Definition}

We now give the definition of a primitive-stable representation, as introduced by Minsky in \cite{minsky_dynamics_2013}. 

\begin{Definition} \label{def-primitive-stable}
Let $\rho : \F_2 \to \mathrm{Isom}(X)$ be a representation. We say that $\rho$ is \emph{primitive-stable} if there exist two constants $C \geq 1$ and $D \geq 0$ such that for any primitive element $\gamma \in \F_2$, the orbit map $\tau_\rho$ restricted to $L_\gamma$ is a $(C,D)$-quasi-isometric embedding. 
\end{Definition}

Minsky defined primitive-stability for representations with value in $\mathrm{PSL}(2,\C) = \mathrm{Isom}^+(\HH^3)$, but his definition generalizes directly to the more general $\delta$-hyperbolic case. Note that since, as mentioned above, the orbit map $\tau_\rho$ is always Lipschitz, showing that the orbit map is a quasi-isometric embedding (either on the Cayley graph or on primitive-leaves) reduced to showing the left inequality of \eqref{quasi-iso-emb} : $\displaystyle \frac{1}{C}d(x,x')-D \leq d(\tau_\rho(x),\tau_\rho(x'))$. Also notice that the primitive-stability condition only needs to be verified on cyclically reduced primitive elements. \\

Minsky proved, in \cite{minsky_dynamics_2013}, that the set of primitive-stable representations is open in the character variety $\chi(\F_2,\mathrm{PSL}(2,\C))$, that it is invariant under the action of the outer automorphism group $\mathrm{Out}(F_2)$ and that this action is properly discontinuous. Hence primitive-stable representations provide an open domain of discontinuity for the action of $\mathrm{Out}(\F_2)$ on the character variety $\chi(\F_2,\mathrm{PSL}(2,C))$. Moreover, Minsky proved that the set of primitive-stable representations strictly contains the set of convex-cocompact representations, which is the interior of the set of discrete representations. From these properties follows that there exist non-discrete primitive-stable representations. \\
Denote by $\mathcal{PS}(\F_2,X)$ the set of primitive representation from $\F_2$ to $\mathrm{Isom}(X)$.

\subsection{Equivalence}
Lee and Xu on one hand (\cite{lee_bowditchs_2019}), and Series independently (\cite{series_primitive_2019}, \cite{series_primitive_2020}), proved that the set of Bowditch representations and primitive-stable representations
with values in $\mathrm{PSL}(2,\C)$ are equal. \\

The aim of this paper is to generalise this result to the case of representations in the isometry group of a $\delta$-hyperbolic space. 
Our proof and techniques are independent of those of Series and Lee-Xu for the case $\mathrm{PSL}(2,\C)$.

\begin{Theorem}
\label{BQ=PS} Let $X$ be a $\delta$-hyperbolic, geodesic and visibility space. \\ 
The set of Bowditch representations and primitive-stable representations are equal : 
\begin{equation*}
    \mathcal{BQ}(\F_2,X)=\mathcal{PS}(\F_2,X).
\end{equation*}
\end{Theorem} 

In particular, this theorem applies for $X$ being any rank one symmetric space : the real hyperbolic $n$-space (with isometry group $\mathrm{Isom}(X)=\mathrm{SO}_0(n,1)$), the complex hyperbolic $n$-space ($\mathrm{Isom}(X)=\mathrm{SU}(n,1)$), the quaternionic hyperbolic $n$-space ($\mathrm{Isom}(X)=\mathrm{Sp}(n,1)$) and the hyperbolic plane over the Cayley numbers ($\mathrm{Isom}(X)=\mathrm{F}_{4(-20)}$). For details on the classification of real Lie groups, see \cite{knapp_lie_1996}.
Another case of interest may be when $X=\HH^{\infty}$, the infinite dimensional hyperbolic space.\\ 

It is not hard to prove that primitive-stable representations form an open subspace of $\chi(\F_2,\mathrm{Isom}(X))$. We give a proof of this property in section \ref{openness}. The action of the outer automorphism group $\mathrm{Out}(\F_2)$ is properly discontinuous on  the set of Bowditch representation $\mathcal{BQ}(\F_2,X)$ so we obtain the following corollary :

\begin{Corollary}
$\mathcal{BQ}(\F_2,X)$ is an open domain of discontinuity for the action of $\mathrm{Out}(\F_2)$ on $\chi(\F_2,\mathrm{Isom}(X))$. 
\end{Corollary}

\begin{Remark} We can draw an analogy between Theorem \ref{BQ=PS} and a result of Delzant, Guichard, Labourie and Mozes (see \cite{delzant_displacing_2011}, see also a survey of Canary \cite{canary_dynamics_2015}) which says that for $\Gamma$ a hyperbolic group and $X$ a metric space, a representation $\rho : \Gamma \to \mathrm{Isom}(X)$ is displacing if and only if its orbit map is a quasi-isometric embedding. By displacing, the authors mean that the displacement length $l_\rho(\gamma)$ grows linearly with the cyclically reduced word length $\Vert \gamma \Vert$. The Bowditch condition (Definition \ref{def-bowditch}) can be seen as a restriction of the displacing condition on primitive elements, hence we could talk about \emph{primitive-displacing} representations. Primitive-stability is a restriction of the condition of quasi-isometric embedding of the orbit map on primitive leaves. Therefore Theorem \ref{BQ=PS} can be reinterpreted by saying that a representation of $\F_2$ is primitive-displacing if and only if its orbit map is a quasi-isometric embedding on primitive leaves. 
\end{Remark}

\begin{Remark} In a different but related direction, Maloni, Palesi and Tan (see \cite{maloni_character_2015}, see also Palesi \cite{palesi_dynamique_2014}) studied the $\mathrm{SL}(2,\C)$-character variety of the four-punctured sphere. In particular, they defined a Bowditch set $\mathcal{BQ}(\pi_1(S_{0,4}),\HH^3)$ and proved that it is an open domain of discontinuity for the action of the mapping class group. In addition, they provide a characterization similar to the Bowditch conditions \ref{trace-commutateur}., \ref{all-hyperbolic}. and \ref{set-finite}. Maloni, Palesi and Yang also studied the $\mathrm{PGL}(2,\R)$ and $\mathrm{SL}(2,\C)$ character variety of the three-holed projective plane in \cite{maloni_type-preserving_2021} and \cite{maloni_character_2020}.
\end{Remark}

\subsection{Strategy of the proof and outline of the paper}
At first, let's remark that the inclusion $\mathcal{PS}(\F_2,X) \subset \mathcal{BQ}(\F_2,X)$ is not hard to show. For completeness, a proof of this fact is given in Lemma \ref{inclusion-easy} of Section \ref{first-inclusion}. The difficulty of the theorem is all contained in the other inclusion. \\
We want to show, starting from a Bowditch representation $\rho$ from $\F_2$ to $\mathrm{Isom}(X)$, that it is primitive-stable, meaning that all the geodesics $L_\gamma$ (for $\gamma$ primitive) are mapped by the orbit map to uniform quasi-geodesics in $X$. It is almost immediate to see that under the Bowditch hypothesis, the geodesic $L_\gamma$ are mapped to quasi-geodesics, but the constants of quasi-geodesicity might depend on $\gamma$. The main difficulty, and first step of the proof, will be to show that these quasi-geodesics $\tau_\rho(L_\gamma)$ stay in a uniform neighborhood of the axis of $\rho(\gamma)$. Namely, this means that our family of quasi-geodesics $\tau_\rho(L_\gamma)$  satisfies a Morse lemma.  After this major step in the proof, done in section \ref{first-step-proof}, there will only be a little work left in order to show the primitive-stability of the representation $\rho$, and this will be done in section \ref{second-step-proof}. \\

Section \ref{structure-F2} is intended to establish some results about primitivity in $\F_2$. It will be dedicated to studying the structure of primitive words, and more specifically to understand the redundancy of primitive subwords within a primitive word in $\F_2$. One of the main ideas that will by highlighted and exploited is that if a primitive subword is found somewhere in a primitive word of $\F_2$, it can be found everywhere. More precisely, in definition \ref{def-wi} and proposition \ref{wi}, we explain how to decompose $\gamma$ (or maybe a cyclic permutation of $\gamma$) as a concatenation of primitive subwords for different scales, which will correspond to the successive steps in the continued fraction expansion of the slope of $\gamma$ (for the reader unfamiliar with these notions, they are recalled in the beginning of section \ref{structure-F2}). Then, in lemma \ref{magic-len}, we prove that for every primitive word $\gamma \in \F_2$, there exist some specific lengths, such that each subword of $\gamma$ of one of these lengths is "nearly" primitive, in the sense that it is so up to changing its last letter.  \\

Let us now explain what are the key ideas of the main step of the proof, which states (proposition \ref{uniform-neighborhood}) that the quasi-geodesics $\tau_\rho(L_\gamma)$ stay "close" to the axis of $\rho(\gamma)$. We will proceed by contradiction and suppose that we can find a primitive element $\gamma$ such that the associated quasi-geodesic $\tau_\rho(L_\gamma)$ does not stay close to the axis of $\rho(\gamma)$. Then we can find what we will call an \emph{excursion}, that is a path extracted from the quasi-geodesic that remains "far away" from the axis of $\rho(\gamma)$ (subsection \ref{excursion-real-map} and \ref{excursion-orbit-map}). We will next define the notion of a \emph{quasi-loop} (subsection \ref{quasi-loop}), which will be an element $u$ of $\F_2$ such that $\rho(u)$ does not displace the basepoint much, and prove in lemma \ref{excursion->QB} that every "big" excursion corresponds to a quasi-loop. This enables us to find a quasi-loop in the element $\gamma$. The goal will be next to find as many disjoints quasi-loops as possible in $\gamma$ and to do so, we will use the results of section \ref{structure-F2}. Indeed, our quasi-loop is contained in a subword of $\gamma$ whose length is one of the  specific lengths defined in lemma \ref{magic-len}, and thus this subword is primitive. But this primitive subword can be found  everywhere in $\gamma$, therefore with this process we will find our quasi-loop many times in $\gamma$. This will ensure that some proportion of~$\gamma$ does not displace the basepoint much (lemma \ref{decoupe-restes-c}). Finally, using a recursive argument, we will show that we can find an arbitrary big proportion of the word $\gamma$ that does not displace the basepoint much (lemma \ref{trouve-gamma}), which will be in contradiction with the Bowditch hypothesis. \\ 

Finally, in appendix \ref{length-path}, which is independent of the other sections, we develop on some classical properties in $\delta$-hyperbolic spaces to prove proposition \ref{long_ext_banane}, which establish a lower bound on the length of a path in a $\delta$-hyperbolic space which stays "far away" from a geodesic. This proposition will be needed in section \ref{first-step-proof}. 

\subsection*{Acknowledgements}
This work was done during my PhD at the Université de Lille and Université de Strasbourg. I am deeply grateful to my thesis advisor François Guéritaud for introducing me to this subject, and for his constant help and guidance during this work. I would also like to thank my thesis referee Louis Funar and Vincent Guirardel for useful comments and remarks. \\
The author is funded by the European Union (ERC, GENERATE, 101124349). Views and opinions expressed are however those of the author only and do not necessarily reflect those of the European Union or the European Research Council Executive Agency. Neither the European Union nor the granting authority can be held responsible for them.

\section{Structure of primitive elements in the free group of rank two}
\label{structure-F2}

\subsection{Constructing primitive elements} ~\\

In this section, we gather some results about primitive elements in $\F_2$.  Nielsen studied primitivity and automorphisms of free groups (in \cite{nielsen_uber_1918}, \cite{nielsen_isomorphismengruppe_1924}). The reader may also refer to \cite{series_geometry_1985} or \cite{gilman_enumerating_2011}. First recall that a \emph{primitive element} in $\F_2$ is an element which is part of some basis of $\F_2$. Fix once and for all $\{a,b\}$ a free generating set of $\F_2$ (hence $\F_2=\la a,b \ra$). Then, obviously, $a$ and $b$ are primitive elements and so are for example $a^{-1}$, $b^{-1}$, $ab$, $ab^{-1}$, $a^{-1}b^{1}$ and $a^nb$, for all $n \in \Z$.  Also note that primitivity is invariant by conjugacy. We denote by $\mathcal{P}(\F_2)$ the set of primitive elements of $\F_2$. We will also denote by $\mathcal{P}(\Z^2)$ the set of primitive elements of $\Z^2$, that is, again, the set of elements of $\Z^2$ which are part of a basis of $\Z^2$ (or equivalently, the set of elements $(p,q) \in \Z^2$ such that $p$ and $q$ are relatively prime numbers). 
Consider the abelianisation map : 

\begin{equation*}
    \mathrm{Ab} : \F_2 \longrightarrow \Z^2
\end{equation*}
It is a surjective morphism which sends any basis of $\F_2$ to a basis of $\Z^2$ (hence primitive elements of $\F_2$ to primitive elements of $\Z_2$). Moreover, since $\Z^2$ is abelian, the values of $\mathrm{Ab}$ are constant on conjugacy classes, thus we can consider the following map :
\begin{equation*}
    \tilde{\mathrm{Ab}} : \mathcal{P}(\F_2)/\sim \, \, \, \longrightarrow \mathcal{P}(\Z^2)/\pm
\end{equation*} where the quotient on the left hand side is taken up to conjugacy and inversion.
\begin{Proposition}
The map $\tilde{\mathrm{Ab}}$ is a bijection. Equivalently, the map 
\begin{align*}
    \mathrm{Slope} : \, \, \mathcal{P}(\F_2)/\sim & \longrightarrow   \Q \cup \infty \\
    [\gamma] & \longmapsto \frac{p}{q}, \text{ with } (p,q)=\mathrm{Ab}(\gamma)
\end{align*} where the quotient is taken up to conjugacy and inversion,
is a bijection. 
\end{Proposition}

Thus we have identified primitive elements (up to conjugacy and inversion) with rational numbers. Every rational number has a \emph{continued fraction expansion}, meaning that it can be written in the following way : 
\begin{equation*}
    \frac{p}{q} = n_1 + \cfrac{1}{n_2+\cfrac{1}{\ddots+\cfrac{1}{n_r}}}
\end{equation*} with $n_1 \in \Z, n_i \in \N^* \text{ for } i\geq 2$ and $n_r \geq 2$. Denote this expansion by $[n_1,n_2,\hdots,n_r]$. The continued fraction expansion of the slope will play a central role when studying the general structure of primitive elements. We will now give the general structure of a primitive element in $\F_2$. \\

Consider $w \in \F_2$. If $w$ is primitive, then, $w$ is either a word on $\{a,b\}$, on $\{a^{-1},b^{-1}\}$, on $\{a,b^{-1}\}$, or on $\{a^{-1},b\}$. In the first two cases, the slope of $w$ is positive and in the two last ones, negative. Thus up to inversion, $w$ can be written as a word on $\{a,b\}$ (positive slope), or $\{a,b^{-1}\}$ (negative slope). For simplicity, in the following we will only consider positive slope (for negative slope, just change $b$ to $b^{-1}$). We say that a letter $s$ is \emph{isolated} in a (cyclic) word $w$, if between two appearances of $s$ there is at least another letter. We say that a word $w$ in $\{a,b\}$, seen as a cyclic word (not necessarily primitive) is \emph{almost constant} if the two following conditions are satisfied : 
\begin{itemize}
    \item Either the letter $a$ or $b$ is isolated in $w$
    \item After possibly exchanging $a$ and $b$, suppose that $b$ is isolated in $w$. Then the powers of $a$ that arise in $w$ can only be two consecutive integers. 
\end{itemize}
In this case we say that the smallest integer that arises in $w$ as a power of $a$ is the \emph{value} of $w$. \\
In other words, a word $w$ in $a$ and $b$ is \emph{almost constant of value} $n \in \N$ if and only if there exists $s\in \N$ such that, after possibly exchanging $a$ and $b$ and up to conjugacy and inversion, $w$ is of the form~:
\begin{equation*}
    a^{n_1}ba^{n_2}b\hdots a^{n_s}b, \, \, \text{ with } n_i \in \{n,n+1\}, \, \, \forall 1 \leq i \leq s
\end{equation*}
If $w$ is almost constant, we can consider its \emph{derived word} by replacing the blocks $a^{n}b$ by $b$ and the blocks $a^{n+1}b$ by $ab$. The derived word is still a word on $a$ and $b$ and thus can be itself almost constant or not. \\
We say that a word is \emph{characteristic} if it can be derived arbitrarily many times, until a single letter is obtained. The \emph{values} of a characteristic word is the sequence of values of the almost constant derived words obtained at each step.

\begin{Proposition}
 Let $\gamma$ be an element of $\F_2$. 
 Then $\gamma$ is primitive if and only if it is characteristic. 
 Moreover, in this case, the values of the characteristic word $\gamma$ are $n_1,n_2,\hdots,n_r$, where $[n_1,n_2,\hdots,n_r]$ is the continued fraction expansion of the slope of $\gamma$.
\end{Proposition}

For proof, see \cite{series_geometry_1985}. Using this fact, we give as explicit construction of the (conjugacy class of) primitive elements, starting from their slope.

\begin{Definition} \label{def-wi} 
Let $\gamma$ be a primitive element in $\F_2$. Consider $[n_1(\gamma),n_2(\gamma),\hdots,n_{r(\gamma)}(\gamma)]$ the continued fraction expansion of the slope of $\gamma$ and assume that $n_i(\gamma)\geq 0$ (that is $\slope(\gamma)\geq 0$). We define recursively, for $0 \leq i \leq r(\gamma)$ the following elements in $\F_2$ : 
\begin{align*}
    w_0(\gamma) & = a & w'_0(\gamma) & = ab \\
    w_i(\gamma) & = w_{i-1}(\gamma)^{n_i(\gamma)-1}w'_{i-1}(\gamma) & w'_i(\gamma) & =w_{i-1}(\gamma)^{n_i(\gamma)}w'_{i-1}(\gamma) 
\end{align*}
Denote by $l_i(\gamma)$ and $l'_i(\gamma)$ the word lengths of $w_i(\gamma)$ and $w'_i(\gamma)$ respectively.
\end{Definition}

As defined, the $w_i(\gamma)$ are the building blocks of $\gamma$. In order to reduce the amount of notation, and when there will be no ambiguity on $\gamma$, we will omit the dependence on $\gamma$ in the notation and write $r,n_i,w_i, w'_i, l_i,l'_i,$. We can check the following : 

\begin{Proposition} \label{wi}
The elements $w_i$ defined previously satisfy :
\begin{enumerate}
    \item For all $0 \leq i \leq r$, $w_i$ and $w'_i$ are primitive and, for $i \geq 1$ their continued fraction expansions are respectively $[n_1,\hdots,n_i]$ and $[n_1,\hdots,n_i+1]$
    \item $w_r=\gamma$ (up to conjugacy and inversion). In particular, $\gamma$ (or its inverse) has a conjugate which is a positive word in $\{w_i,w'_i\}$.
    \item For all $0 \leq i \leq r$, $\{w_i,w'_i \}$ is a free basis of $\F_2$.
\end{enumerate}
\end{Proposition}
\begin{proof}
\begin{enumerate}
    \item For $i=0$, $w_0=a$ and $w'_0=ab$ are both trivially primitive.
    \begin{itemize}
        \item For $i=1$, $w_1=a^{n_1}b$ and $w'_0=a^{n_1+1}b$ are again both primitive and their continued fraction expansions are respectively $[n_1]$ and $[n_1+1]$.
        \item Suppose that both $w_{i-1}$ and $w'_{i-1}$ are primitive and that their continued fraction expansion are respectively $[n_1,\hdots,n_{i-1}]$ and $[n_1,\hdots,n_{i-1}+1]$. Then $w_{i-1}$ and $w'_{i-1}$ can be derived $i-1$ times, to obtain the elements $a$ and $ab$ (or $b$ and $ba$, depending on the parity of $i$). Thus, since $w_i$ and $w'_i$ are positive words on $\{w_{i-1},w'_{i-1} \}$, they can also be derived $i-1$ times, and the $(i-1)-st$-derived elements we obtain are $a^{n_i}b$ and $a^{n_i+1}b$ (or $b^{n_i}a$ and $b^{n_i+1}a$). Those last ones can be derived one more time to obtain $b$ and $ba$ (or $a$ and $ab$) and we have proved that $w_i$ and $w'_i$ are primitive with continued fraction expansion $[n_1,\hdots,n_i]$ and $[n_1,\hdots,n_i+1]$.
    \end{itemize}
    \item It follows directly from the previous point knowing that $w_r$ and $\gamma$ are both primitive with the same slope.
    \item This is an induction on $i$ using the basic fact that if $\{a,b\}$ is a basis of $\F_2$, then so are $\{a,ab\}$ and $\{a,ba\}$.
    \begin{itemize}
        \item The previous argument immediately justifies that $\{w_0,w'_0\}$ is a basis of $\F_2$. 
        \item Suppose that $\{w_{i-1},w'_{i-1}\}$ is a basis of $\F_2$, then, by the same argument as before, so is $\{w_{i-1},w_{i-1}w'_{i-1}\}$ and also, by induction $\{w_{i-1},w_{i-1}^{n_i-1}w'_{i-1}\}$. Now, we deduce the same way that $\{w_{i-1}^{n_i-1}w'_{i-1},w_{i-1}^{n_i}w'_{i-1}\}$ is a basis of $\F_2$ and thus that $\{w_i,w'_i\}$ is a basis of $\F_2$. 
    \end{itemize}
\end{enumerate}
\end{proof}

\begin{Remark} \label{l_i ineq}
Using the recursive definitions of $w_i$ and $w'_i$, we draw the following equalities : 
\begin{align}
    \text{ For all } \quad 2 \leq i \leq r \text{ and for } i=1 \text{ if  } n_1\geq 1, \qquad  &  l_i=(n_i-1)l_{i-1}+l'_{i-1} \, \, \, \text{   and } \, \, \, l'_i=n_il_{i-1}+l'_{i-1} \label{eq1} \\
    \text{ For } i=0 \text{ and for } i=1 \text{ with } n_1=0, \qquad & l_i=1 \text{ and } l'_i=2 \label{eq3} \\
    \text{ For all } \quad 1 \leq i \leq r, \qquad  & l'_i=l_i+l_{i-1} \label{eq2}
\end{align}
We deduce the following inequalities : 

\begin{align}
    \text{ For all } \quad 0 \leq i \leq r, \qquad & l_i<l'_i & \text{ using \eqref{eq2} and $l_{i-1}>0$ } \\
    \text{ For all } \quad 1 \leq i \leq r \text{ and for } i=0 \text{ if } n_1 \geq 1, \qquad & l_i < l_{i+1} & \text{ using \eqref{eq1}, $l'_{i-1} > l_{i-1}$ and $n_i \geq 1$ } \\
    \text{ For all } \quad 2 \leq i \leq r \text{ and for } i=1 \text{ if } n_1 \geq 1, \qquad & l'_i < 2l_i & \text{ using \eqref{eq2} and $l_{i-1} < l_i$ } \label{l'i<2li} \\
   \text{ For } i=0 \text{ and for } i=1 \text{ with } n_1=0, \qquad & l'_i=2l_i & \text{ using \eqref{eq3}}  \\ 
 \text{For all } 0 \leq i \leq r, \qquad & i \leq l_i & \text{by induction} \label{i-leq-li}
\end{align}
 Then $n_il_{i-1} < l_i \leq  (n_i+1)l_{i-1}$ (using \eqref{eq1}, \eqref{eq3} and $l_{i-1} < l'_{i-1} \leq 2l_{i-1}$), thus :
\begin{equation}
\text{ For all } 1 \leq i \leq r, \qquad
    \boxed{n_i<\frac{l_i}{l_{i-1}} \leq n_i+1}.
    \label{ineg-li}
\end{equation}

We also have $l_i-l_{i-1} \geq l'_{i-1}-l_{i-1} = l_{i-2}$, but $l_{i-1} \leq (n_{i-1}+1)l_{i-2}$, so $l_i-l_{i-1} \geq \frac{1}{n_{i-1}+1}l_{i-1}$, and then we deduce :
\begin{align}
 \text{For all } 2 \leq i \leq r, \qquad   \boxed{1+\frac{1}{n_{i-1}+1} \leq  \frac{l_i}{l_{i-1}}} \label{min-li}
\end{align}
\end{Remark}

\subsection{Some usefull lemmas}

Recall that the notation $|u|$ stands for the word length of the element $u \in \F_2$. For $u$ an element in $\F_2$ and $k$ an integer smaller that $|u|$, we will denote by $\overset{\curvearrowleft k}{u}$ its $k$-th cyclic permutation, that is, if $u=s_1 \hdots s_{|u|}$, the element $\overset{\curvearrowleft k}{u}=s_{k+1} \hdots s_{|u|}s_1 \hdots s_k$. We will also denote by $\p[k]{u}$ the prefix of length $k$ of $u$ and $\s[k]{u}$ the suffix of length $k$ of $u$. We have $u=\p[k]{u}\s[|u|-k]{u}$ and $\overset{\curvearrowleft k}{u}=\s[|u|-k]{u}\p[k]{u}$.\\
Sometimes for the sake of simplicity we will no longer specify the integer $k$ and write  $\p{u},\s{u},\cycl{u}$. Finally, we will write $\sm{u}$ to refer to a subword of $u$.

\begin{Lemma}\label{perm-cycl} 

Let $\gamma$ be a primitive element in $\mathbb{F}_2$ and $1 \leq i \leq r(\gamma)$. In particular, up to conjugacy, the element $\gamma$ can be written on the alphabet $\{ w_i(\gamma),w'_i(\gamma) \}$ (see Proposition \ref{wi}).

Now take $\cycl{w_i(\gamma)}$ any cyclic permutation of $w_i(\gamma)$. Then there exist $\cycl{w'_i(\gamma)}$ a cyclic permutation of $w'_i(\gamma)$ and $\cycl{\gamma}$ a cyclic permutation of $\gamma$ such that $\cycl{\gamma}$ can be written on the alphabet $\{ \cycl{w_i(\gamma)},\cycl{w'_i(\gamma)} \}$. Moreover, the element $\cycl{w_i(\gamma)}$ is either a prefix or a suffix of  $\cycl{w'_i(\gamma)}$. \\
In this case we say that the cyclic permutation $\cycl{w'_i(\gamma)}$ is \emph{adapted} to $\cycl{w_i(\gamma)}$.
\end{Lemma}

\begin{proof}
The case $i=0$ is trivial because $w_0=a$ and then $w_0$ has no non-trivial cyclic permutation. \\
Let $W$ be an element of $\mathbb{F}_2$ that can be written on the alphabet $\{ w_i,w'_i \}$, with $i \geq 1$. Then we can write $W=u_1 \hdots u_r$, with $u_j \in \{w_i,w'_i\}$ for $1 \leq j \leq r$. Consider the $k$-th cyclic permutation of $w_i$ : $\cycl[k]{w_i}=\s[l_i-k]{w_i}\p[k]{w_i}$.\\

We make the proof by distinguishing two cases :
\begin{itemize}
    \item \textbf{Case 1 :} If $k \leq (n_i-1)l_{i-1}$. \\ 
    Recall that, because $i>0$, we have the following recursive formulae :
    \begin{align*}
    w_i & =w_{i-1}^{n_i-1}w'_{i-1} \\
    w'_i & =w_{i-1}^{n_i}w'_{i-1}=w_{i-1}^{n_i-1}w_{i-1}w'_{i-1}
    \end{align*}
    Then in that case, we can say that $\p[k]{w_i}=\p[k]{w'_i}$ and so $ \forall 1 \leq j,j' \leq r, \p[k]{u_j}=\p[k]{u_{j'}}$. \\
    Thus : 
    \begin{align*}
        W & =u_1 u_2 \hdots u_r \\
          & = \p[k]{u_1}\s[|u_1|-k]{u_1} \p[k]{u_2}\s[|u_2|-k]{u_2} \hdots \p[k]{u_r}\s[|u_r|-k]{u_r} \\
        \cycl[k]{W}  & = \s[|u_1|-k]{u_1} \p[k]{u_2} \s[|u_2|-k]{u_2} \hdots \p[k]{u_r}\s[|u_r|-k]{u_r} \p[k]{u_1} \\
          & = \s[|u_1|-k]{u_1} \p[k]{u_1} \s[|u_2|-k]{u_2} \hdots \p[k]{u_{r-1}}\s[|u_r|-k]{u_r} \p[k]{u_r} \\
          & = \cycl[k]{u_1} \hdots \cycl[k]{u_r}
    \end{align*}
    We have proved that in that case, $\cycl[k]{W}$ can be written on the alphabet $\{ \cycl[k]{w_i},\cycl[k]{w'_i}\}$. \\
    
 Now let us show that $\cycl[k]{w_i}$ is a suffix of $\cycl[k]{w'_i}$.  We have :
    \begin{align*}
        \cycl[k]{w_i} & = \s[l_i-k]{w_i}\p[k]{w_i}\\
        \cycl[k]{w'_i} & = \s[l'_i-k]{w'_i}\p[k]{w'_i}=\s[l'_i-k]{w'_i}\p[k]{w_i}
    \end{align*}
    Since $w_i$ is a suffix of $w'_i$ and $l'_i-k \geq l_i-k$, it implies that $\s[l_i-k]{w_i}$ is a suffix of $\s[l'_i-k]{w'_i}$ and then that $\cycl[k]{w_i}$ is a suffix of $\cycl[k]{w'_i}$. \\
    
    \item \textbf{Case 2 :} If $k > (n_i-1)l_{i-1}$. \\
    Let $k'=k-(n_i-1)l_{i-1}$. \\
    For all $1 \leq j \leq r$, there exists $\varepsilon_j \in \{n_i-1,n_i \}$ such that $u_j = w_{i-1}^{\varepsilon_j}w'_{i-1}$. \\
    Then $u_j=w_{i-1}^{\varepsilon_j}\p[k']{w'_{i-1}}\s[l'_{i-1}-k']{w'_{i-1}}$ and so :
    \begin{align*}
        W & = u_1 u_2 \hdots u_r \\
          & = w_{i-1}^{\varepsilon_1}\p[k']{w'_{i-1}}\s[l'_{i-1}-k']{w'_{i-1}} \hspace{0,2cm} w_{i-1}^{\varepsilon_2}\p[k']{w'_{i-1}}\s[l'_{i-1}-k']{w'_{i-1}} \hspace{0,2cm} \hdots \hspace{0,2cm}  w_{i-1}^{\varepsilon_r}\p[k']{w'_{i-1}}\s[l'_{i-1}-k']{w'_{i-1}}
    \end{align*}
    Now let $k_j=\varepsilon_j l_{i-1}+k'$. \\ Then
    $
     k_j = \left\{
    \begin{array}{ll}
        k & \text{ if } \varepsilon_j=n_i-1 \\
        k+l_{i-1} & \text{ if } \varepsilon = n_i
    \end{array}
        \right. 
    = \left\{
    \begin{array}{ll}
        k  & \text{ if } u_j=w_i \\
        k+l_{i-1}  & \text{ if } u_j=w'_i
    \end{array}
\right.
    $ \\
    Therefore $\s[l'_{i-1}-k']{w'_{i-1}}w_{i-1}^{\varepsilon_j}\p[k']{w'_{i-1}}=\cycl[k_j]{u_j}$. Finally : 
    
    \begin{equation*}
        \cycl[k_1]{W} = \cycl[k_2]{u_2} \hdots \cycl[k_r]{u_r} \cycl[k_1]{u_1}
    \end{equation*}
    
    We have proved that in that case, $\cycl[k_1]{W}$ can be written on the alphabet $\{ \cycl[k]{w_i}, \cycl[k+l_{i-1}]{w'_{i-1}} \}$. \\
    
    Now let us show that $\cycl[k]{w_i}$ is a prefix or a suffix of  $\cycl[k+l_{i-1}]{w'_i}$. We have : 
     \begin{align*}
        \cycl[k]{w_i} & = \s[l'_{i-1}-k']{w'_{i-1}}w_{i-1}^{n_i-1}\p[k']{w'_{i-1}}\\
        \cycl[k+l_{i-1}]{w'_i} & = \s[l'_{i-1}-k']{w'_{i-1}}w_{i-1}^{n_i}\p[k']{w'_{i-1}}
    \end{align*}
    We first handle the case $i=1$. In that case, $w'_0=ab, l'_0=2$. There are three possibilities~: either $k'=0$, and then $\cycl[k]{w_i}$ is a prefix of $\cycl[k+l_{i-1}]{w'_i}$, or $k'=2$, and then $\cycl[k]{w_i}$ is a suffix of~$\cycl[k+l_{i-1}]{w'_i}$, or $k'=1$, and then $\cycl[k]{w_i}=ba^{n_i-1}a$ and $\cycl[k+l_{i-1}]{w'_i} = ba^{n_i}a$ so $\cycl[k]{w_i}$ is a prefix of $\cycl[k+l_{i-1}]{w'_i}$. \\
    
    For the case $i \geq 2$, we can use the recursive formulae : 
    \begin{align*}
    w_{i-1} & =w_{i-2}^{n_{i-1}-1}w'_{i-2} \\
    w'_{i-1} & =w_{i-2}^{n_{i-1}}w'_{i-2}
    \end{align*}
    
    We distinguish two cases :
    \begin{itemize}
        \item \underline{If $k' \leq (n_{i-1}-1)l_{i-2}$} : then $\p[k']{w'_{i-1}}$ is a prefix of $w_{i-1}$ so $\cycl[k]{w_i}$ is a prefix of  $\cycl[k+l_{i-1}]{w'_i}$.
        \item \underline{If $k' > (n_{i-1}-1)l_{i-2}$ } : then $l'_{i-1}-k'<l'_{i-1}-(n_{i-1}-1)l_{i-2}=l_{i-2}+l'_{i-2}$. So $\s[l'_{i-1}-k']{w'_{i-1}}$ is a suffix of $w_{i-2}w'_{i-2}$. 
        \begin{itemize}
            \item If $n_{i-1}>1$, $w_{i-2}w'_{i-2}$ is a suffix of $w_{i-1}$ so $\s[l'_{i-1}-k']{w'_{i-1}}$ is a suffix of $w_{i-1}$ and so  $\cycl[k]{w_i}$ is a suffix of $\cycl[k+l_{i-1}]{w'_i}$.
            \item If $n_{i-1}=1$, we have $w_{i-1}=w'_{i-2}$ and $w'_{i-1}=w_{i-2}w'_{i-2}$. If $i>2$, we have $w'_{i-1}=w_{i-2}w_{i-3}w_{i-2}$ and so a suffix of $w_{i-2}w'_{i-2}$ is also a suffix of $w'_{i-1}w_{i-1}$, so $\s[l'_{i-1}-k']{w'_{i-1}}$ is a suffix of $w'_{i-1}w_{i-1}$ and so $\cycl[k]{w_i}$ is a suffix of  $\cycl[k+l_{i-1}]{w'_i}$. If $i=2$, we have $w_{i-1}=ab$ and $w'_{i-1}=aab$, so $w'_{i-1}w_{i-1}=aabab$. Moreover $w_{i-2}w'_{i-2}=aab$ so   
            \begin{align*}
            \cycl[k]{w_i} & = \s[l'_{i-1}-k']{aab}(ab)^{n_i-1}\p[k']{aab}\\
            \cycl[k+l_{i-1}]{w'_i} & = \s[l'_{i-1}-k']{aab}(ab)^{n_i}\p[k']{aab}
            \end{align*}
            And we notice that again in that case, $\cycl[k]{w_i}$ is either a prefix or a suffix of $\cycl[k+l_{i-1}]{w'_i}$.
        \end{itemize}
        
        \end{itemize}
\end{itemize}
\end{proof}

The following Proposition is the most important of this section. It says that subwords of some specific lengths ($l_i(\gamma)$) of primitive elements are themselves primitive, after a possible change of letter. 

\begin{Proposition} \label{magic-len}
Let $\gamma$ be a primitive element of $\mathbb{F}_2$ and $u$ any subword of $\gamma$ (or of a cyclic permutation of $\gamma$) of length $l_i(\gamma)$, for $0 \leq i \leq r(\gamma)$. Then, after possibly changing its last letter, $u$ is in fact a cyclic permutation of $w_i(\gamma)$.
\end{Proposition}

\begin{proof}
For the sake of simplicity, we omit the dependence on $\gamma$ in the notations in the proof. \\

First of all, we deal separately with the cases $i=0$ and $i=1$. 
\begin{itemize}
    \item If $i=0$, then $u=a$ or $u=b$, and thus it is trivial.
    \item If $i=1$, then $u=a^{n_1+1}$ or $u=a^kba^{n_1-k}$, with $0 \leq k \leq n_1$. In the first case, after changing the last letter $a$ into $b$, we obtain $u=a^{n_1}b=w_1$. In the second case, no change is needed.  \\
\end{itemize}

From now on, we suppose that $i\geq 2$. In that case, the recursive formulae $w_i=w_{i-1}^{n_i-1}w'_{i-1}$ and $w'_i=w_{i-1}^{n_i}w'_{i-1}$ apply. \\
The word $\gamma$ can be written on the alphabet $\{w_i,w'_i \}$, and $l'_i \geq l_i$ so the subword $u$ of $\gamma$ shall take one of the following forms : $\s{w_i}\p{w_i}, \s{w_i}\p{w'_i},\s{w'_i}\p{w_i}, \s{w'_i}\p{w'_i} \text{ or } \sm{w'_i}$. \\
Furthermore, since $w'_i=w_{i-1}w_i$ and $|\s{w'_i}| \leq l_i$, we deduce that $\s{w'_i}=\s{w_i}$. Thus, $u$ can actually be reduced to one of the following three forms :

\begin{enumerate}
    \item \label{gentil} $u=\s{w_i}\p{w_i}$
    \item \label{moyen} $u=\sm{w'_i}$
    \item \label{rec} $u=\s{w_i}\p{w'_i}$
\end{enumerate}
We deal with each case separately. 

\begin{enumerate}

    \item The case \ref{gentil} is actually immediate because $|\s{w_i}\p{w_i}|=l_i$, so $|\p{w_i}\s{w_i}|=l_i$, which requires $\p{w_i}\s{w_i}=w_i$ and then $u=\s{w_i}\p{w_i}$ is a cyclic permutation of $w_i$. \\
    \item Now let us deal with the case \ref{moyen}. Recall $w'_i=w_{i-1}^{n_i}w'_{i-1}=w_{i-1}w_{i-1}^{n_i-1}w'_{i-1}$ (because $n_i \geq 1$). Then a subword of length $l_i$ of $w'_i$ must be of the form : $u=\sm{w'_i}=\s{w_{i-1}}w_{i-1}^{n_i-1}\p{w'_{i-1}}$. Therefore, up to cyclic permutation we have $\cycl{u}=w_{i-1}^{n_i-1}\p{w'_{i-1}}\s{w_{i-1}}$, with $|\p{w'_{i-1}}\s{w_{i-1}}|=l'_{i-1}$. Moreover, for $i \geq 2$, we have $\s{w_{i-1}}=\s{w'_{i-1}}$ and then $\p{w'_{i-1}}\s{w_{i-1}}=\p{w'_{i-1}}\s{w'_{i-1}}=w'_{i-1}$. Thus $\cycl{u}=w_{i-1}^{n_i-1}\p{w'_{i-1}}\s{w_{i-1}}=w_{i-1}^{n_i-1}w'_{i-1}=w_i$. \\
    
    \item For the case \ref{rec}, we use a two-step induction on $i \geq 0$. More precisely, we will show by induction on $i\geq 0$, that if $u$ is a subword of $\gamma$ of the form $u=\s{w_i}\p{w'_i}$ with $|u|=l_i$, then, after possibly changing its last letter, that is after possibly changing the last letter of $\p{w'_i}$, we obtain $\phat{w'_i}\s{w_i}=w_i$ or $\phat{w'_i}\s{w_i}=w_{i-1}^{n_i}w_{i-2}$ (the latter case can only occur if $i \geq 2$), with $\phat{w'_i}$ the word obtained from $\p{w'_i}$ after the change of letter. Thus, we deduce that after possibly changing its last letter, $u$ is actually a cyclic permutation of $w_i$. \\
    In the following, the notation $\phat{w}$ stands for the word obtained from $\p{w}$ by changing its last letter. \\ 
    
    \textbf{Initial cases :} \\
    The trivial case $i=0$ has already been mentioned at the beginning of the proof. \\
    Suppose $i=1$. Then $\s{w_1}=\s{a^{n_1}b}$. If $|\s{w_1}| \geq 1$, then there exists an integer $k$ such that $\s{w_1}=a^kb$ and $\p{w'_1}=a^{n_1-k}$ and so $\p{w'_1}\s{w_1}=a^{n_1}b=w_1$. If $|\s{w_1}|= 0$, then $\p{w'_1}\s{w_1}=\p{w'_1}=a^{n_1+1}$. So after possibly changing the last letter of $\p{w'_1}$ into $b$, we have  $\phat{w'_1}\s{w_1}=a^{n_1}b=w_1$. \\
    
    \textbf{Induction step :} We now fix $i\geq 2$\\
    Let $k=|\p{w'_i}|$. Let us distinguish two cases :
    
    \begin{enumerate}
       
        \item If $k \leq (n_i-1)l_{i-1}(=l_i-l'_{i-1})$ : then $\p{w'_i}=\p{w_i}$,  so as before in the proof we have $\p{w'_i}\s{w_i}=\p{w_i}\s{w_i}=w_i$ because $|\p{w_i}\s{w_i}|=l_i$. Thus $\boxed{\p{w'_i}\s{w_i}=w_i}$ \\
        
        \item If $k \geq (n_i-1)l_{i-1}(=l_i-l'_{i-1})$, then $\p{w'_i}=w_{i-1}^{n_i-1}\p{w_{i-1}w'_{i-1}}$ on one hand and $|\s{w_i}|\leq l'_{i-1}$ so $\s{w_i}=\s{w'_{i-1}}$ on the other hand. Thus $\boxed{\p{w'_i}\s{w_i}=w_{i-1}^{n_i-1}\p{w_{i-1}w'_{i-1}}\s{w'_{i-1}}}$. \\
        
        We now need to understand the word $\p{w_{i-1}w'_{i-1}}\s{w'_{i-1}}$, with $|\p{w_{i-1}w'_{i-1}}\s{w'_{i-1}}|=l'_{i-1}$. \\
        
        \begin{enumerate}
            \item \underline{If $|\p{w_{i-1}w'_{i-1}}|\leq (n_{i-1}-1)l_{i-2}$}, then $\p{w_{i-1}w'_{i-1}}=\p{w_{i-2}^{n_{i-1}-1}}=\p{w'_{i-1}}$, so \\  $\p{w_{i-1}w'_{i-1}}\s{w'_{i-1}}=\p{w'_{i-1}}\s{w'_{i-1}}=w'_{i-1}$ because $|\p{w_{i-1}w'_{i-1}}\s{w'_{i-1}}|=l'_{i-1}$. In this case we obtain  $\boxed{\p{w_{i-1}w'_{i-1}}\s{w'_{i-1}}=w'_{i-1}}$.\\
            
            \item \underline{If $(n_{i-1}-1)l_{i-2} \leq |\p{w_{i-1}w'_{i-1}}|\leq n_{i-1}l_{i-2}$}, then $\p{w_{i-1}w'_{i-1}}=w_{i-2}^{n_{i-1}-1}\p{w'_{i-2}}$ on one hand, and $l'_{i-2} \leq |\s{w'_{i-1}}|\leq l'_{i-2}+l_{i-2} $ so $\s{w'_{i-1}}=\s{w_{i-2}}w'_{i-2}$ on the other hand. \\
            Thus $\p{w_{i-1}w'_{i-1}}\s{w'_{i-1}}=w_{i-2}^{n_{i-1}-1}\p{w'_{i-2}}\s{w_{i-2}}w'_{i-2}$, with $|\p{w'_{i-2}}\s{w_{i-2}}|=l_{i-2}$. Here we use our induction hypothesis : \\
            
            \begin{itemize}
                \item If, after possibly changing the last letter of $\p{w'_{i-2}}$, we have  $\phat{w'_{i-2}}\s{w_{i-2}}=w_{i-2}$, then, after possibly changing the last letter of $\p{w_{i-1}w'_{i-1}}$, $$\phat{w_{i-1}w'_{i-1}}\s{w'_{i-1}}=w_{i-2}^{n_{i-1}-1}\phat{w'_{i-2}}\s{w_{i-2}}w'_{i-2}=w_{i-2}^{n_{i-1}-1}w_{i-2}w'_{i-2}=w'_{i-1}.$$ Thus, as in the previous case, $\boxed{\phat{w_{i-1}w'_{i-1}}\s{w'_{i-1}}=w'_{i-1}}.$
                
                \item If, after possibly changing the last letter of $\p{w'_{i-2}}$, we have $\phat{w'_{i-2}}\s{w_{i-2}}=w_{i-3}^{n_{i-2}}w_{i-4}$ (recall that this case can only occur if $i\geq 4$, as stated at the beginning of the induction), then, after possibly changing the last letter of $\p{w_{i-1}w'_{i-1}}$,
                \begin{align*}
                \phat{w_{i-1}w'_{i-1}}\s{w'_{i-1}} & = w_{i-2}^{n_{i-1}-1}\phat{w'_{i-2}}\s{w_{i-2}}w'_{i-2}=w_{i-2}^{n_{i-1}-1}w_{i-3}^{n_{i-2}}w_{i-4}w'_{i-2}\\
                &=w_{i-2}^{n_{i-1}-1}w_{i-3}^{n_{i-2}}w_{i-4}w_{i-3}w_{i-2} =w_{i-2}^{n_{i-1}-1}w_{i-3}^{n_{i-2}}w'_{i-3}w_{i-2} \\
                & = w_{i-2}^{n_{i-1}-1}w'_{i-2}w_{i-2}=w_{i-1}w_{i-2}.
                \end{align*}
                Thus, in that case, we obtain  $\boxed{\phat{w_{i-1}w'_{i-1}}\s{w'_{i-1}}=w_{i-1}w_{i-2}}$. \\
                
            \end{itemize}
            
            \item \underline{If $n_{i-1}l_{i-2} \leq |\p{w_{i-1}w'_{i-1}}|\leq l'_{i-1}$}, then $\p{w_{i-1}w'_{i-1}}=w_{i-2}^{n_{i-1}-1}\p{w'_{i-2}w_{i-2}}$ on one hand and $|\s{w'_{i-1}}| \leq l'_{i-2} $ so $\s{w'_{i-1}}=\s{w'_{i-2}}$ on the other hand. \\
            Thus $\boxed{\p{w_{i-1}w'_{i-1}}\s{w'_{i-1}}=w_{i-2}^{n_{i-1}-1}\p{w'_{i-2}w_{i-2}}\s{w'_{i-2}}}$ with $|\p{w'_{i-2}w_{i-2}}\s{w'_{i-2}}|=l'_{i-2}+l_{i-2}$. We still need to understand $\p{w'_{i-2}w_{i-2}}\s{w'_{i-2}}$ with $|\p{w'_{i-2}w_{i-2}}\s{w'_{i-2}}|=l'_{i-2}+l_{i-2}$. \\ 
            
            \begin{itemize}
                \item \underline{If $|\p{w'_{i-2}w_{i-2}}| \geq l'_{i-2}$}, then $\p{w'_{i-2}w_{i-2}}=w'_{i-2}\p{w_{i-2}}$ on one hand and \\ $|\s{w'_{i-2}}| \leq l_{i-2} $ so $\s{w'_{i-2}}=\s{w_{i-2}}$. Now we can compute  $\p{w'_{i-2}w_{i-2}}\s{w'_{i-2}}=w'_{i-2}\p{w_{i-2}}\s{w_{i-2}}=w'_{i-2}w_{i-2}$ and the last equality stands because \\ $|\p{w_{i-2}}\s{w_{i-2}}|=l_{i-2}$. Thus $\boxed{\p{w'_{i-2}w_{i-2}}\s{w'_{i-2}}=w'_{i-2}w_{i-2}}.$
                \item \underline{If $|\p{w'_{i-2}w_{i-2}}| \leq l'_{i-2}$}, then we have $l_{i-2}\leq |\p{w'_{i-2}w_{i-2}}| \leq l'_{i-2}$ because \\  $|\s{w'_{i-2}}|\leq l'_{i-2}$. We also deduce that $l_{i-2} \leq |\s{w'_{i-2}}| \leq l'_{i-2} $. \\ 
                If $i=2$, then $l_{i-2}=l_0=1, l'_{i-2}=l'_0=2, w'_0w_0=aba$ and $w'_0=ab$. Then $\p{w'_0w_0}\s{w'_0}=aab$ or $\p{w'_0w_0}\s{w'_0}=abb$. So, after possibly changing the last letter of $\p{w'_0w_0}$, we have $\phat{w'_0w_0}\s{w'_0}=aab=w_0w'_0$. \\
                If $i>2$, we can write $w'_{i-2}=w_{i-3}^{n_{i-2}}w'_{i-3}$ and so $\p{w'_{i-2}w_{i-2}}=w_{i-3}^{n_{i-2}}\p{w'_{i-3}}$ and $\s{w'_{i-2}}=\s{w_{i-3}}w_{i-2}$. Thus $\p{w'_{i-2}w_{i-2}}\s{w'_{i-2}}=w_{i-3}^{n_{i-2}}\p{w'_{i-3}}\s{w_{i-3}}w_{i-2}$ with $|\p{w'_{i-3}}\s{w_{i-3}}|=l'_{i-3}$. But, as seen in the case \ref{moyen}, for $i>3, \p{w'_{i-3}}\s{w_{i-3}}=w'_{i-3}$. Then we compute $\p{w'_{i-2}w_{i-2}}\s{w'_{i-2}}=w_{i-3}^{n_{i-2}}\p{w'_{i-3}}\s{w_{i-3}}w_{i-2}=w_{i-3}^{n_{i-2}}w'_{i-3}w_{i-2}=w'_{i-2}w_{i-2}$. \\
                For $i=3$, recall that $w_{i-3}=w_0=a$ and $w'_{i-3}=w'_0=ab$. If $|\s{w_0}|=0$, then $\p{w'_{i-3}}\s{w_{i-3}}=\p{w'_0}\s{w_0}=\p{w'_0}=w'_0$. Thus $\p{w'_{i-2}w_{i-2}}\s{w'_{i-2}}=w_{i-3}^{n_{i-2}}\p{w'_{i-3}}\s{w_{i-3}}w_{i-2}=w_{i-3}^{n_{i-2}}w'_{i-3}w_{i-2}=w'_{i-2}w_{i-2}$. However, if $|\s{w_0}|=1$, then $\p{w'_{i-3}}\s{w_{i-3}}=\p{w'_0}\s{w_0}=\p{w'_0}a=a^2$, and then, after changing the last letter of $\p{w'_0}$ from $a$ to $b$, we obtain $\phat{w'_0}\s{w_0}=ba$. Thus, let us now compute $\phat{w'_{i-2}w_{i-2}}\s{w'_{i-2}}$ : \\
                $\phat{w'_{i-2}w_{i-2}}\s{w'_{i-2}}=w_{i-3}^{n_{i-2}}\phat{w'_{i-3}}\s{w_{i-3}}w_{i-2}=w_{i-3}^{n_{i-2}}baw_{i-2}=a^{n_1}baw_1=w_1w'_1$.\\
                \end{itemize}
                
            Thus, we conclude from those different cases that, after possibly changing the last letter of $\p{w'_{i-2}w_{i-2}}$, we have either $\phat{w'_{i-2}w_{i-2}}\s{w'_{i-2}}=w'_{i-2}w_{i-2}$, or \\ $\phat{w'_{i-2}w_{i-2}}\s{w'_{i-2}}=w_{i-2}w'_{i-2}$. We can now compute $\phat{w_{i-1}w'_{i-1}}\s{w'_{i-1}}$ : \\
            \begin{itemize}
                \item If $\phat{w'_{i-2}w_{i-2}}\s{w'_{i-2}}=w'_{i-2}w_{i-2}$, then 
            $$ \phat{w_{i-1}w'_{i-1}}\s{w'_{i-1}}=w_{i-2}^{n_{i-1}-1}\phat{w'_{i-2}w_{i-2}}\s{w'_{i-2}} = w_{i-2}^{n_{i-1}-1}w'_{i-2}w_{i-2} = w_{i-1}w_{i-2}.$$ 
                \item If $\phat{w'_{i-2}w_{i-2}}\s{w'_{i-2}}=w_{i-2}w'_{i-2}$, then 
            $$ \phat{w_{i-1}w'_{i-1}}\s{w'_{i-1}}=w_{i-2}^{n_{i-1}-1}\phat{w'_{i-2}w_{i-2}}\s{w'_{i-2}}=w_{i-2}^{n_{i-1}-1}w_{i-2}w'_{i-2}= w_{i-2}^{n_{i-1}}w'_{i-2}=w'_{i-1}.$$
            \end{itemize}
        \end{enumerate}
        
        Thus, we showed that after possibly changing the last letter of  $\p{w_{i-1}w'_{i-1}}$, we have $\phat{w_{i-1}w'_{i-1}}\s{w'_{i-1}}=w'_{i-1}$ or $\phat{w_{i-1}w'_{i-1}}\s{w'_{i-1}}=w_{i-1}w_{i-2}$. We are now ready to compute $\phat{w'_i}\s{w_i}$ :
        
        \begin{itemize}
            \item If $\phat{w_{i-1}w'_{i-1}}\s{w'_{i-1}}=w'_{i-1}$, then
            $$\phat{w'_i}\s{w_i}=w_{i-1}^{n_i-1}\phat{w_{i-1}w'_{i-1}}\s{w'_{i-1}}=w_{i-1}^{n_i-1}w'_{i-1}=w_i.$$
            \item If $\phat{w_{i-1}w'_{i-1}}\s{w'_{i-1}}=w_{i-1}w_{i-2}$, then 
            $$\phat{w'_i}\s{w_i}=w_{i-1}^{n_i-1}\phat{w_{i-1}w'_{i-1}}\s{w'_{i-1}}=w_{i-1}^{n_i-1}w_{i-1}w_{i-2}=w_{i-1}^{n_i}w_{i-2}.$$
        \end{itemize}
    \end{enumerate}
    Thus, we have proven what we announced for the case \ref{rec}. 
\end{enumerate}
Hence the lemma is proved. 
\end{proof}

We end this section with a short Lemma which "counts" the number of occurrences of $w_i(\gamma)$ in a subword $u$ of $\gamma$. 

\begin{Lemma} \label{bloc}
Let $\alpha>4$. Let $\gamma$ be a primitive element of $\mathbb{F}_2$ and $u$ any subword of $\gamma$ (or of a cyclic permutation of $\gamma$). Let $i \in \{1,\hdots,r(\gamma)\}$ and suppose that $|u| \geq \alpha l_i(\gamma)$. Let $\overset{\curvearrowleft}{w_i(\gamma)}$ be any cyclic permutation of $w_i(\gamma)$ and $\overset{\curvearrowleft}{w'_i(\gamma)}$ a cyclic permutation of $w'_i(\gamma)$ adapted to $w_i(\gamma)$ (see lemma \ref{perm-cycl}). Then, there is at least $\frac{\alpha -4}{2}$ occurrences of $w_i(\gamma)$ and $w'_i(\gamma)$ in $u$.
\end{Lemma}

\begin{proof} Recall that $l_i(\gamma)\leq l'_i(\gamma) \leq 2 l_i(\gamma)$. Then we have the inequality $|u| \geq \frac{\alpha}{2}l'_i(\gamma)$. \\
The element $u$ is a subword of $\gamma$ and, by Lemma \ref{perm-cycl}, $\gamma$ can be written on the alphabet $\{ \overset{\curvearrowleft}{w_i(\gamma)},\overset{\curvearrowleft}{w'_i(\gamma)} \}$ so $u$ can be written in the following way: $u=\p{u}u_1 \hdots u_r \s{u}$, with $u_k \in \{ \overset{\curvearrowleft}{w_i(\gamma)},\overset{\curvearrowleft}{w'_i(\gamma)} \}$ for $1 \leq k \leq r$ and $\p{u}$ and $\s{u}$ being respectively a prefix and a suffix of $u$ such that $|\p{u}|,|\s{u}| \leq |\overset{\curvearrowleft}{w'_i(\gamma)}|=l'_i(\gamma)$. 
Thus $|u|-|\p{u}|-|\s{u}| \geq \frac{\alpha }{2}l'_i(\gamma)-2l'_i(\gamma)=\frac{\alpha-4}{2}l'_i(\gamma)$. \\
Furthermore $|u|-|\p{u}|-|\s{u}| = \sum_{k=1}^r |u_k| \leq rl'_i(\gamma)$ because $|u_k| \leq l'_i(\gamma)$. \\
We deduce $r\geq \frac{\alpha -4}{2}$, hence Lemma \ref{bloc}.
\end{proof}

\section{An important example of a uniform quasi-geodesicity setting} 
\label{UQG-setting}

In this section, we are going to study a uniform quasi-geodesicity setting in the space $X$. The space $X$ is supposed to be a $\delta$-hyperbolic, geodesic and visibility space. When given any isometry $A$ in the metric space $X$, we have already defined in the introduction its \emph{displacement length}, and this notion has been used to define \emph{Bowditch representations} (see definition \ref{def-bowditch}). We can also consider its \emph{stable length}, which is defined by $l_S(A)= \underset{n \to \infty} \lim \frac{1}{n}d(A^n o,o)$. We can check that the stable length is well-defined and invariant under the choice of the basepoint $o$ in~
$X$. Moreover, it satisfies $l_S(A^n)=nl_S(A)$ (whereas this is not true in general for the displacement length). Suppose now that $A$ is a hyperbolic isometry of $X$, then by definition, the map from $\Z$ to $X$ that sends $n$ to $A^no$ is a quasi-isometry. Moreover, $A$ has two fixed points in the boundary $\partial X$, one attracting and the other repelling, denoted by $A^+$ and $A^-$ respectively.  \\

When we are given two isometries $A$ and $B$, we can consider the set $\mathcal{W}(A,B)$ of (finite) words on $A$ and $B$. For $G \in \mathcal{W}(A,B)$, we denote by $|G|$ its word length, that is the minimal number of letters ($A$ and $B$) needed to write $G$. We also consider $\mathcal{H}(A,B)=\{A,B\}^\Z$ the set of bi-infinite words on $A$ and $B$, that is, $H=(H_n)_{n\in\Z} \in \mathcal{H}(A,B)$ if and only if for all $n \in \Z$, $H_n \in \{A,B\}$. When we have a bi-infinite word $H=(H_n)_{n\in \Z} \in \mathcal{H}(A,B)$, we associate to it a bi-infinite sequence of finite words $G=(G_n)_{n \in \Z}$ in the following way : 
\begin{align*}
    G_n & =\left\{ \begin{array}{ll}
       H_0H_1\hdots H_{n-1}  &  \text{ for } n>0 \\
       I_d  &  \text{ for } n=0 \\ 
       H_{-1}^{-1}\hdots H_{n}^{-1} & \text{ for } n<0
    \end{array}
    \right.
\end{align*}
Hence $G_n \in \mathcal{W}(A,B)$ for all $n \geq 0$ and $G_n \in \mathcal{W}(A^{-1},B^{-1})$ for all $n \leq 0$. Moreover, the word length of $G_n$ is  $|G_n|=|n|$ for all $n \in \Z$ and the following recursive formula holds for all $n \in \Z$ : $G_{n+1}=G_nH_n$. We denote by $\mathcal{G}(A,B)$ the set of bi-infinite sequences of finite words associate to bi-infinite words in $\mathcal{H}(A,B)$. \\

In this section, we want to study a particular class of bi-infinite words $H=(H_n)_{n \in \Z}$ and their associate bi-infinite sequences of words $G=(G_n)_{n\in \Z}$. Let us fix an integer $N \geq 1$ and define $\mathcal{H}_N(A,B)$ to be the subset of $\mathcal{H}(A,B)$ consisting of the bi-infinite words $H=(H_n)_{n\in \Z}$ which satisfy the following condition : \\
If $n_1 < n_2$ are two integers in $\Z$ such that $H_{n_1}=H_{n_2}=B$ and 
for all $n_1 <n<n_2$, $H_n=A$, then $n_2-n_1-1 \geq N$. \\
Thus the bi-infinite words in $\mathcal{H}_N(A,B)$ are precisely those for which the appearances of $B$ are isolated and the powers of $A$ are always greater than $N$. We denote by $\mathcal{G}_N(A,B)$ the set of bi-infinite sequences $G=(G_n)_{n\in \Z} \in \mathcal{G}(A,B)$ associate to bi-infinite words in $\mathcal{H}_N(A,B)$. \\

Fix $o$ a basepoint in $X$. Starting from a bi-infinite word $H=(H_n)_{n\in \Z} \in \mathcal{H}(A,B)$ and its associate bi-infinite sequence $G=(G_n)_{n \in \mathbb{Z}} \in \mathcal{G}(A,B)$, we define the sequence of points in $X$ : $x_n=G_no, \forall n \in \Z$. The goal of this section is to study the uniform quasi-geodesicity of sequences of points defined by the elements of $\mathcal{G}_N(A,B)$, that is the existence of two reals $\lambda >0$ and $k \geq 0$ such that for all $n,m \in \mathbb{Z}$, we have : $\frac{1}{\lambda}|n-m|-k \leq d(x_n,x_m)\leq \lambda |n-m| + k$. The sequence $(x_n)_{n \in \Z}$ is a \emph{$(\lambda,k,L)$-local-quasi-geodesic} if we have $\frac{1}{\lambda}|n-m|-k \leq d(x_n,x_m)\leq \lambda |n-m| + k$ whenever $|n-m|\leq L$. Precisely we prove the following lemma :

\begin{Lemma} \label{unif-quasi-geod}
Let $X$ be a $\delta$-hyperbolic, geodesic and visibility space, and $o \in X$ any basepoint. Pick $A$ and $B$ two hyperbolic isometries of $X$ and suppose that $B(A^+) \neq A^-$. Then, there exists $\lambda > 0, k\geq0$ and $N \in \mathbb{N}^*$, such that $\forall G=(G_n)_{n \in \mathbb{Z}} \in \mathcal{G}_N(A,B)$, the sequence of points $x_n=G_no$ is a $(\lambda,k)$-quasi-geodesic. 
\end{Lemma}

For this purpose, we will use the Local-Global Lemma, which enables us to pass from local-quasi-geodesicity to global quasi-geodesicity, under the assumption of hyperbolicity. We recall it hereafter :

\begin{Lemma}[Local-Global, \cite{coornaert_geometrie_1990}, Chapter 3, Theorem 1.4] \label{local-global}
Let $X$ be a geodesic $\delta$-hyperbolic space. For all pairs $(\lambda,k)$, with $\lambda \geq 1$ and $k \geq 0$, there exists a real number $L$ and a pair $(\lambda',k')$ such that every $(\lambda,k,L)$-local-quasi-geodesic is a $(\lambda',k')$-quasi-geodesic (global). Moreover, $\lambda',k'$ and $L$ only depend on $\delta, \lambda$ and $k$.
\end{Lemma}

\begin{proof}
$\bullet$ \textbf{Step 1 : Quasi-isometry on a period} \\
The goal is at first to show that there exists two constants $\lambda > 0$ and $k \geq 0$, only depending on $\delta, A, B$ and $o$, such that the following inequality is satisfied : 
\begin{equation}\label{step1}
    \frac{1}{\lambda}|A^nBA^m|-k \leq d(A^nBA^mo,o) \qquad \text{ for all } n,m \geq 0
\end{equation}

By hypothesis, the two points at infinity $B(A^+)$ and $A^-$ are distinct, so we can consider a geodesic, called $\Lambda$, with endpoints $B(A^+)$ and $A^-$. Such a geodesic exists because $X$ is supposed to be a visibility space. Now consider $p$ a projection map on $\Lambda$, that is $p : X \to \Lambda$ satisfying $\forall x \in X, d(x,p(x))=d(x,\Lambda)=\underset{y \in \Lambda} \inf \, d(x,y)$ (such a map exists but is not necessarily unique). Since $(A^{-n}o)_{n \in \N}$ is a (half) quasi-geodesic with endpoint $A^-$ and $\Lambda$ is a geodesic with $A^-$ as one of its endpoints, we have, by stability of quasi-geodesics in $\delta$ hyperbolic spaces, the existence of a constant $K_1 >0$ (only depending on $\delta,A,B$ and $o$) such that $\{A^{-n}o\}_{n \in \N}$ and the half geodesic $[p(o),A^{-})$ remain in the $K_1$-neighborhood of each other. We deduce the following inequality :
\begin{equation} \label{K1}
    d(A^{-n}o,p(A^{-n}o)) \leq K_1, \text{ for all } n \in \N
\end{equation}
With the same argument, namely that the (half) geodesic $(BA^mo)_{m \in \N}$ and $\Lambda$ share the same endpoint $B(A^+)$, we deduce the existence of constant $K_2 >0$ (only depending on $\delta, A, B$ and $o$) such that 
\begin{equation} \label{K2}
    d(BA^mo,p(BA^mo)) \leq K_2,  \text{ for all } m \in \N
\end{equation}

Then we can draw the following inequalities :
\begin{align*}
    d(A^nBA^mo,o)& =d(BA^mo,A^{-n}o) \text{ because $A^n$ is an isometry } \\
                & \geq d(p(BA^mo),p(A^{-n}o)) -d(p(BA^mo),BA^mo)-d(p(A^{-n}o),A^{-n}o) \\
                 &  \geq d(p(BA^mo),p(A^{-n}o)) -K_1 -K_2 \text{ by inequalities \ref{K1} and \ref{K2}} \\
\end{align*}

But since $A^{-n}o \underset{n\to \infty} \longrightarrow A^-$, we also have $p(A^{-n}o) \underset{n\to \infty} \longrightarrow A^-$, and in the same way, since $BA^mo \underset{m\to \infty} \longrightarrow B(A^+)$, we deduce $p(BA^mo) \underset{m\to \infty} \longrightarrow B(A^+)$. Then, for $n$ and $m$ sufficiently large, $p(A^{-n}o)$ belongs to $[p(o),A^-) \cap [p(Bo),A^-)$ and $p(BA^mo)$ belongs to $[p(o),B(A^+)) \cap [p(Bo),B(A^+))$. This shows that for $n$ and $m$ sufficiently large, the four points $p(A^{-n}o), p(BA^mo),p(Bo)$ and $p(o)$ are aligned in one of the two following orders on the geodesic $\Lambda$ : $p(A^{-n}o),p(o),p(Bo), p(BA^mo)$ or $p(A^{-n}o),p(Bo),p(o), p(BA^mo)$. In the first case $$d(p(A^{-n}o),p(BA^mo))=d(p(A^{-n}o),p(o))+d(p(o),p(Bo))+d(p(Bo),p(BA^mo))$$ and in the second one : $$d(p(A^{-n}o),p(BA^mo))=d(p(A^{-n}o),p(o))-d(p(o),p(Bo))+d(p(Bo),p(BA^mo))$$
so, in every case, for $n$ and $m$ sufficiently large :
\begin{align*}
   d(p(A^{-n}o),p(BA^mo)) & \geq d(p(A^{-n}o),p(o))-d(p(o),p(Bo))+d(p(Bo),p(BA^mo)). 
\end{align*}
On an other hand, 
\begin{align*}
    d(p(A^{-n}o),p(o)) & \geq d(A^{-n}o,o) - d(A^{-n}o,p(A^{-n}o))-d(p(o),o) \\ 
    & \geq d(A^{-n}o,o)-K_1-d(p(o),o) \text{ by inequality \ref{K1}}
\end{align*}
and similarly :
\begin{align*}
    d(p(BA^mo),p(Bo)) & \geq d(BA^mo,Bo) - K_2  -d(p(Bo),Bo) \text{ by inequality \ref{K2}}.
\end{align*}
We can now finish our sequence of inequalities :
\begin{align*}
    d(A^nBA^mo,o)& \geq d(p(BA^mo),p(A^{-n}o)) -K_1 -K_2 \\
    & \geq d(p(A^{-n}o),p(o))-d(p(o),p(Bo))+d(p(Bo),p(BA^mo)) -K_1-K_2 \\
    & \geq d(A^{-n}o,o)-d(p(o),o)+d(BA^mo,Bo)-d(p(Bo),Bo)-d(p(o),p(Bo))-2K_1-2K_2 \\
    & = d(A^no,o)+d(A^mo,o)-d(p(o),o)-d(p(Bo),Bo)-d(p(o),p(Bo))-2K_1-2K_2 \\
    & \geq (n+m)l_S(A)-d(p(o),o)-d(p(Bo),Bo)-d(p(o),p(Bo))-2K_1-2K_2
\end{align*}
In the last inequality, we used the basic fact that $d(A^no,o) \geq nl_S(A)$, where $l_S(A)$ denotes the stable length of the isometry $A$. 
Since $(n+m)l_S(A)=(n+m+1)l_S(A)-l_S(A)=|A^nBA^m|l_S(A)-l_S(A)$, we have proved the inequality \ref{step1} for $n$ and $m$ sufficiently large, prescribing $\lambda=\frac{1}{l_S(A)}$ (recall $l_S(A)>0$ when $A$ is hyperbolic), and $k=l_S(A)+ d(p(o),o)+d(p(Bo),Bo)+d(p(o),p(Bo))+2K_1+2K_2$. But there is only a finite number of value of $A^nBA^m$, for $n$ and $m$ smaller than a fixed constant, so the inequality \ref{step1} is still true for all $n,m \in \N$, after possibly changing the value of $\lambda$ and $k$.  \\

$\bullet$ \textbf{Step 2 : From local to global quasi-isometry}\\

We may now conclude using the local-global lemma (Lemma \ref{local-global}). \\ 

Let $L>0$ and $(\lambda',k')$ such as in Lemma \ref{local-global}, with $\lambda$ and $k$ defined in the first step. Fix $N=\lfloor L \rfloor +1$. Then every interval of length smaller than $L$ is of length smaller than $N$. Now choose $G$ a sequence in $\mathcal{G}_N(A,B)$ which is associate to a bi-infinite words $H=(H_n)_{n\in \Z} \in \mathcal{H}_N(A,B)$. Thus, the subwords of $H$ of length smaller than $L$ are of the form $A^nBA^m$ or $A^n$, with $n,m \in \mathbb{N}$. Therefore, by Step 1, the sequence of points $(x_n)_{n\in \Z}$ is a  $(\lambda,k,L)$-local-quasi-geodesic. So, by the local-global lemma \ref{local-global}, there exists $\lambda' \geq 1,k'\geq 0$ (only depending on $\lambda$ and $k$, that is on $\delta, A, B$ and $o$), such that $(x_n)_{n\in \Z}$ is a  $(\lambda',k')$-quasi-geodesic (global). Thus, the proposition \ref{unif-quasi-geod} is proved. 

\end{proof}

\section{Bowditch and primitive-stable representations of $\F_2$ in a $\delta$-hyperbolic space}
\label{Bowditch}

\subsection{A first property of Bowditch representations} ~\\ 

Let $(X,d)$ be a $\delta$-hyperbolic, geodesic and visibility space, and $o \in X$ a basepoint. The stable length (defined in the previous section) is the right notion to determine whether an isometry is hyperbolic or not : it can be shown that an isometry $A$ is hyperbolic if and only if $l_S(A)>0$, whereas this equivalence is not true in general when considering the displacement length. Finally, we can compare the displacement and the stable length, with the following inequality : \begin{equation}\label{comp-l-ls}
    l_S(A) \leq l(A) \leq l_S(A)+16\delta
\end{equation}
The left inequality follows directly from the definitions, hence is true in any metric space, whereas the right inequality is really a feature of $\delta$-hyperbolicity. 
The proofs of all these facts about the stable length can be found in \cite{coornaert_geometrie_1990}, Chapter 10.6.\\

\begin{Remark}
Since $l_\rho(\gamma) = \underset{o \in X} \inf d(\rho(\gamma)o,o)$, we deduce that a Bowditch representation of constants $(C,D)$ satisfies, for any basepoint $o \in X$ : 
\begin{equation*}
   \forall \gamma \in \mathcal{P}(\mathbb{F}_2), \hspace{1cm} \frac{1}{C} \Vert \gamma \Vert -D \leq d(\rho(\gamma)o,o)
\end{equation*}

Moreover, because of the inequalities \eqref{comp-l-ls}, we can deduce that a Bowditch representation of constants $(C,D)$ also satisfies : 
\begin{equation*}
    \forall \gamma \in \mathcal{P}(\mathbb{F}_2), \hspace{1cm} \frac{1}{C} \Vert \gamma \Vert -D-2\delta  \leq  l_S(\rho(\gamma))
\end{equation*}
Hence, in definition \ref{def-bowditch}, we could also use the stable length instead of the displacement length. 
\end{Remark}

Now we establish the useful fact that the image of primitive elements by a Bowditch representation are hyperbolic isometries. 
\begin{Lemma} \label{gamma-hyp}
Let $\rho$ be a Bowditch representation of constants $(C,D)$. Then, for every primitive element $\gamma$ in $\F_2$, $ \frac{1}{C}\Vert \gamma \Vert \leq l_S(\rho(\gamma))$. \\ 
In particular, for every primitive element $\gamma$, $\rho(\gamma)$ is hyperbolic and $\rho$ is also a Bowditch representation of constant $(C,0)$.
\end{Lemma}

\begin{proof} Suppose that $\gamma$ is cyclically reduced. The primitivity hypothesis on $\gamma$ gives the existence of another primitive element, $\delta \in \mathbb{F}_2$, such that we have both $\{\gamma,\delta \}$ is a free basis of $\F_2$, and $\Vert \gamma^n\delta \Vert=n \Vert \gamma \Vert + | \delta | $. Thus, for all $n \in \mathbb{N}$, the element $\gamma^n\delta$ is primitive. The Bowditch inequality applied to  $\gamma^n\delta$ gives :

\begin{align*}
    \frac{1}{C}\Vert \gamma^n\delta \Vert -D & \leq l_\rho(\gamma^n\delta) \leq d(\rho(\gamma^n\delta)o,o) \text{ by definition of the displacement length} \\
    & \leq d(\rho(\gamma^n)o,o) + d(\rho(\delta)o,o) \text{ by the triangle inequality} 
\end{align*}
Using that $\Vert \gamma^n\delta \Vert=n \Vert \gamma \Vert + | \delta | $, and after dividing by $n$, we obtain :
\begin{align*}
    \frac{1}{C}\Vert \gamma \Vert  + \frac{|\delta|}{nC}-\frac{D}{n} \leq \frac{1}{n}d(\rho(\gamma)^no,o)+\frac{1}{n}d(\rho(\delta)o,o)
\end{align*}
Now, let $n$ tends to infinity and use the definition of the stable length : 
\begin{align*}
    \frac{1}{C}\Vert \gamma \Vert  \leq l_S(\rho(\gamma)).
\end{align*}
which is indeed the desired inequality. \\
Thus, the stable length of $\rho(\gamma)$ is positive, we deduce that $\rho(\gamma)$ is hyperbolic.
At last, since the stable length is always smaller than or equal to the displacement length (see inequality \eqref{comp-l-ls}), we also deduce the inequality $\displaystyle \frac{1}{C} \Vert \gamma \Vert \leq l_\rho(\gamma)$, which finishes the proof. 
\end{proof}

\subsection{The inclusion $\mathcal{PS}(\F_2,X) \subset \mathcal{BQ}(\F_2,X)$} ~\\
\label{first-inclusion}

Recall that we defined in definition \ref{def-primitive-stable} primitive-stability. We will now see in this section a first inclusion between primitive-stable representations and Bowditch representations. It is quite easy to check that primitive-stable representations are in particular Bowditch :  

\begin{Lemma} \label{inclusion-easy}
Let $\rho : \F_2 \to \mathrm{Isom}(X)$ be a primitive-stable representation. Then $\rho$ is a Bowditch representation. 
\end{Lemma}

\begin{proof}
Let $o \in X$ be a basepoint. Then, there exist two constants $C$ and $D$ such that $\rho$ is primitive-stable with constants $(C,D)$. Let $\gamma \in \mathcal{P}(\F_2)$ be a cyclically reduced primitive element in $\F_2$ and $n \in \N$. The elements 1 and $\gamma^n$ both belong to the geodesic $L_\gamma$ in the Cayley graph of $\F_2$, therefore~:
\begin{align*}
 & \frac{1}{C}\Vert \gamma^n \Vert -D  \leq d(\tau_\rho(\gamma^n),\tau_\rho(1)) =d(\rho(\gamma)^no,o), \\
 \text{then, dividing by $n$, } \qquad & \frac{1}{C}\Vert \gamma \Vert -\frac{D}{n} \leq \frac{1}{n}d(\rho(\gamma)^no,o), \\
 \text{and taking the limit when $n \to \infty$,} \qquad & \frac{1}{C}\Vert \gamma \Vert \leq l_S(\rho(\gamma)). \\
 \text{Using the inequality \eqref{comp-l-ls},} \qquad  & \frac{1}{C}\Vert \gamma \Vert \leq l(\rho(\gamma)).
\end{align*}
Hence $\rho$ is a Bowditch representation. 
\end{proof}

\subsection{A lemma on Bowditch representations}

\begin{Lemma} \label{hyp-gen}
Let $X$ be a $\delta$-hyperbolic, geodesic and visibility space, and $\rho : \mathbb{F}_2 \to \mathrm{Isom}(X)$ a Bowditch representation. Fix $\{a,b\}$ a free basis of $\mathbb{F}_2$ and denote by $A=\rho(a),B=\rho(b)$ the images of the generators by $\rho$. \\
Then $B(A^+) \neq A^-$ (where $A^+$ and $A^-$ refer respectively to the attracting and repelling fixpoints of $A$).
\end{Lemma}

\begin{proof}
Before starting the proof, recall that we have shown that $\rho(a)$ and $\rho(b)$ are hyperbolic isometries (because $a$ and $b$ are primitive elements, see lemma \ref{gamma-hyp}), therefore $A^+$ and $A^-$ are well-defined. \\
Let us chose some basepoint $o \in X$. Then, because $A$ is hyperbolic, the sequence $(A^n)_{n \in \Z}$ is a quasi-isometry with repelling fixpoint $A^-$. Furthermore, the sequence $(BA^m)_{m \in \Z}$ is again a quasi-isometry, with attracting fixpoint $B(A^+)$. Now suppose by absurdity that $B(A^+)=A^-$. The stability of quasi-geodesics in $\delta$-hyperbolic spaces then gives the existence of a constant $K>0$ such that the half-geodesics $(A^{-n})_{n \in \N}$ and $(BA^m)_{m \in \N}$ stay at a distance $K$ of each other. Thus, we deduce the existence, for all $n \in \N$, of an integer $\phi(n) \in \N$ such that $d(A^{-n}o,BA^{\phi(n)}o) \leq K$. But the element $a^nba^{\phi(n)} \in \F_2$ is primitive, so by the Bowditch hypothesis, we have the following inequality :
\begin{align*}
    \frac{1}{C}\Vert a^nba^{\phi(n)} \Vert -D \leq d(\rho(a^nba^{\phi(n)})o,o) = d(A^nBA^{\phi(n)}o,o) = d(BA^{\phi(n)}o,A^{-n}o) 
\end{align*}
Here the right hand side of the inequality is bounded by $K$, and the left hand side tends to infinity because $\Vert a^nba^{\phi(n)}\Vert =n+\phi(n)+1$, this is a contradiction.
\end{proof}

\subsection{Openness of the set of primitive-stable representations} ~\\
\label{openness}

Here we want to prove that we can deform primitive-stable representations, in other words that the set of primitive-stable representations is open in the character variety. Recall that the primitive-stability condition is invariant under conjugacy, hence the notion of primitive-stability is well-defined in the character variety. Note that although we write the proof in the case of the group $\F_2$ since our focus is on $\F_2$, the proof works the same more generally for higher rank free groups $\F_n$ with $n \geq 2$.

\begin{Proposition}
$\mathcal{PS}(\F_2,X)$ is open in the character variety $\chi(\F_2,\mathrm{Isom}(X))$.
\end{Proposition}

\begin{proof}
Let $\rho : \F_2 \longrightarrow \mathrm{Isom}(X)$ be a primitive-stable representation.  Denote by $(C,D)$ the two constants of primitive-stability of $\rho$. We want to find an open neighborhood of $\rho$ in $\mathrm{Hom}(\F_2,\mathrm{Isom}(X))$ consisting only of primitive-stable representations. Our open set will be of the following type : \\
For $L>0$ and $\varepsilon>0$ two positive constants, define $$\mathcal{V}_\rho(L,\varepsilon) = \{\rho' : \F_2 \to \mathrm{Isom}(X) \, | \, \forall u \in \F_n \quad |u| \leq L \implies d(\rho(u)o,\rho'(u)o) < \varepsilon \}. $$
Recall that $\mathrm{Hom}(\F_2,\mathrm{Isom}(X))$ is endowed with the compact-open topology, then $\mathcal{V}_\rho(L,\varepsilon)$ is an open subset of  $\mathrm{Hom}(\F_2,\mathrm{Isom}(X))$. In the following, we will use the local-global lemma, (which we have recalled previously in Lemma \ref{local-global}). \\

Now let us fix $\varepsilon=1$ (we could have chosen any other value for $\varepsilon$). The local-global lemma gives the existence of three constants $L,C',D'$ such that any $(C,D+1,L)$-local-quasi-geodesic is a $(C',D')$-quasi-geodesic. Consider $\mathcal{V}_\rho(L,1)$. It is an open neighborhood of $\rho$. We will now show that $\mathcal{V}_\rho(L,1)$ consists only of primitive-stable representations. Indeed, take $\rho' \in \mathcal{V}_\rho(L,1)$ and let $\gamma$ be a primitive element in $\F_2$ and $u$ and $v$ two integer points on the geodesic $L_\gamma$ of the Cayley graph of $\F_2$ (recall that $L_\gamma$ is the geodesic of the Cayley graph generated by $\gamma$). Then
\begin{align*}
    \frac{1}{C}d(u,v)-D & \leq d(\tau_\rho(u),\tau_\rho(v)) \qquad \text{ because 
    $\rho$ is primitive-stable } \\
    & \leq d(\rho(u)o,\rho(v)o) \qquad \text{ because $u$ and $v$ are integer points} \\ 
    & \leq d(\rho(v^{-1}u)o,o) \qquad \text{ because $\rho(v^{-1})$ is an isometry } \\
    & \leq d(\rho(v^{-1}u)o,\rho'(v^{-1}u)o)+d(\rho'(v^{-1}u)o,o) \qquad  \text{ by the triangle inequality} \\
    & \leq 1 + d(\rho'(u)o,\rho'(v)o) \qquad \text{ because $|v^{-1}u|=d(u,v) \leq L$ and $\rho' \in \mathcal{V}_\rho(L,1)$}
\end{align*}
So we deduce the inequality :
\begin{equation}
    \frac{1}{C}d(u,v)-D-1 \leq d(\rho'(u)o,\rho'(v)o)
\end{equation}
which shows that $\rho'(L_\gamma)$ is a $(C,D+1,L)$ local-quasi-geodesic, hence a $(C',D')$ quasi-geodesic from the local-global lemma. This shows that $\rho'$ is primitive-stable and thus that $\mathcal{V}_\rho(L,1)$ is an open-neighborhood of $\rho$ consisting only of primitive-stable representations. Its image under the projection to $\chi(\F_2,\mathrm{Isom}(X))$ is again an open neighborhood consisting of primitive-stable representations and thus $\mathcal{PS}(\F_2,X)$ is open.
\end{proof}

\section{From Bowditch's hypothesis to uniform tubular neighborhoods}
\label{first-step-proof}
The purpose of this section is to show Proposition \ref{uniform-neighborhood}, which is the heart of the proof that a Bowditch representation is primitive-stable. \\

Before stating the proposition and starting the proof, recall that when $A$ is a hyperbolic isometry of $X$, it defines two points in the boundary of $X$, $A^+$ and $A^-$, respectively attracting and repelling fixpoints of the action of $A$ on $\partial X$. Let's denote by $\mathrm{Axis}(A)$ the union of all the geodesics of $X$ joining the two points $A^+$ and $A^-$. Since $X$ is a visibility space, this set in by assumption non-empty. When $X=\HH^n$ the usual hyperbolic space of dimension $n$, the geodesic joining $A^+$ and $A^-$ is unique and corresponds to the usual definition of the axis of the hyperbolic isometry $A$. The set $\mathrm{Axis}(A)$ in invariant under $A$ : indeed, for every geodesic $\ell$ joining $A^+$ and $A^-$ in $X$, $A(\ell)$ is still a geodesic because $l$ is a geodesic and $A$ an isometry. Now using the fact that the endpoints of $\ell$ are the fixpoints at infinity of $A$, we deduce that the endpoints of $A(\ell)$ are also $A^+$ and $A^-$, so $A(\ell) \subset \mathrm{Axis}(A)$, thus $\mathrm{Axis}(A)$ is $A$-invariant. From the $A$-invariance of $\mathrm{Axis}(A)$ also follows the $A$-invariance of the map $d(\cdot , \mathrm{Axis}(A) )$. 

For a subset $Y$ of $X$ and $K>0$, denote by $N_K(Y)$ the $K$-neighborhood of $Y$, that is $N_K(Y)=\{x \in X \, : \, d(x,Y) \leq K \}$. Fix $\ell$ any geodesic of $X$ joining $A^+$ and $A^-$. Then we have the following lemma :
\begin{Lemma} \label{nbh-axis}
\begin{itemize}
    \item $N_K(\ell) \subset N_K(\mathrm{Axis}(A))$
    \item There exists a constant $C(\delta)$, depending only on the hyperbolic constant $\delta$, such that~ \\ $ \displaystyle N_K(\mathrm{Axis}(A)) \subset N_{K+C(\delta)}(\ell)$.
\end{itemize}
\end{Lemma}

\begin{proof} \label{neighborhood}
\begin{itemize}
    \item The first point is immediate because $\ell \subset \mathrm{Axis}(A)$.
    \item 
The second point basically follows from the Morse lemma. Since $X$ is $\delta$-hyperbolic, there exists a constant $C(\delta)$, depending only on $\delta$, such that any two geodesics with the same endpoints remain at a distance $C(\delta)$ of each other. \\
Then, if $x \in N_K(\mathrm{Axis}(A))$, there exists $y \in \mathrm{Axis}(A)$ such that $d(x,y) \leq K$. But since $y \in \mathrm{Axis}(A)$, in particular, $y$ belongs to a geodesic with endpoints $A^+$ and $A^-$, let's denote it by $\ell_y$. Thus $\ell_y$ and $\ell$ remain at a distance $C(\delta)$ of each other and thus $x$ is at distance at most $K+C(\delta)$ of $\ell$.
\end{itemize}
\end{proof}

Now this section is dedicated to proving Proposition \ref{uniform-neighborhood}, which shows that the Morse lemma is satisfied for the primitive elements of a Bowditch representation, meaning that the orbit map restricted to primitive leaves stays in a uniform tubular neighborhood of the axis of primitive elements in $X$. 

\begin{Proposition}
\label{uniform-neighborhood}
Let $\rho : \mathbb{F}_2 \to \mathrm{Isom}(X)$  be a Bowditch representation. The orbit map restricted to primitive leaves stays in a uniform tubular neighborhood of the axis of primitive elements in $X$. Precisely :
\begin{equation*}
    \exists K>0, \hspace{1cm} \forall \gamma \in \mathcal{P}(\mathbb{F}_2), \hspace{1cm} \tau_\rho(L_{\gamma}) \subset N_K(\mathrm{Axis}(\rho(\gamma)))
\end{equation*}
\end{Proposition}

Recall that $L_\gamma$ denotes the (geodesic) axis of $\gamma$ in the Cayley graph of $\mathbb{F}_2$, and that for any primitive element $\gamma$ in $\F_2$, $\rho(\gamma)$ is hyperbolic by Lemma \ref{gamma-hyp} so $\mathrm{Axis}(\rho(\gamma))$ is well-defined.\\

\begin{proof}
Pick $\rho$ a Bowditch representation and let $C >0, C' >0$ be two constants such that
$$ \forall \gamma \in \mathcal{P}(\mathbb{F}_2), \hspace{0,5cm} \frac{1}{C}\Vert \gamma \Vert \leq l(\rho(\gamma)) \hspace{0,5cm} \text{ and } \hspace{0,5cm} \forall u \in \mathbb{F}_2, \hspace{0,5cm} d(\rho(u)o,o) \leq C' |u|. $$ 
Note that such constants automatically satisfy $CC' \geq 1$. \\
Let us proceed by contradiction and suppose there exists a sequence $(\gamma_n)_{n \in \mathbb{N}}$ of cyclically reduced primitive elements of $\F_2$ satisfying the following hypothesis :

\begin{equation*}
\sup \left\{ d(x,\mathrm{Axis}(\rho(\gamma_n))) : x \in \tau_\rho(L_{\gamma_n}) \right\} \underset{n \to \infty} \longrightarrow +\infty \tag{H} \label{hyp-absurde}
\end{equation*}
We fix such a sequence $(\gamma_n)_{n \in \mathbb{N}}$ for all that follows. Without loss of generality, we can assume that the elements $\gamma_n$ are pairwise distinct and that $\Vert \gamma_n \Vert  \to \infty$. ~ \\
Indeed, for all $N \in \mathbb{N}$, the set $\Gamma_N = \{ n \in \mathbb{N} : \gamma_n=\gamma_N \}$ is finite : if this was not true, there would exist a subsequence $(\gamma_{\sigma(n)})_{n \in \mathbb{N}}$ such that $\sup \{ d(x,\mathrm{Axis}(\rho(\gamma_N))) : x \in \tau_\rho(L_{\gamma_N}) \}=\sup \{ d(x,\mathrm{Axis}(\rho(\gamma_{\sigma(n)}))) : x \in \tau_\rho(L_{\gamma_{\sigma(n)}}) \}$ and this would contradict the hypothesis \eqref{hyp-absurde}. Therefore, after passing to a subsequence, we can assume that the elements $(\gamma_n)_{n \in \mathbb{N}}$ are pairwise distinct. \\
Moreover, for all $A>0$, $\{ \gamma \in \mathbb{F}_2 : |\gamma| \leq A \}$ is finite, so, since the elements $(\gamma_n)_{n \in \mathbb{N}}$ are pairwise distinct and cyclically reduced, we also have the finiteness of the set $\{n \in \mathbb{N} : \Vert \gamma_n \Vert \leq A \}$ for all $A >0$. Then for $n$ sufficiently large, $\Vert \gamma_n\Vert \geq A$, hence  $\Vert \gamma_n\Vert  \underset{ n \to \infty} \longrightarrow \infty $.

\subsection{Continued fraction expansion of $\gamma_n$} ~\\
The element $\gamma_n$ is primitive, thus corresponds to a rational (see section \ref{structure-F2}), we can then write its continued fraction expansion : $\gamma_n = [N_1^n,N_2^n,\hdots,N_{r(n)}^n]$.

\begin{Lemma} \label{forme-frac-continue}
Up to subsequence, $r(n) \to + \infty$ and for all $i \in \mathbb{N}$, $(N_i^n)_{n \in \mathbb{N} \hspace{0,1cm} : \hspace{0,1cm} r(n) \geq i}$ is bounded. 
\end{Lemma}

\begin{proof}

Suppose there exists $i \in \mathbb{N}$ such that $(N_i^n)_n$ is defined for an infinity of $n$ and is not bounded. Then consider the smallest such $i \in \mathbb{N}$. For all $1 \leq j < i$, the sequence $(N_j^n)_n$ is a bounded sequence of integers so after passing to subsequence we assume that there exists an integer $N_j$ such that for all $1 \leq j < i$ and for all $n \in \mathbb{N}$ such that $r(n) \geq j$, $N_j^n=N_j$. Thus $\slope(\gamma_n) = [N_1,N_2,\hdots,N_{i-1},N_i^n,\hdots,N_{r(n)}^n]$. We set $u=w_{i-1}(\gamma_n)$ and $v=w'_{i-1}(\gamma_n)$ (see definition \ref{def-wi} in section \ref{structure-F2}). Therefore, $\slope(u)=[N_1,N_2,\hdots,N_{i-1}]$, $\slope(v)=[N_1,N_2,\hdots,N_{i-1}+1]$ (by Proposition \ref{wi}) and $u$ and $v$ do not depend on the integer $n$. Moreover, $u$ and $v$ form a free basis of $\mathbb{F}_2$ such that (up to  cyclic permutation and inversion) $\gamma_n$ is a positive word on $u$ and $v$ (by Proposition \ref{wi}). Denote by $U=\rho(u)$ and $V=\rho(v)$ their images by $\rho$, then by the lemma \ref{hyp-gen}, we conclude $V(U^+)\neq U^-$. Then, by considering the bi-infinite word obtained by concatenating infinitely many copies of $\gamma_n$, or equivalently the bi-infinite word obtained by following the geodesic $L_{\gamma_n}$ in the Cayley graph, we can see $\rho_{|L_{\gamma_n}}$ as an element of $\mathcal{G}(U,V)$ (the definition is given at the beginning of section \ref{UQG-setting}). Thus, we define $N \in \mathbb{N}^*$ as in the lemma \ref{unif-quasi-geod} (depending on $\delta,U,V$ and the basepoint $o$) and since by hypothesis $N_i^n \to +\infty$, $\rho_{|L_{\gamma_n}}$ is a sequence of $\mathcal{G}_N(U,V)$ for $n$ sufficiently large. Then, using the lemma \ref{unif-quasi-geod}, we obtain the existence of two constants $\lambda>0$ and $k\geq 0$ (only depending on $\delta,U,V$ and the basepoint $o$) such that $\tau_{\rho |L_{\gamma_n}}$ is a  $(\lambda,k)$-quasi-geodesic. The Morse lemma now gives the existence of a constant $K>0$ only depending on $\lambda$ and $k$ such that $\tau_\rho(L_{\gamma_n})$ remains in the $K$-neighborhood of $\mathrm{Axis}(\rho(\gamma_n))$. This contradicts our hypothesis \eqref{hyp-absurde} on $\rho$ for $n$ sufficiently large. Hence, for all $i \in \mathbb{N}$, $(N_i^n)_n$ is bounded. \\

Let us now justify that $r(n) \to + \infty$. If $r(n)$ stays bounded, $r(n)\leq R$, then for all $1 \leq i \leq R$, $(N_i^n)_n$ is bounded by what has been previously done and so the word length of $\gamma_n$ is also bounded, which is false. Thus $r(n) \to +\infty$. \\
In particular, we deduce that under the assumption \eqref{hyp-absurde}, the sequence  $(N_i^n)_n$ is always well-defined for $n$ sufficiently large ($n$ such that $r(n) \geq i$).\\
\end{proof}

\subsection{Uniform bounds on the lengths $l_i(\gamma_n)$} ~\\
\label{bornes-li}
Using the notations of the previous sections and the definition of $l_i(\gamma_n)$ given in Definition \ref{def-wi}, we have, using the inequalities \ref{ineg-li} and \ref{min-li} of the section \ref{structure-F2} together with the upper bound $N_i^n \leq N_i$ : 

\begin{equation*}
    \forall n \in \mathbb{N}, \forall 0<i\leq r(n), \hspace{1cm} 1 + \frac{1}{N_{i-1}+1} \leq \frac{l_i(\gamma_n)}{l_{i-1}(\gamma_n)} \leq N_i+1
\end{equation*}
We deduce, since for any integer $n$, $l_0(n)=1$, that for any integer $i$, there exists a positive constant $L_i >0$ such that :

\begin{equation} \label{encadr-li(n)}
    \forall n\in \mathbb{N}, \qquad \forall 0 \leq i \leq r(n), \hspace{1cm} i \leq l_i(\gamma_n) \leq L_i 
\end{equation}

\subsection{Excursions for real map}
\label{excursion-real-map}

\begin{Definition} 
An \emph{excursion} is the data of two reals $a \leq b$ and of a map $E : [a,b] \to \mathbb{R}$ which satisfies : 
\begin{itemize}
    \item $E$ is continuous on $[a,b]$
    \item $E(a)=E(b)$
    \item $\forall t \in [a,b], E(t) \geq E(a)$
\end{itemize}
We define the \emph{length of excursion} of $E$ as the non-negative real $b-a$. \\
Furthermore, the map $E : [a,b] \to \mathbb{R}$ is said to be a \emph{$K$-excursion} if $E$ is an excursion such that $E(a)=K$.  
\end{Definition}

\begin{Definition}
Let $E : [a,b] \to \mathbb{R}$ be an excursion. We say that $E'$ is a \emph{sub-excursion} of $E$ if there exists a subinterval  $[c,d] \subset [a,b]$ such that $E'=E_{|[c,d]}$ and $E'$ is an excursion. \\
Furthermore, $E'$ is said to be a \emph{$K$-sub-excursion} of $E$ if $E'(c)=K$.
\end{Definition}

The goal of this section is Lemma \ref{temps-excursion}, which shows that an excursion always has sub-excursions of any prescribed length up to a factor of 2. 
\begin{Remark}
Trivially, if $E$ is an excursion, $E$ is a sub-excursion of itself and for all $c \in [a,b], E_{|[c,c]}$ is also a sub-excursion of $E$.
\end{Remark}

\begin{Lemma}
\label{K-sous-exc}
Let $E : [a,b] \to \mathbb{R}$ be an excursion. We set  $K_{\min} =\min E = E(a)=E(b)$ and $K_{\max} = \max E$.
Then, for all $K \in [K_{\min},K_{\max}]$, there exists a $K$-sub-excursion of $E$.
\end{Lemma}

\begin{proof}
Let $K_{\min} \leq K \leq K_{\max}$. Choose $c \in [a,b]$ such that $E(c)=K_{\max}$. We denote $X_L=E^{-1}(K) \cap [a,c]$ and $X_R=E^{-1}(K) \cap [c,b]$. The sets $X_L$ and $X_R$ are closed (by continuity of $E$) and non-empty (by the intermediate value theorem) so we can consider
\begin{align*}
    x_K=\max X_L \text{ et } y_K=\min X_R
\end{align*}
Then $E_{|[x_K,y_K]}$ is a $K$-sub-excursion of $E$.
\end{proof}

\begin{Lemma} \label{exc-ferme}
Let $l>0$ and $E$ be an excursion of length $l$. Let $T_E$ be the set of all lengths of excursion of sub-excursion of $E$, that is : $$T_E = \{0 \leq l' \leq l : \text{there exists a sub-excursion of $E$ of length $l'$} \} $$
Then $D_E$ is a closed subset of $[0,l]$.
\end{Lemma}

\begin{Remark}
By the previous remark, we always have $0 \in T_E, l \in T_E$.
\end{Remark}

\begin{proof} 
Let $E : [a,b] \to \mathbb{R}$ be an excursion of length $l$, which means that $E$ is continuous, $E(a)=E(b)$, $\forall t \in [a,b], E(t) \geq E(a)$ and $b-a=l$. \\
Let $(l_n)_{n\in \mathbb{N}}$ be a sequence of $T_E$ such that $l_n \to l_\infty \in [0,l]$. \\
Let $(a_n)_{n\in \mathbb{N}}$ and $(b_n)_{n \in \mathbb{N}}$ be two sequences of $[a,b]$ such that $E : [a_n,b_n] \to \mathbb{R}$ is a sub-excursion of length $l_n$. Up to subsequence, since $[a,b]$ is compact, we can assume that $a_n \to a_\infty \in [a,b]$ and $b_n \to b_\infty \in [a,b]$. Moreover, using the continuity of $E$, $\forall n \in \mathbb{N}, E(a_n)=E(b_n)$ and $\forall t \in [a_n,b_n] E(t) \geq E(a_n)$, we obtain $E(a_\infty)=E(b_\infty)$ and $\forall t \in [a_\infty,b_\infty], E(t) \geq E(a_\infty)$. Finally, $l_\infty=\underset{n}\lim \, l_n = \underset{n} \lim \, (b_n-a_n) = b_\infty-a_\infty$ so $E : [a_\infty,b_\infty] \to \mathbb{R}$ is indeed a sub-excursion of length $l_\infty$.
\end{proof}

\begin{Lemma} \label{exc-plus-petite}
Let $l >0$ and $E$ an excursion of length $l$. Then there exists a sub-excursion of $E$ of length $l'>0$ such that $\displaystyle \frac{l}{2} \leq l' < l$.
\end{Lemma}

\begin{proof}
Let $E : [a,b] \to \mathbb{R}$ be an excursion of length $l$ (then $E$ is continuous, $E(a)=E(b)$, $\forall t \in [a,b], E(t) \geq E(a)$ and $b-a=l$). \\
We distinguish two cases : 
\begin{itemize}
    \item \textbf{$1^{\text{st}}$ case} : There exists $t \in ]a,b[, E(t)=E(a)$ : \\
    Then $E : [a,t] \to \mathbb{R}$ and $E : [t,b] \to \mathbb{R}$ are two sub-excursions of length $t-a$ and $b-t$ respectively. But either $t-a \geq \frac{b-a}{2}=\frac{l}{2}$ or $b-t \geq \frac{b-a}{2}=\frac{l}{2}$ so one of these two sub-excursion is in fact of length $\frac{l}{2} \leq l'< l$. \\
    
    \item \textbf{$2^{\text{nd}}$ case} : For all $t \in ]a,b[, E(t)>E(a)$ : \\
    Let $c,d \in ]a,b[$ such that $d-c \geq \frac{l}{2}$. The map $E$ is continuous on the segment $[c,d]$. Denote $\delta = \underset{[c,d]}\min E$. Then $\delta>E(a)$. Let $h=\frac{1}{2}(E(a)+\delta)$ and define :
    \begin{align*}
        a_h & = \max \{ a' \in [a,c] : E(a')=h \} \\
        b_h & = \min \{ b' \in [d,b] : E(b')=h \}
    \end{align*}
    The set $\{ a' \in [a,c] : E(a')=h \} $ is non-empty (because $E(c) \geq \delta > h > E(a) $ and $E$ is continuous) and closed, so $a_h$ is well-defined. Likewise, $b_h$ is well-defined. \\
    Therefore, we have :
    \begin{itemize}
        \item For all $t \in [c,d], E(t) > h$ because $h<\delta=\underset{[c,d]}\min E$
        \item For all $t \in [a_h,c], E(t) \geq h$ : indeed, if there  was $t \in [a_h,c]$ such that $E(t)<h$, then on one hand $t \in ]a_h,c[$, and on the other hand, since $E(c)>h$, by the intermediate value theorem, there would exist $a' \in [t,c]$ such that $E(a')=h$ and $a'>a_h$, which is impossible because of the choice of $a_h$.
        \item For all $t \in [d,b_h], E(t) \geq h$ : the same argument as above works. \\
    \end{itemize}
    Therefore, $E : [a_h,b_h] \to \mathbb{R}$ is an excursion of length $l'=b_h-a_h$ which satisfies $$\frac{l}{2} \leq d-c \leq b_h-a_h=l'< b-a =l$$
\end{itemize}

\end{proof}

\begin{Lemma} \label{temps-excursion}
Let $l>0$ and $E$ be an excursion of length $l$. Let $T_E$ be the  set of all lengths of excursions of $E$. Fix $0 < a < \frac{l}{2}$. Then $T_E \cap [a,2a( \neq \emptyset$.
\end{Lemma}

\begin{proof}
$T_E \cap [2a,l]$ is closed (by the lemma \ref{exc-ferme}) and non-empty (because $l \in T_E\cap [2a,l]$). Denote $l'=\min T_E\cap [2a,l]$. By the lemma \ref{exc-plus-petite}, there exists $l'' \in T_E$ such that $\frac{l'}{2} \leq l'' < l'$. Then $l'' <2a$ because $l'' < l'=\min T_E \cap [2a,l]$ and $l'' \geq \frac{l'}{2}\geq \frac{2a}{2}=a$. Therefore, $l'' \in T_E\cap[a,2a($.
\end{proof}

\subsection{Excursions of the orbit map} ~\\
\label{excursion-orbit-map}
Let $\gamma$ be a primitive element of $\F_2$. Recall that $L_\gamma$ is the geodesic of the Cayley graph of $\F_2$ generated by $\gamma$. We want to study the following map :
$E_\gamma : L_\gamma \longrightarrow \R_+$ such that $E_\gamma(u)=d(\tau_\rho(u),\mathrm{Axis}(\rho(\gamma)))$. 

\begin{Lemma} \label{E_gamma}
$E_\gamma$ is Lipschitz-continuous (hence continuous) and $\gamma$-invariant. 
\end{Lemma}

\begin{proof} It is a general fact that the distance map to any subspace of a metric space is 1-Lipschitz-continuous, because of the triangle inequality. Since the orbit map $\tau_{\rho | L_\gamma}$ is Lipschitz-continuous, we deduce the Lipschitz-continuity of $E_\gamma$. \\

The $\gamma$-invariance of $E_\gamma$ follows from the $\gamma$-invariance of $\mathrm{Axis}(\rho(\gamma))$. 
\begin{align*}
    E_\gamma(\gamma u) & =d(\tau_\rho(\gamma u),\mathrm{Axis}(\rho(\gamma)))=d(\rho(\gamma)\tau_\rho(u),\mathrm{Axis}(\rho(\gamma))) \\
     & =d(\tau_\rho(u),\mathrm{Axis}(\rho(\gamma))) \text{ because } \mathrm{Axis}(\rho(\gamma)) \text{ is } \rho(\gamma)\text{-invariant} \\
    &=E_\gamma(u)
\end{align*} 
\end{proof}

Since $L_\gamma$ is a geodesic in the Cayley graph of $\F_2$, it is isometric to $\R$, therefore we can think of $E_\gamma$ as a map from $\R$ to $\R$. Thus, we can apply the language of excursions defined previously. 

\begin{Definition} Let $\gamma$ be a primitive element in $\mathbb{F}_2$. Let $[u,v] \subset L_\gamma$ be a segment of the geodesic~$L_\gamma$. We say that $[u,v]$ is an \emph{excursion} if the map $E_{\gamma |[u,v]}$ is an excursion. \\
Let $K \geq 0$. We say that $[u,v]$ is a \emph{$K$-excursion} if the map $E_{\gamma |[u,v]}$ is an excursion such that $E_\gamma(u)=K$. In this case, we call \emph{length of excursion} of $[u,v]$ the length of excursion of  $E_{\gamma|[u,v]}$, that is the non-negative real $d(u,v)$. \\
At last, we say that $\gamma$ has an excursion (respectively a $K$-excursion) if there exists $[u,v] \in L_\gamma$ such that $[u,v]$ is an excursion (respectively a $K$-excursion). 
\end{Definition}

\begin{Lemma} \label{excursions-grandes}
There exist two sequences of positive reals $(K_n)_{n\in \mathbb{N}}$ and $(l_n)_{n \in \mathbb{N}}$, such that $K_n \to \infty$, $l_n \to \infty$ and, up to subsequence,  for all $n \in \mathbb{N}$, $\gamma_n$ has a $K_n$-excursion of length~$l_n$. 
\end{Lemma}

\begin{proof}
Let
\begin{align*}
    K_{\max,n} = \underset{[1,\gamma_n]}\max \, E_{\gamma_n} \text{ et } K_{\min,n} = \underset{[1,\gamma_n]}\min \, E_{\gamma_n} 
\end{align*}
In particular we have $K_{\max,n} = \max E_{\gamma_n}$ and $K_{\min,n} = \min E_{\gamma_n}$ since $E_{\gamma_n}$ is $\gamma_n$-invariant (see lemma \ref{E_gamma}). The hypothesis \ref{hyp-absurde} on the sequence $(\gamma_n)_{n \in \mathbb{N}}$ means that $K_{\max,n} \to \infty$. \\

\textbf{Fact} : For all $K_{\min,n} \leq K \leq K_{\max,n}$, $\gamma_n$ has a $K$-excursion. 
\begin{proof}
Indeed, $\gamma_n$ has a $K_{\min,n}$-excursion (by $\gamma_n$ invariance of $E_{\gamma_n}$) so by Lemma \ref{K-sous-exc}, $\gamma_n$ has a $K$-excursion.
\end{proof}

\begin{itemize}
    \item If $(K_{\min,n})_{n \in \mathbb{N}}$ is not bounded, then up to subsequence, we can assume that $K_{\min,n} \to \infty$. By definition of $K_{\min,n}$, there exists $u_n \in [1,\gamma_n]$ such that $E_{\gamma_n}(u_n)=K_{\min,n}$, and so $[u_n,\gamma_n u_n]$ is a  $K_{\min,n}$-excursion of length $|\gamma_n|$. By setting $K_n=K_{\min,n}$ and $l_n=|\gamma_n|$, we then have $K_n \to \infty, l_n \to \infty$ and $\gamma_n$ has a $K_n$-excursion of length $l_n$. \\
    
    \item If $(K_{\min,n})_{n \in \mathbb{N}}$ is bounded, then there exists $K>0$ such that for all $n \in \mathbb{N}, K_{\min,n} \leq K$. Let $K_n=\frac{K_{\max,n}}{2}$. Then, $K_n \to \infty$ and moreover, for $n$ sufficiently large, $K_{\min,n} \leq K \leq K_n < K_{\max,n}$. So, by the above fact, $\gamma_n$ has $K_n$-excursions. Now let us justify this excursion can be chosen in such a way that its length $l_n$ satisfies $l_n \to \infty$. \\
    Let $n \in \mathbb{N}$, there exists $u_n \in [1,\gamma_n]$ such that $E_{\gamma_n}(u_n)=K_{\max,n} > K_n$ ($E_{\gamma_n}$ is continuous), so there exists a $K_n$-excursion containing $u_n$. Denote it by $[v_n,w_n]$ and set $l_n=d(v_n,w_n)$. Now we are going to justify that $l_n \to \infty$. \\
    Let $x_n=\tau_\rho(u_n)$, we have $d(x_n,\mathrm{Axis}(\rho(\gamma_n)))=K_{\max,n}=2K_n$. \\
    Define $\partial N_{K_n}(\mathrm{Axis}(\rho(\gamma_n)))= \{y \in X : d(y,\mathrm{Axis}(\rho(\gamma_n)))=K_n \}$ and let $y_n$ be a projection of $x_n$ on $\partial N_{K_n}(\mathrm{Axis}(\rho(\gamma_n)))$. Then $y_n$ satisfies : $y_n \in \partial N_{K_n}(\mathrm{Axis}(\rho(\gamma_n)))$ and $\forall y \in \partial N_{K_n}(\mathrm{Axis}(\rho(\gamma_n))), d(x_n,y_n) \leq d(x_n,y)$. Since the map $d(\cdot,\mathrm{Axis}(\rho(\gamma_n)))$ is 1-Lipschitz-continuous, we have : 
    \begin{align*}
        & |d(x_n,\mathrm{Axis}(\rho(\gamma_n)))-d(y_n,\mathrm{Axis}(\rho(\gamma_n)))| \leq d(x_n,y_n), \\
        \text{hence} \qquad & K_n \leq d(x_n,y_n).
    \end{align*}

    In addition, because $[v_n,w_n]$ is a $K_n$-excursion, $\tau_\rho(v_n)$ and $\tau_\rho(w_n)$ belong to $\partial N_{K_n}(\mathrm{Axis}(\rho(\gamma_n)))$ so $2d(x_n,y_n) \leq d(x_n,\tau_\rho(v_n))+d(x_n,\tau_\rho(w_n))$. Then :
    \begin{align*}
        2K_n & \leq 2d(x_n,y_n) \leq d(x_n,\tau_\rho(v_n))+d(x_n,\tau_\rho(w_n)) \\
        & \leq C'd(u_n,v_n) +C'd(u_n,w_n) \qquad
        \text{ because $\tau_\rho$ is $C'$-Lipschitz-continuous }\\
        & = C'd(v_n,w_n) \qquad \text{ because } u_n \in [v_n,w_n] \\
        & = C'l_n
    \end{align*}
    We conclude by using that $K_n \to \infty$.

\end{itemize}

\end{proof}

\subsection{Quasi-loops}
\label{quasi-loop}

\begin{Definition}
Let $\varepsilon > 0$ and $w \in \mathbb{F}_2$ (not necessarily primitive). We say that $w$ is an $\emph{$\varepsilon$-quasi-loop} $ if we have the following inequality :
$$ d(\rho(w)o,o) \leq \varepsilon ~ |w|$$
\end{Definition}

Thinking of $\varepsilon$ as very small, an $\varepsilon$-quasi-loop is an element that does not displace the points much. Note that the definition of a quasi-loop depends on the representation $\rho$. \\

Let $\gamma$ be a primitive element of $\mathbb{F}_2$ and $u \in L_\gamma$. We denote by $\lfloor u \rfloor$ the integer point in $L_\gamma$ just before $u$ (if $u$ is an integer point in $L_\gamma$, $\lfloor u \rfloor=u$) and $\lceil u \rceil$ the integer point of $L_\gamma$ just after $u$ (thus $\lfloor u \rfloor$ and $\lceil u \rceil$ are the endpoints of an edge of length $1$ in the Cayley graph and $u$ belongs to this edge). \\ 

Recall that we have fixed a (Bowditch) representation $\rho : \F_2 \to \mathrm{Isom}(X)$ and that the notions of a $K$-excursion and of a $\varepsilon$-quasi-loop depend on $\rho$.

\begin{Lemma} \label{excursion->QB}
Let $\varepsilon > 0$. There exist $l_\varepsilon > 0$ and $K_\varepsilon >0$ such that for all primitive elements $\gamma$, for all $K \geq K_\varepsilon, l \geq l_\varepsilon$, if $[u,v]$ is a $K$-excursion of length $l$, then the element  $w=\lfloor u \rfloor^{-1}\lfloor v \rfloor$ (which is a subword of $\gamma$) is an $\varepsilon$-quasi-loop. 
\end{Lemma}

\begin{proof}
Let $\varepsilon'=\frac{\varepsilon}{2}$. \\
Let $\gamma$ be a primitive element in $\F_2$ and $[u,v] \subset L_\gamma$ such that $[u,v]$ is a $K$-excursion of length $l$. Then $d(\tau_\rho(u),\mathrm{Axis}(\rho(\gamma)))=d(\tau_\rho(v),\mathrm{Axis}(\rho(\gamma)))=K$, for all $t \in [u,v],$ \, $d(\tau_\rho(t),\mathrm{Axis}(\rho(\gamma)))\geq K$ and $d(u,v) = l$. Choose $\ell_\gamma$ a geodesic in $X$ with endpoints $\rho(\gamma)^+$ and $\rho(\gamma)^-$, where $\rho(\gamma)^+$ and $\rho(\gamma)^-$ are respectively the attracting and repelling fixpoints of the hyperbolic isometry $\rho(\gamma)$. Using Lemma \ref{neighborhood}, we conclude that $d(\tau_\rho(u),\ell_\gamma) \leq K + C(\delta)$ and $d(\tau_\rho(v),\ell_\gamma) \leq K + C(\delta)$, where $C(\delta)$ is the constant introduced in the lemma. In addition, since for all $t \in [u,v],$ \,
$d(\tau_\rho(t),\mathrm{Axis}(\rho(\gamma)))\geq K$, we can also conclude from lemma \ref{neighborhood} that for all $t \in [u,v]$, \, $d(\tau_\rho(t),\ell_\gamma)\geq K$. \\

First, let us show that for $K$ and $l$ large enough, $d(\tau_\rho(u),\tau_\rho(v))\leq \varepsilon' d(u,v)=\varepsilon' l$ : \\

We denote by $d=d(\tau_\rho(u),\tau_\rho(v))$ and $L=\mathrm{length}(\tau_\rho([u,v]))$. Then, since $\tau_\rho$ is piecewise geodesic on $[u,v]$, we can apply the proposition  \ref{long_ext_banane} to obtain, after denoting $D(\delta)=\max(C(\delta),\delta)$ :
\begin{enumerate}
    \item If $d\leq2K+6\delta$, then $L \geq (2^{\frac{d}{2\delta}-\frac{D(\delta)}{\delta}-5}-2)\delta$
    \item If $d>2K+6\delta$, then there exists an integer $n \geq 2$ such that :
    \begin{equation*}
    \left\{ 
    \begin{array}{ll} 
        L  \geq (n-1) (2^{\frac{K}{\delta}-3}-2)\delta \\ 
        d  \leq 18n\delta +2K +2C(\delta)
    \end{array}
    \right. 
    \end{equation*}
\end{enumerate}
On the other hand, we also have that $L \leq C'd(u,v)=C'l$. Indeed, this is a consequence of the $C'$-Lipschitz-continuity $\tau_\rho$. Therefore :
\begin{enumerate}
    \item If $d \leq 2K+6\delta$, we have : 
    \begin{align*}
         (2^{\frac{d}{2\delta}-\frac{D(\delta)}{\delta}-5}-2)\delta & \leq L \leq C'l \\
        \text{ so      } 2^{\frac{d}{2\delta}-\frac{D(\delta)}{\delta}-5} & \leq \frac{C'}{\delta}l+2 \\
        \text{ then   } \qquad d & \leq 2\delta\log_2 \left( \frac{C'}{\delta}l+2 \right)+2D(\delta)+10\delta \\
        \text{ But } \qquad \hspace{0.2cm} &  \frac{2\delta\log_2(\frac{C'}{\delta}l+2)+2D(\delta)+10\delta}{l} \underset{l \to +\infty}\longrightarrow 0, \\
    \end{align*}
    so there exists $l_\varepsilon >0$ (depending only on $C', \delta$ and $\varepsilon$) such that :
    $$ \text{ If } l \geq l_\varepsilon \text{ then } d \leq \varepsilon'l.$$ 
    
    \item If $d > 2K+6\delta$, we have : 
    \begin{align*}
       & (n-1)(2^{\frac{K}{\delta}-3}-2)\delta \leq L \leq C'l 
        \text{ so } l \geq (n-1)(2^{\frac{K}{\delta}-3}-2)\frac{\delta}{C'}. \\
        \text{ On the other hand } &d  \leq 18n\delta +2K +2C(\delta), \\
        \text{ then } & \frac{d}{l} \leq \frac{18n\delta+2K+2C(\delta)}{(n-1)(2^{\frac{K}{\delta}-3}-2)\frac{\delta}{C'}}=\frac{18\delta + \frac{2K+2C(\delta)}{n}}{(1-\frac{1}{n})(2^{\frac{K}{\delta}-3}-2)\frac{\delta}{C'}}. \\
        \text{ But } &  n \geq 2 \text{ so } 1-\frac{1}{n} \geq \frac{1}{2} \text{ and } \frac{2K+2C(\delta)}{n} \leq K+C(\delta), \\
        \text{therefore } & \frac{d}{l} \leq \frac{2C'}{\delta}\frac{18\delta+K+C(\delta)}{2^{\frac{K}{\delta}-3}-2} \underset{K \to +\infty} \longrightarrow 0.
    \end{align*}
    so there exists $K_\varepsilon >0$ (depending only on $C', \delta$ and $\varepsilon$) such that if $K \geq K_\varepsilon$, then \, $d\leq \varepsilon'l.$
\end{enumerate}
 Thus we have shown that if $K \geq K_\varepsilon$ and $l \geq l_\varepsilon$, we have in every case $d\leq \varepsilon'l$. \\

Now, let us show that this implies that $w=\lfloor u \rfloor^{-1}\lfloor v \rfloor$ is an $\varepsilon$-quasi-loop. 

    \begin{align*}
        d(\rho(w)o,o) & = d(\rho(\lfloor u \rfloor^{-1}\lfloor v \rfloor)o,o) = d(\rho(\lfloor u \rfloor)o,\rho(\lfloor v \rfloor)o) = d(\tau_\rho(\lfloor u \rfloor),\tau_\rho(\lfloor v \rfloor))\\
        & \leq d(\tau_\rho(\lfloor u \rfloor),\tau_\rho(u)) + d(\tau_\rho(u),\tau_\rho(v)) + d(\tau_\rho(v),\tau_\rho( \lceil v \rceil) \text{ by the triangle inequality}\\
        & \leq 2C'+d \text{ because $u$ and $\lfloor u \rfloor$, resp. $v$ and $\lceil v \rceil $, are at a distance less than 1 in the Cayley graph } \\
        & \leq 2C'+\varepsilon'd(u,v) \text{ by what have been previously done } \\
          & \leq 2C'+ \eps' (d(\lf u \rf,\lf v \rf)+d(\lf v \rf,v)-d(\lf u \rf,u)) \text{ because } \lf u \rf,u,\lf v \rf, v \text{ are aligned in this order on } L_\gamma \\
        & \leq 2C'+\eps'+\eps'|w|.
    \end{align*}
    Let us further assume that $l \geq \frac{2C'}{\eps'}+2$. Then in particular $|w|=d(\lf u \rf,\lf v \rf)=d(\lf u \rf, u)+d(u,v)-d(\lf v \rf, v) \geq l-1 \geq \frac{2C'}{\eps'}+1$. Therefore :  
    
    \begin{equation*}
    d(\rho(w)o,o) \leq 2C'+\eps'+\eps'|w| \leq 2\eps '|w|=\eps |w|.
    \end{equation*}
    Thus, after possibly changing $l_\varepsilon$ to $\max(l_\varepsilon,\frac{4C'}{\varepsilon}+2)$, we have shown that if $l \geq l_\varepsilon$ and $K \geq K_\varepsilon$, then $w$ is an $\varepsilon$-quasi-loop.

\end{proof}

\subsection{Induction step} ~\\ 
Let $n \in \mathbb{N}$ and $0 \leq i \leq r(n)$. Recall that the notations $w_i(\gamma_n)$ and $w'_i(\gamma_n)$ have been defined in definition \ref{def-wi}, and refer to some specific subwords of (a cyclic permutation of) $\gamma_n$, corresponding to a truncation of the continued fraction expansion of the slope of $\gamma_n$. The integers $l_i(\gamma_n)=|w_i(\gamma_n)|$ and $l'_i(\gamma_n)=|w'_i(\gamma_n)|$ refer to their lengths and $r(n)$ is the depth of the continued fraction expansion of the slope $\gamma_n$. In order to reduce the amount of notations, we write more simply $w_i(n),w'_i(n)$ and $l_i(n),l'_i(n)$. \\

Let $\cycl{w_i(n)}$ be a cyclic permutation of $w_i(n)$ and $\cycl{w'_i(n)}$ a cyclic permutation of $w'_i(n)$ adapted to $\cycl{w_i(n)} $ (see Lemma \ref{perm-cycl}). Recall that, by Lemma \ref{perm-cycl}, there exists a cyclic permutation of $\gamma_n$ that can be written on the alphabet $\{\cycl{w_i(n)},\cycl{w'_i(n)}\}$. Thus, a subword $u$ of $\gamma_n$ can be written in the following way :
\begin{equation*}
    u=p w^1 \hdots w^{b} s
\end{equation*}
with :
\begin{itemize}
    \item $b \in \mathbb{N}$
    \item $p$ is a suffix of either $\cycl{w_i(n)}$ or $\cycl{w'_i(n)}$
    \item $s$ is a prefix of either $\cycl{w_i(n)}$ or $\cycl{w'_i(n)}$
    \item For all $1 \leq k \leq b$, $w^k \in \{\cycl{w_i(n)},\cycl{w'_i(n)}\}$ \\
\end{itemize}

Recall that the constants $C$ and $C'$ have been chosen at the beginning of the proof of Proposition \ref{uniform-neighborhood} and satisfy~:
$$ \forall \gamma \in \mathcal{P}(\mathbb{F}_2), \hspace{0,5cm} \frac{1}{C} \Vert \gamma \Vert  \leq l(\rho(\gamma)), \qquad \forall u \in \mathbb{F}_2, \hspace{0,5cm} d(\rho(u)o,o) \leq C' |u| \qquad \text{and} \qquad CC'\geq 1. $$ 
In the following lemma, we find an $\varepsilon$-quasi-loop in a cyclic permutation of $w_i(n)$ which occupy at least half of its length. We also ask that the remainder of the cyclic permutation of $w_i(n)$ which is not in the quasi-loop is sufficiently large (to be able later on to continue the process of finding quasi-loop inside) and that the length of $w_i(n)$ is not too big (to be able to control the number of $w_i(n)$ we can find). 

\begin{Lemma} \label{trouve-QB}
Let $\displaystyle 0 < \varepsilon < \frac{1}{C}$ and $\alpha >6$. Fix $\displaystyle r_0= 4+\frac{2\varepsilon}{C'}$. Let $r \geq r_0$. \\
There exists a constant $R >0$ and two integers $n_0 \in \N, i \in \N$ such that, for all integer $n \geq n_0$ the following properties are satisfied~: 

\begin{enumerate}
    \item $1 \leq i \leq r(n)$
    \item $l_{i-1}(n) \geq r$
    \item $\displaystyle l_i(n) \leq \frac{R}{\alpha}$
    \item There exists a cyclic permutation of $w_i(n)$, denoted by $\cycl{w_i(n)}$ such that $\cycl{w_i(n)}=v_1v_2$, with $v_1$ and $v_2$ two elements of $\F_2$ satisfying the following properties :
    \begin{enumerate}
        \item $v_1$ is an $\varepsilon$-quasi-loop
        \item $|v_1| \geq |v_2|$
        \item $|v_2| \geq r$
    \end{enumerate}
\end{enumerate}
\end{Lemma}

\begin{proof}
First of all, let's consider two sequences $(K_n)_{n \in \N}$ and $(l_n)_{n \in \N}$ as in the lemma \ref{excursions-grandes}. Then $K_n \to \infty$, $l_n \to \infty$ and for all $n\in \N$, $\gamma_n$ has a  $K_n$-excursion of length $l_n$. Now, let's introduce the constants $K_{\frac{\eps}{2}}$ and $l_{\frac{\eps}{2}}$ given by the lemma \ref{excursion->QB} for $\frac{\eps}{2}$. In order to simplify the notation, we still denote by $K_\eps$ and $l_\eps$ these two constants. Then, there exists an integer $n_1 \in \N$ such that for all integers $n \geq n_1$, we have $l_n \geq l_\varepsilon$ and $K_n \geq K_\varepsilon$. Now, let $\displaystyle r'= \frac{C'r}{\frac{1}{C}-\varepsilon}$. Then, because $CC' \geq 1$, we have $\displaystyle r'=\frac{C'r}{\frac{1}{C}-\varepsilon} \geq CC'r \geq r$. Let $l(r',\varepsilon)=\max \, (r'+1,2l_\varepsilon)$. Since $r(n) \to \infty$ by the lemma \ref{forme-frac-continue}, we deduce the existence of an integer $n_2 \in \N$ such that for all integers $n \geq n_2$, we have $r(n) \geq l(r',\varepsilon)$. Let us consider $i$ the smallest integer such that $i \geq l(r',\varepsilon)$. Then for all $n \geq n_2$, we have $l(r',\varepsilon) \leq i \leq r(n)$ (because $r(n)$ is an integer). Since $l_n \to \infty$, we can find $n_3 \in \N$ such that for all $n \geq n_3$, we have $L_i < l_n$ (recall that $L_i$ is a constant introduced in \ref{bornes-li} which satisfies $l_i(n) \leq L_i$ for all $n \in \N$). Then, set $n_0=\max(n_1,n_2,n_3)$ and let's summarise the inequalities that are true for all integers $n \geq n_3$ : 
$$l_n \geq l_\varepsilon, \, K_n \geq K_\varepsilon, \, r(n) \geq r'+1, \, r(n) \geq 2l_\varepsilon, \, l_n > L_i \text{ and  } l(r',\varepsilon) \leq i \leq r(n). $$ 
Finally, we set $R=\alpha L_i$. \\

Now that all these constants have been introduced, we show that this choice of $n_0, i$ and $R$ satisfies the property requested in the lemma. Let $n \geq n_0$. We can easily check the first three properties : 
\begin{enumerate}
    \item $\begin{aligned}[t]
        r(n) &  \geq i \geq l(r',\varepsilon) \geq r'+1 \text{ because } l(r',\varepsilon)=\max \, (r',2l_\varepsilon-1) \\
        & \geq r+1 \geq r_0 \text{ because $r$ is chosen larger than $r_0$ } \\
        & > 4 \text{ because } r_0=4+\frac{2\varepsilon}{C'}, \\ 
     & \text{hence we have $1 \leq i \leq r(n)$.}
    \end{aligned}$\hfill \mbox{}
    \item $\begin{aligned}[t]
        l_{i-1}(n) & \geq i-1 \text{ by the inequalities of \ref{bornes-li}} \\
        & \geq r \text{ as already seen above }
    \end{aligned}$ \hfill \mbox{}
    \item $\begin{aligned}[t]
        l_i(n) &  \leq L_i \text{ by the inequalities of  \ref{bornes-li}} \\
         & = \frac{R}{\alpha} \text{ by definition of $R$.} 
    \end{aligned}$
\end{enumerate}
Now let us show the fourth property. \\

The element $\gamma_n$ has a $K_n$-excursion of length $l_n$. But $l_i(n) \leq L_i < l_n$, so we can use lemma \ref{temps-excursion} with $a=\frac{l_i(n)}{2}$ to show the existence of a $K'_n$-sub-excursion of length $l'_n \in [a,2a($. Denote it by $[u',v']$. Then we have $d(u',v')=l'_n, \, K'_n \geq K_n $ and $\displaystyle \frac{l_i(n)}{2} \leq l'_n < l_i(n)$. \\
Then, since $K_n \geq K_\varepsilon$, we deduce $K'_n \geq K_\varepsilon$. In addition :
\begin{align*}
    l'_n & \geq \frac{l_i(n)}{2} \geq \frac{i}{2} \text{ by the inequalities of \ref{bornes-li}} \\
    & \geq \frac{l(r',\varepsilon)}{2} \text{ because } i\geq l(r',\varepsilon) \\
    & \geq \frac{2l_\varepsilon}{2} \text{ by definition of $l(r',\varepsilon)$} \\
    & =l_\varepsilon.
\end{align*}
We apply Lemma \ref{excursion->QB} to $[u',v']$ in order to show that the subword $v_1:=\lf u' \rf^{-1} \lf v' \rf$ of (a cyclic permutation of) $\gamma_n$ is an $\frac{\eps}{2}$-quasi-loop. 
Let us now look at the length of $v_1$. 
\begin{align*}
   \text{We have } \qquad & |v_1|=d(\lf u' \rf, \lf v' \rf)=d(\lf u' \rf,u')+d(u',v')-d(\lf v' \rf,v') \\
    \text{ Hence } \qquad & l'_n-1<|v_1| < l'_n +1, \\
    \text{ so } \qquad & \frac{l_i(n)}{2}-1<|v_1|<l_i(n)+1. \\
    \text{ Then } \qquad & \frac{l_i(n)}{2}-\frac{1}{2} \leq |v_1| \leq l_i(n)
\end{align*}
Note that, after possibly deleting the last letter of $v_1$ or adding a letter at the end of $v_1$, we can in fact assume that the resulting word, which we denote by $v_1$ again, satisfies $\frac{l_i(n)}{2}\leq v_1 \leq l_i(n)-1$. Indeed : 
\begin{itemize}
    \item If $|v_1|=l_i(n)$, write $v_1=v'_1s$, with $|s|=1$. Then $|v'_1|=l_i(n)-1$, and we also have $|v'_1| \geq \frac{l_i(n)}{2}$, because $l_i(n) \geq 2$ (since $i\geq1$).
    \item If $|v_1|<\frac{l_i(n)}{2}$, write $v'_1=v_1s$, with $|s|=1$. Then $|v'_1|=|v_1| +1 < \frac{l_i(n)}{2}+1 \leq l_i(n)$, so $|v_1|+1 \leq l_i(n)-1$, and $\frac{l_i(n)}{2} \leq \frac{l_i(n)}{2}+\frac{1}{2}<|v_1|+1$. 
\end{itemize}

\begin{Lemma} \label{QB-last-letter}
    Let $\eps>0$. Let $w \in \F_2$ be an $\eps$-quasi-loop of length $|w| \geq \frac{C'}{\eps}$. Then, after deleting the last letter of $w$ or adding a letter at the end of $w$, the resulting word is an $2\eps$-quasi-loop. 
\end{Lemma}

\begin{proof}
    \begin{itemize}
        \item If $w'=ws$, with $|s|=1,|w'|=|w|+1$ :
        \begin{align*}
            d(\rho(w')o,o)=d(\rho(ws)o,o) & \leq d(\rho(ws)o,\rho(w)o)+d(\rho(w)o,o) 
            \\
            & \leq d(\rho(s)o,o)+\eps|w| \qquad \text{ because $w$ is an $\eps$-quasi-loop} \\
            & \leq C'+\eps|w| \qquad \text{ since } |s|=1 \\ 
            & \leq 2\eps|w| \qquad \text{ because } C' \leq \eps|w|
        \end{align*}
        \item If $w=w's$, with $|s|=1, |w'|=|w|-1$ :
        \begin{align*}
            d(\rho(w')o,o) & \leq d(\rho(w')o,\rho(w's)o)+d(\rho(w's)o,o) \\
            & \leq d(\rho(s)o,o)+d(\rho(w)o,o) \\
            & \leq C'+\eps|w| \qquad \text{ because $w$ is an $\eps$-quasi-loop and $|s|=1$} \\
            & \leq 2\eps |w| \qquad \text{ because } C' \leq \eps|w|
        \end{align*}
    \end{itemize}
\end{proof}
Thus by lemma \ref{QB-last-letter}, $v_1$ is an $\eps$-quasi-loop. 

Since $|v_1| \leq l_i(n)-1 \leq |\gamma_n|$, there exists a subword $v_3$ of (a cyclic permutation of) $\gamma_n$ of length $l_i(n)$ that can be written $v_3=v_1v_2$, with $1 \leq |v_2| \leq |v_3|=l_i(n)$. Since $v_3$ is of length $l_i(n)$, we can again use the Proposition \ref{magic-len} to ensure that after possibly changing the last letter of $v_3$, that is the last letter of $v_2$ (because $v_2$ is non empty), $v_3$ is in fact a cyclic permutation of $w_i(n)$. Then, noting again $v_3$ and $v_2$ after this potential change of letter, there exists a cyclic permutation of $w_i(n)$, which we denote $\cycl{w_i(n)}$, such that $\cycl{w_i(n)}=v_1v_2$. Recall that we have already shown that $v_1$ is an $\varepsilon$-quasi-loop, and because $|v_1| \geq \frac{l_i(n)}{2}$, we have $|v_1|\geq |v_2|$. So we still have to show that $|v_2| \geq r$ to finish the proof of the fourth point. We proceed as follows :  
\begin{align*}
    \frac{1}{C}l_i(n) & \leq d(\rho(\cycl{w_i(n)})o,o) \qquad \text{ by the Bowditch hypothesis, because $\cycl{w_i(n)}$ is primitive } \\ 
    & \leq d(\rho(v_1v_2)o,o) \qquad \text{ since } \cycl{w_i(n)}=v_1v_2 \\
    & \leq d(\rho(v_1v_2)o,\rho(v_1)o) + d(\rho(v_1)o,o) \qquad \text{ by the triangle inequality } \\
    & = d(\rho(v_2)o,o)+d(\rho(v_1)o,o) \qquad \text{ because $\rho(v_1)$ is an isometry} \\
    & \leq C'|v_2| + \varepsilon|v_1| \qquad \text{ since $v_1$ is an $\varepsilon$-quasi-loop} \\
    & \leq C'|v_2| + \varepsilon l_i(n) \qquad \text{ because } |v_1|\leq l_i(n). \\
    \text{ Therefore } |v_2| & \geq \frac{1}{C'}(\frac{1}{C}-\varepsilon)l_i(n) \geq \frac{1}{C'}(\frac{1}{C}-\varepsilon) r' \qquad \text{ because } l_i(n) \geq r', \\
    & \geq \frac{1}{C'}(\frac{1}{C}-\varepsilon) \frac{C'}{\frac{1}{C}-\varepsilon}r=r, \qquad \text{which finishes the proof that }|v_2| \geq r
\end{align*}
\end{proof}

The following lemma aims, when given a sufficiently large subword of some $\gamma_n$,
to write it as a concatenation of subwords being either $\varepsilon$-quasi-loops or
sufficiently large "remainders", and such that the proportion of the word in an $\varepsilon$-quasi-loop is at least $c$, where $c$ is a constant between 0 and $\frac{1}{4}$, fixed in advance. It will be used recursively in the next lemma.

\begin{Lemma} \label{decoupe-restes-c}
Let $\displaystyle 0<\varepsilon < \frac{1}{C}$ and $\displaystyle 0<c<\frac{1}{4}$. Fix $\displaystyle r_0=4+\frac{2\varepsilon}{C'}$ and let $r \geq r_0$. There exists a constant $R>0$ and an integer $n_0 \in \N$, such that, given any integer $n \geq n_0$ and any subword $u$ of $\gamma_n$ such that $|u| \geq R$, then there exists a positive integer $q \in \N^*$, a subset $QL \subset \{1,\hdots,q\}$ and $q$ words $u_1, \hdots,u_q \in \F_2$ such that :
\begin{enumerate}
    \item $u=u_1\hdots u_q$
    \item For all $k \in QL$, $u_k$ is an $\varepsilon$-quasi-loop
    \item $\displaystyle \sum_{k \notin QL} |u_k| \leq (1-c)|u|$
    \item For all $k \notin QL$, $|u_k| \geq  r$
\end{enumerate}
\end{Lemma}

\begin{proof}
Let $b=\frac{8c+2}{1-4c}$. Then $c=\frac{b-2}{4b+8}$ and for all $b' \geq b$, we have $\frac{b'-2}{4b'+8}\geq c$. In addition, since $0 < c < \frac{1}{4}$, we have $b>2$. Let $\alpha=2b+4$, we have $\alpha>8$. \\
Now let us introduce the constants $R>0$, $i$ and $n_0$ given by the lemma \ref{trouve-QB}. Let $n \geq n_0$ be an integer and $u$ a subword of $\gamma_n$ such that $|u| \geq R$. Then, the lemma \ref{trouve-QB} states that $1 \leq i \leq r(n)$, \, $l_{i-1}(n) \geq r$, \, $l_i(n) \leq \frac{|u|}{\alpha}$, \, and there exists a cyclic permutation of $w_i(n)$, denoted by $\cycl{w_i(n)}$, that decomposes into the form $\cycl{w_i(n)}=v_1v_2$, with $v_1$ an $\varepsilon$-quasi-loop and $|v_1| \geq |v_2| \geq r$. Let $\cycl{w'_i(n)}$ be a cyclic permutation of $w'_i(n)$ adapted to $\cycl{w_i(n)}$ (see the lemma \ref{perm-cycl}). \\

Then we can write a decomposition of $u$ under the form : 
$u=p w_1 \hdots w_{b'} s $, 
 with $b \in \N$, $p$ a suffix of $\cycl{w_i(n)}$ or $\cycl{w'_i(n)}$, $s$ a prefix of $\cycl{w_i(n)}$ or $\cycl{w'_i(n)}$, and for all $1 \leq k \leq b'$, $w_k \in \{\cycl{w_i(n)},\cycl{w'_i(n)}\}$. \\

The lemma \ref{trouve-QB} ensures that $l_i(n) \leq \frac{|u|}{\alpha}$. Therefore we can use the lemma \ref{bloc} to conclude that $b' \geq \frac{\alpha-4}{2}=b >2$. Namely there is at least three central blocs in the decomposition ($b'$ is an integer). \\
Denote $p'=pw_1$ and $s'=w_{b'}s$. We have :  
\begin{align*}
    |p'| & \geq |w_1| \text{ because } |p'|=|p| + |w_1| \\
            & \geq l_i(n) \text{ since } l'_i(n) \geq l_i(n) \\
            & \geq l_{i-1}(n) \text{ because the sequence $(l_i(n))_i$ is increasing} \\
            & \geq r \text{ as provided by the lemma \ref{trouve-QB}}.
\end{align*}

We also obtain $|s'| \geq l_i(n) \geq r $.  \\
Therefore the word $u$ can be written : $u = p' w_2 \hdots w_{b'-1}s'$. \\

Moreover, by the lemma \ref{perm-cycl}, $\cycl{w_{i(n)}(n)}$ is either a prefix or a suffix of $\cycl{w'_{i(n)}(n)}$ so there exists a word $w$ such that $\cycl{w'_i(n)}=w\cycl{w_i(n)}$ or $\cycl{w'_i(n)}=\cycl{w_i(n)}w$. In addition, $\cycl{w_{i(n)}(n)}=v_1v_2$, so $\cycl{w'_i(n)}=wv_1v_2$ or $\cycl{w'_i(n)}=v_1v_2w$. Then, for all $k \in \{2,\hdots,b'-1\}$, $w_k \in \{\cycl{w_i(n)},\cycl{w'_i(n)} \}$ so $w_k$ is a concatenation of $w,v_1$ and $v_2$. In addition, $v_1$ is an $\varepsilon$-quasi-loop (as provided by the lemma \ref{trouve-QB}), and we have :
\begin{align*}
    |w| & = l'_i(n)-l_i(n) \text{ because } |\cycl{w'_i(n)}|=|w|+|\cycl{w_i(n)}| \\
        & =l_{i-1}(n) \text{ because } l'_i(n)=l_i(n)+l_{i-1}(n) \\
        & \geq r \text{ by the lemma \ref{trouve-QB} }, \\
   \text{ on the other hand, } |v_2| & \geq r \text{ still by lemma \ref{trouve-QB}}.
\end{align*}
Since we have previously shown that $|p'| \geq r, |s'|\geq r$, we can indeed write a decomposition of $u$ into the form $u=u_1\hdots u_q$ (with $u_1=p'$ and $u_q=s'$) such that there exists a subset $QL \subset \{1,\hdots,q\}$, such that for all $k \in QL, u_k$ is an $\varepsilon$-quasi-loop and for all $k \notin QL$, $|u_k| \geq r$. Moreover, $\# QL = b'-2$ since each bloc $w_2, \hdots, w_{b'-1}$ contains the $\varepsilon$-quasi-loop $v_1$ exactly once. It remains to show that $\sum_{k \notin QL} |u_k| \leq (1-c)|u|$. In order to do so, let us find a lower bound on the total length of the $\varepsilon$-quasi-loops :

\begin{align*}
        \underset{k \in QL} \sum |u_k| &  = (b'-2)|v_1| \qquad \text{ since the quasi-loop $v_1$ appears exactly $b'-2$ times in our decomposition, } \\
                & \geq (b'-2)\frac{l_i(n)}{2} \qquad \text{ because $|v_1| \geq |v_2|$ by the lemma \ref{trouve-QB}} 
\end{align*}

\begin{align*}
    \text{ But } |u| & = |p|+ \sum_{k=1}^{b'} |w_k| + |s| \qquad \text{ because } u=p w_1 \hdots w_{b'} s \\
                    & \leq |p| + b' \max \{l'_i(n),l_i(n) \} + |s| \qquad  \text{ because } w_k \in \{\cycl{w_i(n)},\cycl{w'_i(n)}\} \\ 
                    & \leq (b'+2)\max \{l'_i(n),l_i(n) \}  \qquad \text{ since $p$ (resp. $s$) is a suffix (resp. prefix) of $\cycl{w_i(n)}$ or $\cycl{w'_i(n)}$} \\ 
                    & \leq (b'+2)l'_i(n) \qquad \text{ because } l_i(n) \leq l'_i(n) \\
                    & \leq 2(b'+2)l_i(n) \qquad \text{because } l'_i(n) \leq 2l_i(n) 
\end{align*}
\begin{align*}
    \text{ Therefore } \underset{k \in QL} \sum |u_k| & \geq \frac{b'-2}{2}l_i(n) \geq \frac{b'-2}{4(b'+2)}|u| \\
    & \geq \frac{b-2}{4(b+2)}|u| \qquad \text{ because we have shown at the beginning of the proof that } b' \geq b \\
    & = c|u|  \qquad \text{ by definition of $b$.} 
\end{align*}        
The last inequality can we rewritten as follows : $\displaystyle \sum_{k \notin QL} |u_k| \leq (1-c)|u|$, which completes the proof. 
\end{proof}

\subsection{Final contradiction and conclusion} ~\\ 

Now, we are able to find a primitive element $\gamma$ (from the sequence $(\gamma_n)_{n\in \mathbb{N}}$) which contains a very large proportion of quasi-loops.

\begin{Lemma} \label{trouve-gamma}
Let $0 < \varepsilon < \frac{1}{C}$ and $1 - \frac{1}{C'}(\frac{1}{C}-\varepsilon) <  \lambda < 1$. There exists a primitive element $\gamma$ such that $\gamma$ contains  $\varepsilon$-quasi-loops that occupy at least a proportion $\lambda$ of $\gamma$.
\end{Lemma}

\begin{proof} 
Let $0 < c < \frac{1}{4}$ and $r_0=4+\frac{2\varepsilon}{C'}$. Precisely, we will show the following property by recursion on $k \in \N$ : \\

For any integer $k \in \N$, for any real $r \geq r_0$, there exists an integer $n_0 \in \N$, such that for all $n \geq n_0$, there exists an integer $q \in \N^*$, a subset $QL \subset \{1,\hdots,q\}$, and some elements $u_1, \hdots, u_q \in \F_2$ satisfying the following properties : 
\begin{enumerate}
    \item $\gamma_n = u_1 \hdots u_q$
    \item For all $i \in QL$, $u_i$ is an $\varepsilon$-quasi-loop
    \item $\displaystyle \sum_{i \notin QL} |u_i| \leq (1-c)^k|\gamma_n|$
    \item For $i \notin QL$, $|u_i| \geq r$ \\
\end{enumerate}

\begin{itemize}
    \item For $k=0$, it's trivial, it is sufficient to choose $n$ large enough so that $|\gamma_n| \geq r$, $q=1, QL=\emptyset$ and thus the properties are satisfied.
    \item Suppose that this is true for some $k$. Let $r \geq r_0$. Let us introduce the constants $R>0$ and $n_0 \in \N$ given by Lemma \ref{decoupe-restes-c}. Now, let us apply the recursion hypothesis to $r_1 = \max (R, r_0)$. Then, there exists an integer $n_1 \in \N$ such that for all $n \geq n_1$, there exists $q \in \N^*, QL \subset \{1,\hdots,q\}$ and $u_1, \hdots, u_q$ such that $\gamma_n=u_1 \hdots u_q$, for all $i \in QL$, $u_i$ is an $\varepsilon$-quasi-loop, $\sum_{i \notin QL} |u_i| \leq (1-c)^k|\gamma_n|$ and for all $i \notin QL, |u_i| \geq r_1 \geq R$. This is still true for all integers $n \geq \max(n_0,n_1)$. Then, since for all $i \notin QL$, $u_i$ is a subword of $\gamma_n$, with $n \geq n_0$, and $|u_i| \geq R$, we can apply Lemma \ref{decoupe-restes-c} to each $u_i \notin QL$. That is, for all $i \notin QL$, there exists an integer $q_i \in \N^*$, a subset $QL_i \subset \{1,\hdots q_i \}$, and $q_i$ elements $u_{i,1}, \hdots, u_{i,q_i} \in \F_2$ such that : $u_i=u_{i,1}\hdots u_{i,q_i}$, for all $j \in QL_i$, $u_{i,j}$ is an $\varepsilon$-quasi-loop, $\sum_{j \notin QL_i} \leq (1-c)|u_i|$ and for all $j \notin QL_i, |u_{i,j}| \geq r$. Then we have : 
    \begin{enumerate}
        \item  \[ \gamma_n=\prod_{i=1}^q \left\{ \begin{array}{cc}
            u_i  & \text{ if } i \in QL \\
            u_{i,1}\hdots u_{i,q_i} & \text{ if } i \notin QL 
        \end{array} \right. \] (here the product denotes the concatenation)
        \item For all $i \in QL, u_i$ is an $\varepsilon$-quasi-loop and for all $i \notin QL$, for all $j \in QL_i$, $u_{i,j}$ is an $\varepsilon$-quasi-loop.
        \item We have 
        \begin{align*}
            \sum_{i \notin QL} \sum_{j \notin QL_i} |u_{i,j}| \leq \sum_{i \notin QL} (1-c) |u_i| = (1-c) \sum_{i \notin QL} |u_i| \leq (1-c)(1-c)^k|\gamma_n|=(1-c)^{k+1}|\gamma_n|.
        \end{align*}
        \item For all $i \notin QL$, for all $j \notin QL_i, |u_{i,j}| \geq r$
    \end{enumerate}
\end{itemize}
which completes the proof of the recursion. \\ 
Thus, since $0 < 1-c < 1$ and $0 < \lambda < 1$, there exists an integer $k$ such that $(1-c)^k < 1-\lambda$ (simply choose $k=\lceil \frac{\ln(1-\lambda)}{\ln(1-c)} \rceil$), which completes the proof of the lemma.
\end{proof}

Recall that the constants $C$ and $C'$ satisfy : for all primitive elements $\gamma$, $\frac{1}{C}\Vert \gamma \Vert \leq l(\rho(\gamma))\leq C'|\gamma|$. Then, in particular $CC' \geq 1$. Let $0 < \varepsilon < \frac{1}{C}$. Therefore
$$ \frac{1}{C'}(\frac{1}{C}-\varepsilon) < \frac{1}{C'}\frac{1}{C } \leq 1 \hspace{0,5cm} \text{ so } \hspace{0,5cm} 0 < 1-\frac{1}{C'}(\frac{1}{C}-\varepsilon) < 1. $$

\begin{Lemma} \label{contrad}
Let $0 < \varepsilon < \frac{1}{C}$ and $1 - \frac{1}{C'}(\frac{1}{C}-\varepsilon) <  \lambda < 1$. Let $\gamma$ be a primitive element of $\mathbb{F}_2$ which contains disjoint $\varepsilon$-quasi-loops which occupy at least a proportion $\lambda$ of $\gamma$. Then $$d(\rho(\gamma)o,o) < \frac{1}{C} |\gamma|. $$  
\end{Lemma}

\begin{proof}
Since we can find disjoint $\varepsilon$-quasi-loops in $\gamma$ which occupy at least a proportion $\lambda$ of  $\gamma$, there exists $p \in \mathbb{N}, QL \subset \{1,\hdots,p \}$ and some elements $u_1, \hdots, u_p$  such that we can write $\gamma$ in the following way : $\gamma=u_1\hdots u_p$, with for all $i \in QL, u_i$ is an $\varepsilon$-quasi-loop and $\sum_{i \in QL} |u_i| \geq \lambda |\gamma|$. Then we have : 
\begin{equation} \label{prop-pas-quasi-boucle}
    \sum_{i\notin QL} |u_i| = |\gamma| - \sum_{i \in QL} |u_i| \leq |\gamma| - \lambda |\gamma|=(1-\lambda)|\gamma|
\end{equation}

Thus :
\begin{align*}
    d(\rho(\gamma)o,o) & \leq \sum_{i=1}^p d(\rho(u_i)o,o) \text{ by the triangle inequality }\\
    & = \sum_{i \in QL} d(\rho(u_i)o,o) + \sum_{i \notin QL} d(\rho(u_i)o,o) \\
    & \leq \sum_{i \in QL} \varepsilon |u_i| + \sum_{i \notin QL} d(\rho(u_i)o,o) \text{ because } \forall i \in I, u_i \text{ is an $\varepsilon$-quasi-loop} \\
    & \leq \varepsilon \sum_{i\in QL} |u_i| + \sum_{i \notin QL} C'|u_i| \text{ because } \forall u \in \mathbb{F}_2, d(\rho(u)o,o) \leq C'|u| \\
    & \leq \varepsilon |\gamma| + C' \sum_{i \notin QL} |u_i| \\
    & \leq \varepsilon |\gamma| + C'(1-\lambda)|\gamma| \text{ by the inequality } \eqref{prop-pas-quasi-boucle}\\
    & < \varepsilon |\gamma|+ (\frac{1}{C}-\varepsilon)|\gamma| \text{ by the hypothesis on } \lambda \\ 
    & < \frac{1}{C} |\gamma| 
\end{align*}
\end{proof}

This finishes the proof of proposition \ref{uniform-neighborhood}. Indeed : \\ 
Let $0 < \varepsilon < \frac{1}{C}$ and $1 - \frac{1}{C'}(\frac{1}{C}-\varepsilon) <  \lambda < 1$. 
 Lemma \ref{trouve-gamma} gives the existence of a (cyclically reduced) primitive element $\gamma$ and some  $\varepsilon$-quasi-loops in $\gamma$ which occupy at least a proportion $\lambda$ of $\gamma$ and then Lemma \ref{contrad} ensures that $d(\rho(\gamma)o,o) < \frac{1}{C} |\gamma| = \frac{1}{C}\Vert \gamma \Vert$. But since $\gamma$ is primitive, the Bowditch hypothesis (combined with Lemma \ref{gamma-hyp}) states that $\frac{1}{C} \Vert \gamma \Vert \leq d(\rho(\gamma)o,o)$, which is a contradiction. 
\end{proof}

\section{From uniform tubular neighborhoods and Bowditch's hypothesis to primitive-stability} 
\label{second-step-proof}

This section is dedicated to finish the proof of theorem \ref{BQ=PS}, that is that a Bowditch representation is primitive-stable. Pick once and for all a Bowditch representation $\rho$, with constants $(C,D)$,  and let $C'$ be a constant such that for all $\gamma$ in $\F_2$, $d(\rho(\gamma)o,o) \leq |\gamma|$. In the section \ref{first-step-proof}, we prove the existence of a constant $K >0$ such that for all primitive elements $\gamma \in \F_2$, we have the inclusion $\tau_\rho(L_\gamma) \subset N_K(\mathrm{Axis}(\rho(\gamma)))$. (Recall that $L_\gamma$ denotes the axis of $\gamma$ in the Cayley graph of $\F_2$ and $\mathrm{Axis}(\rho(\gamma))$ the axis of the hyperbolic isometry $\rho(\gamma)$, see the beginning of section \ref{first-step-proof}.) For every $\gamma$ in $\mathcal{P}(\F_2)$, pick~$\ell_\gamma$ some geodesic joining the two attracting and repelling point of $\rho(\gamma)$, $\rho(\gamma)^+$ and $\rho(\gamma)^-$. Then~$\ell_\gamma \subset \mathrm{Axis}(\rho(\gamma))$ and by Lemma \ref{nbh-axis}, there exists a constant $C(\delta)$ such that  $N_K(\mathrm{Axis}(\rho(\gamma))) \subset N_{K+C(\delta)}(\ell_\gamma)$. Then, noting $K_\delta=K+C(\delta)$, we obtain that for all primitive elements $\gamma \in \F_2$, we have $\tau_\rho(L_\gamma) \subset N_{K_\delta}(\ell_\gamma)$.   \\ 
Let $p_o$ be some projection of the basepoint $o$ on $\ell_\gamma$. For a point $p$ on the geodesic $\ell_\gamma$, we define the real  $H_\gamma(p)=\pm d(p,p_o)$. The sign plus or minus is determined according to which side of $p_o$ the point $p$ is located on. Thus $H_\gamma$ is an isometry between $\ell_\gamma$ and $\R$ sending $p_o$ to 0. \\ 
We begin by the following lemma : 

\begin{Lemma} \label{right-order}
Let $\gamma$ be a primitive element in $\F_2$ and pick an integer $0\leq i \leq r(\gamma)$. Let $g,g'$ and $g''$ be three points on $L_\gamma$, aligned in this order, such that $d(g,g')=d(g',g'')=l_i(\gamma)$ (hence $d(g,g'')=2l_i(\gamma)$). Denote  $x=\rho(g)o,x'=\rho(g')o,x''=\rho(g'')o$ and choose $p,p'$ and $p''$ respectively  projections of $x,x'$ and $x''$ on the geodesic $\ell_\gamma$. \\
Suppose that $l_i(\gamma)>C(4C'+24\delta+2K_\delta+D)$, then $p,p'$ and $p''$ are aligned in this order on $\ell_\gamma$.
\end{Lemma}

\begin{proof}
We prove this lemma by contraposition. Suppose that the points $p,p'$ and $p''$ are not aligned in this order on the geodesic $l_\gamma$. This means that the reals $H_\gamma(p')-H_\gamma(p)$ and $H_\gamma(p'')-H_\gamma(p')$ are of opposite signs. Without loss of generality, suppose that $H_\gamma(p) \leq H_\gamma(p')$ and $H_\gamma(p') \geq H_\gamma(p'')$. Now consider all the integer points on the segment $[g,g'']$ :
$g_0=g,g_1,\hdots,g_{l_i(\gamma)}=g',g_{l_i(\gamma)+1},\hdots,g_{2l_i(\gamma)}=g''$. For $0 \leq j \leq 2l_i(\gamma)$, denote $x_j=\rho(g_j)o$ and choose $p_j$ a projection of $x_j$ on $\ell_\gamma$ (choose $p_0=p,p_{l_i(\gamma)}=p'$ and $p_{2l_i(\gamma)}=p''$). Therefore, because of our hypothesis on $p,p'$ and $p''$, there exists $0 \leq j \leq l_i(\gamma)-1$ such that, $H_\gamma(p_j) \leq H_\gamma(p_{j+l_i(\gamma)})$ and $H_\gamma(p_{j+1}) \geq H_\gamma(p_{j+l_i(\gamma)+1})$. Hence :
\begin{align*}
    d(p_j,p_{j+l_i(\gamma)}) & = H_\gamma(p_{j+l_i(\gamma)})-H_\gamma(p_j) \qquad \text{ because } H_\gamma(p_j) \leq H_\gamma(p_{j+l_i(\gamma)}) \\
    & = H_\gamma(p_{j+l_i(\gamma)})-H_\gamma(p_{j+l_i(\gamma)+1})+H_\gamma(p_{j+l_i(\gamma)+1})-H_\gamma(p_{j+1})+H_\gamma(p_{j+1})-H_\gamma(p_j) \\
    & \leq d(p_{j+l_i(\gamma)},p_{j+l_i(\gamma)+1})+d(p_{j+1},p_j) \text{ because } H_\gamma(p_{j+l_i(\gamma)+1})-H_\gamma(p_{j+1}) \leq 0 \\
    & \leq d(\rho(g_{j+l_i(\gamma)})o,\rho(g_{j+l_i(\gamma)+1})o)+12\delta+d(\rho(g_{j+1})o,\rho(g_j)o) + 12\delta \text{ by Lemma \ref{alternative2}} \\
    & \leq 2C'+24\delta \qquad \text{ because } d(g_{j+l_i(\gamma)},g_{j+l_i(\gamma)+1})=d(g_j,g_{j+1})=1.
\end{align*}
But since $d(g_j,g_{j+l_i(\gamma)})=l_i(\gamma)$, we have, by Bowditch's hypothesis and Lemma \ref{magic-len} :
\begin{align*}
    \frac{1}{C}l_i(\gamma) -D \leq d(\rho(g_j)o,\rho(g_{j+l_i(\gamma)})o)+2C'.
\end{align*}
Now recall that we have proven that $\tau_\rho(L_\gamma)$ remains in the $K_\delta$-neighborhood of $\ell_\gamma$, then   $$d(\rho(g_j)o,\rho(g_{j+l_i(\gamma)})o)\leq d(p_j,p_{j+l_i(\gamma)})+2K_\delta.$$
Thus, we can bound $l_i(\gamma)$ :
$l_i(\gamma) \leq C(2C'+24\delta+2K_\delta+2C'+D).$
\end{proof}

Let us now prove that $\rho$ is primitive-stable. By contradiction, suppose that it is not. Then for all $n \in \N$, we can find a primitive element $\gamma_n$ and two points $g_n$ and $h_n$ on $L_{\gamma_n}$ such that $d(\rho(g_n)o,\rho(h_n)o) \leq \frac{1}{n} d(g_n,h_n)-1$. Let $x_n=\rho(g_n)o$ and $y_n=\rho(h_n)o$. We have that $d(g_n,h_n) \geq n$.  \\ 
We can make the assumption that the elements $\gamma_n$ are pairwise distinct. Indeed, if the sequence ${(\gamma_n)}_n$ only takes finitely many values, then, up to subsequence, we can suppose that $\gamma_n=\gamma$ for some primitive element $\gamma$. But $\rho(\gamma)$ is an hyperbolic isometry so there exist two constants $C_\gamma$ and $D_\gamma$ (depending on $\gamma$) such that $\tau_\rho(L_\gamma)$ is a $(C_\gamma,D_\gamma)$-quasi-geodesic. Then, since $g_n$ and $h_n$ belong to $L_\gamma$, we have : 
\begin{align*}
    \frac{1}{C_\gamma}d(g_n,h_n)-D_\gamma & \leq d(\rho(g_n)o,\rho(h_n)o) \leq \frac{1}{n}d(g_n,h_n)-1 \\
    \text{so} \qquad \frac{1}{C_\gamma}-\frac{D_\gamma}{d(g_n,h_n)} & \leq \frac{1}{n}-\frac{1}{d(g_n,h_n)}, \\
    \text{then, taking the limit when $n \to \infty$,} \quad \frac{1}{C_\gamma} & \leq 0, \quad \text{which is absurd}.
\end{align*}
Thus we can suppose that the elements $\gamma_n$ are pairwise distinct and therefore $|\gamma_n| \to \infty$. Denote by $[N_1(\gamma_n),\hdots,N_{r(\gamma_n)}(\gamma_n)]$ the continued fraction expansion of $\gamma_n$. As in the proof of the previous section (\ref{first-step-proof}), we can prove the following lemma. 

\begin{Lemma} \label{Ni-borné}
For all $i \in \N^*$, there exists a constant $C_i >0$ such that for all $n \in \N^*$, whenever $N_i(\gamma_n)$ is well defined (that is $r(\gamma_n) \geq i$), we have $N_i(\gamma_n) \leq C_i$. Moreover, up to subsequence, $r(\gamma_n) \to \infty$. 
\end{Lemma}

\begin{proof}
The proof is essentially the same as for Lemma \ref{forme-frac-continue}.
\end{proof}

Now fix an increasing map $\psi : \mathbb{N}^* \longrightarrow \mathbb{R}^*_+$ satisfying $1 \leq \psi(n) \leq n, \, \forall n \in \mathbb{N}^* ;  \, \, \psi(n) \underset{n \to \infty} \longrightarrow + \infty$ and $\psi(n)=o(n)$ (for example, take $\psi(n)=\sqrt{n}$). \\ 
We set $X_n=\{ 0 \leq i \leq r(\gamma_n) : l_i(\gamma_n) \leq \psi(n) \}$. For $n \geq 1$, $X_n$ is non-empty because we always have $0 \in X_n$. Thus the integer $I_n=\max X_n$ is well-defined.

\begin{Lemma} \label{l_I_n to inf}
Up to subsequence, $I_n \underset{n \to \infty} \longrightarrow + \infty$. Moreover, $l_{I_n}(\gamma_n) \underset{n \to \infty} \longrightarrow + \infty$. 
\end{Lemma}

\begin{proof}
If the sequence $(I_n)_n$ were bounded, it would admit a constant subsequence. Then let us fix an integer $I$ such that, up to subsequence, $I_n = I$ for all $n \in \N^*$. By Lemma \ref{Ni-borné}, $r(\gamma_n) \underset{n \to \infty} \longrightarrow + \infty$ so for $n$ sufficiently large, $r(\gamma_n) \geq I+1$. \\
Therefore, using again lemma \ref{Ni-borné} and up to passing to subsequence, we can assume that there exists $N_1,\hdots, N_I$ some positive integers such that for all $1 \leq j \leq I, N_j(\gamma_n)=N_j$. As a consequence, the sequences  $(l_I(\gamma_n))_{n \in \mathbb{N}^*}$ and $(l'_I(\gamma_n))_{n \in \mathbb{N}^*}$ are constants, and we denote $l_I=l_I(\gamma_n),l'_I=l'_I(\gamma_n)$. Therefore we have for all $n \in \mathbb{N}^*, l_I = l_{I_n}(\gamma_n) \leq \psi(n) \leq l_{{I_n}+1}(\gamma_n)=l_{I+1}(\gamma_n)$, so $l_{I+1}(\gamma_n) \underset{n \to \infty}\longrightarrow + \infty$. But $$l_{I+1}(\gamma_n)=(N_{I+1}(\gamma_n)-1)l_I + l'_I \leq (N_{I+1}(\gamma_n)+1)l_I $$ so we deduce that $N_{I+1}(\gamma_n) \underset{n \to \infty} \longrightarrow +\infty$, contradicting Lemma \ref{Ni-borné}. \\
The fact that $l_{I_n}(\gamma(n)) \to \infty$ is now immediate knowing that $l_{I_n}(\gamma_n) \geq I_n +1$ (see the inequalities of remark \ref{l_i ineq}).
\end{proof}

In order to simplify notations, denote $l_{I_n}=l_{I_n}(\gamma_n)$. \\ 
Now consider the segment $[g_n,h_n]$ in $L_{\gamma_n}$ and let us cut it out in subsegments of length $l_{I_n}$, except maybe the last segment that must be of length smaller that $l_{I_n}$. Precisely, consider the Euclidean division of the integer $d(g_n,h_n)$ by $l_{I_n}$ : $d(g_n,h_n)=q_n l_{I_n}+r_n$, with $0 \leq r_n < l_{I_n}$, and set $g_{0,n}=g_n,g_{1,n},\hdots,g_{q_n,n}$ points on $L_{\gamma_n}$ such that $d(g_{k,n},g_{k+1,n}) = l_{I_n}$, $d(g_{q_n,n},h_n)=r_n\leq l_{I_n}$. Moreover, since $l_{I_n} \leq \psi(n) \leq \psi(d(g_n,h_n)) \leq d(g_n,h_n)$, we conclude that $q_n \geq 1$. Now consider $x_{k,n}=\rho(g_{k,n})o$ for $0 \leq k \leq q_n$ the corresponding point in $X$ and finally $p_{k,n}=p(x_{k,n})$ its projection on $\ell_{\gamma_n}$, a geodesic joining the attracting and repelling points of $\rho(\gamma_n)$. 
On one hand, we have the following inequalities :
\begin{align*}
    d(x_{0,n},x_{q_n,n})& \leq d(x_{0,n},y_n)+d(y_n,x_{q_n,n}) \qquad \text{ by the triangle inequality} \\
    & \leq d(x_n,y_n)+d(\rho(h_n)o,\rho(g_{q_n,n})o) \qquad \text{by the definitions of $x_{0,n}$, $x_{q_n,n}$}\\
    & \leq \frac{1}{n}d(g_n,h_n)-1+d(\rho(h_n)o,\rho(g_{q_n,n})o) \qquad \text{ by hypothesis on the points $g_n$ and $h_n$} \\
    & \leq \frac{1}{n}d(g_n,h_n)-1+C'd(g_{q_n,n},h_n) \qquad \text{because $\tau_\rho$ is $C'$-Lipschitz-continuous} \\
    & \leq \frac{1}{n}d(g_n,h_n)-1+C'l_{I_n} \qquad \text{ since } d(g_{q_n,n},h_n)=r_n \leq l_{I_n}.
\end{align*}
On the other hand, since $x_{0,n}$ and $x_{q_n,n}$ belong to $N_{K_\delta}(\ell_{\gamma_n})$, we have that : 
\begin{align*}
    d(p_{0,n},p_{q_n,n}) \leq d(x_{0,n},x_{q_n,n}) +2K_\delta 
\end{align*}
and by Lemma \ref{right-order}, $p_{0,n},p_{1,n},\hdots,p_{q_n,n}$ are aligned in this order on $l_{\gamma_n}$, hence  
\begin{align*}
    d(p_{0,n},p_{q_n,n}) &= \sum_{k=1}^{q_n} d(p_{k-1,n},p_{k,n}) 
\end{align*}
Combining Lemma \ref{magic-len} and the Bowditch hypothesis :
\begin{align*}
    \frac{1}{C}l_{I_n}-D & \leq d(\rho(g_{i-1,n})o,\rho(g_{i,n})o)+2C', \forall 1 \leq i \leq q_n \\
    & = d(x_{i-1,n},x_{i,n})+2C', \forall 1 \leq i \leq q_n \\
    & \leq d(p_{i-1,n},p_{i,n})+2K_\delta+2C', \forall 1 \leq i \leq q_n
\end{align*}
Therefore, by summing :
\begin{align*}
    \frac{q_n}{C}l_{I_n}-Dq_n & \leq \sum_{i=1}^{q_n} d(p_{i-1,n},p_{i,n}) +q_n(2K_\delta+2C') \\
    \frac{q_n}{C}l_{I_n} & \leq \sum_{i=1}^{q_n} d(p_{i-1,n},p_{i,n})+q_n(D+2K_\delta+2C') \\
    & \leq  d(x_{0,n},x_{q_n,n})+2K_\delta + q_n(D+2K_\delta+2C') \\
    & \leq \frac{1}{n}d(g_n,h_n)-1+C'l_{I_n}+2K_\delta+q_n(D+2K_\delta+2C')
\end{align*}

Dividing by $q_nl_{I_n}$ : 
\begin{equation} \label{1/Cleqn}
    \frac{1}{C} \leq \frac{1}{n}\frac{d(g_n,h_n)}{q_nl_{I_n}}+\frac{D+2K_\delta+2C'}{l_{I_n}}+\frac{C'}{q_n}+\frac{2K_\delta-1}{q_nl_{I_n}}
\end{equation}
We now verify that the right hand side on this last inequality tends to zero :
\begin{itemize}
    \item $l_{I_n} \to \infty$ by Lemma \ref{l_I_n to inf}.
    \item We deduce that $q_nl_{I_n} \to \infty$ because $q_n \geq 1$.
    \item $\displaystyle q_n=\frac{d(g_n,h_n)-r_n}{l_{I_n}} \geq \frac{d(g_n,h_n)}{l_{I_n}}-1 \geq \frac{n}{\psi(n)}-1 \to \infty$ because $\psi(n)=o(n)$.
    \item $\displaystyle \frac{d(g_n,h_n)}{q_nl_{I_n}} \leq \frac{(q_n+1)l_{I_n}}{q_n l_{I_n}}=1+ \frac{1}{q_n} \to 1$ because $q_n \to \infty$.
\end{itemize}
Therefore, taking the limit of the inequality \ref{1/Cleqn}, we obtain : $\displaystyle \frac{1}{C} \leq 0 $ which is absurd. Then the representation $\rho$ is primitive-stable. \\ 

\appendix

\section{Length of a path in a $\delta$-hyperbolic space}
\label{length-path}

In this appendix we will be interested in some properties of length of paths in hyperbolic space. The goal is to prove the inequalities of Proposition \ref{long_ext_banane}. Part of the material of this appendix is drawn from or inspired by \cite{coornaert_geometrie_1990}. In particular, the first case of Proposition \ref{long_ext_banane} corresponds to the Lemma 1.8 of Chapter 3 of \cite{coornaert_geometrie_1990}, but the second case of Proposition \ref{long_ext_banane} was not done in this book. We needed some intermediate results to prove the second case so along the way we also included the proof of the first case. \\
Let $X$ be a Gromov-hyperbolic space with hyperbolic constant $\delta$, and suppose that $X$ is geodesic. Denote by $d$ the hyperbolic distance of $X$. We will write $[x,y]$ for some geodesic segment with endpoints $x$ and $y$ in $X$ and $T=[x,y,z]$ for a triangle with vertices $x,y$ and $z$. Recall that we say that a triangle is \emph{$\delta$-thin} if each side of the triangle is included in the $\delta$-neighborhood of the other two. A metric space is \emph{$\delta$-hyperbolic} if every triangle is $\delta$-thin. \\

This first lemma, which is a classical result of hyperbolic geometry, gives the existence of three "close" points in $\delta$-thin triangles. 

\begin{Lemma} \label{triangle-points-proches}
Let $(X,d)$ be a metric space and $T=[x,y,z]$ a $\delta$-thin triangle of $X$. There exists $r \in [x,y], s \in [y,z]$ and $t \in [x,z]$ such that $d(r,s)\leq \delta$ and $d(r,t) \leq \delta$. (Hence in addition $d(s,t) \leq 2\delta$)
\end{Lemma}

\begin{proof}
Consider :
\begin{align*}
    L & = \{ r \in [x,y] : \exists t \in [x,z], d(r,t) \leq \delta \} \\
    R & = \{ r \in [x,y] : \exists s \in [z,y], d(r,s) \leq \delta \} 
\end{align*}
Then :
\begin{itemize}
    \item $R$ and $L$ are non-empty because $x \in L$ and $y \in R$.
    \item $R$ and $L$ are closed because $[x,z]$ and $[z,y]$ are compact.
    \item We have $[x,y]=R \cup L$ since the triangle $[x,y,z]$ is $\delta$-thin. 
\end{itemize}
We deduce that $R \cap L \neq \emptyset$. 
Indeed, if $R \cap L = \emptyset$, we would obtain an open cover of $[x,y]$ with two disjoint non-empty open sets ($[x,y]\backslash L$ and $[x,y] \backslash R$), which would contradict the connectedness of $[x,y]$.  
Thus we deduce the existence of $r \in [x,y]$, $s \in [y,z]$ and $t \in [x,z]$ such that $d(r,s) \leq \delta$ and $d(r,t) \leq \delta$. 

\end{proof}

The following lemma generalises the notion of $\delta$-thin triangle in $\delta$-hyperbolic spaces and is taken from \cite{coornaert_geometrie_1990}. We include the proof for completeness. 

\begin{Lemma} \label{gen-delta-fin}
Let $X$ be a $\delta$-hyperbolic geodesic space and $Y=[x_0,x_1]\cup[x_1,x_2]\cup \hdots \cup [x_{n-1},x_n]$ a chain of $n$ geodesic segments, with $n \leq 2^k$, where $k$ is an integer such that $k \geq 1$. Then, for any point $x$ in a geodesic segment $[x_0,x_n]$, we have $d(x,Y) \leq k\delta$.
\end{Lemma}

\begin{proof}
Let us proceed by induction on $k \geq 1$.
\begin{itemize}
    \item If $k=1$, that is $n=2$, it is the case of a triangle. Since $X$ is $\delta$-hyperbolic, then the triangles are $\delta$-thin, so we have the requested inequality.
    \item Assume that the property is true for some $k\geq 1$ and let us consider $n$ geodesic segments, with $n \leq 2^{k+1}$. After possibly artificially adding points on $Y$, we can assume that $n=2^{k+1}$. Let $x \in [x_0,x_n]$. The triangle with vertices $x_0,x_{\frac{n}{2}}$ and $x_n$ is $\delta$-thin which ensures the existence of $m \in [x_0,x_{\frac{n}{2}}]\cup[x_{\frac{n}{2}},x_n]$ such that $d(x,m)\leq \delta$. Without loss of generality, suppose that $m\in[x_0,x_{\frac{n}{2}}]$. We have $\frac{n}{2}\leq 2^k $ so by induction we can find a point $m' \in [x_0,x_1] \cup \hdots \cup [x_{\frac{n}{2}-1},x_{\frac{n}{2}}]$ such that $d(m,m') \leq k \delta$. Thus, by the triangle inequality, we obtain : $d(x,m')\leq d(x,m)+d(m,m')\leq \delta+k\delta=(k+1)\delta$. \qedhere
\end{itemize}
\end{proof}

In particular, by Lemma \ref{gen-delta-fin}, if we consider a hyperbolic quadrilateral with vertices $x,y,y_1,x_1$ (in this order), then every point of $[x,y]$ is at a distance at most $2\delta$ of a point of $[x,x_1]\cup[x_1,y_1]\cup[y_1,y]$. Thus we have the following alternative :
\begin{itemize}
    \item $\forall z \in [x,y], d(z,[x_1,y_1]) > 2\delta$ and in this case $\forall z \in [x,y], d(z,[x,x_1]\cup[y,y_1]) \leq 2\delta$, or
    \item $\exists z \in [x,y], d(z,[x_1,y_1]) \leq 2 \delta$.
\end{itemize}

Several of the following lemmas depend on this alternative. The following lemma finds three "close" points in the quadrilateral in the first case of this alternative. 

\begin{Lemma}
\label{quadri-hyp}
Let $X$ be a $\delta$-hyperbolic geodesic space, and $x,x_1,y_1,y \in X$ (in this order) be the vertices of a hyperbolic quadrilateral of $X$. We further assume that for all $z \in [x,y], d(z,[x_1,y_1]) >2\delta$. Then there exists $z\in [x,y], r \in [x,x_1]$ and $s \in [y,y_1]$ such that $d(z,r) \leq 2 \delta$ and $d(z,s) \leq 2 \delta$. 
\end{Lemma}

\begin{proof} The proof follows the same lines as the proof of Lemma \ref{triangle-points-proches}. \\

By Lemma \ref{gen-delta-fin}, every point $z \in [x,y]$ is at a distance at most $2\delta$ of $[x,x_1]\cup[x_1,y_1]\cup[y_1,y]$. But, by hypothesis, for all $z \in [x,y], d(z,[x_1,y_1]) >2\delta$, so for all $z \in [x,y]$, there exists $z' \in [x,x_1]\cup[y_1,y]$ such that $d(z,z') \leq 2\delta$. \\

Consider :
\begin{align*}
    L & = \{ z \in [x,y] : \exists r \in [x,x_1], d(z,r) \leq 2\delta \} \\
    R & = \{ z \in [x,y] : \exists s \in [y,y_1], d(z,s) \leq 2\delta \} 
\end{align*}
Then :
\begin{itemize}
    \item $R$ and $L$ are non-empty because $x \in L$ and $y \in R$.
    \item $R$ and $L$ are closed because $[x,x_1]$ and $[y,y_1]$ are compact.
    \item We have $[x,y]=R \cup L$ because, by Lemma \ref{gen-delta-fin}, every point of $[x,y]$ is at a distance at most $2\delta$ of $[x,x_1]\cup[x_1,y_1]\cup[y_1,y]$, and, by hypothesis, every point of $[x,y]$ is at a distance at least $2\delta$ of $[x_1,y_1]$.
\end{itemize}
We deduce that $R \cap L \neq \emptyset$. 
Indeed, if $R \cap L = \emptyset$, we would have an open cover of $[x,y]$ in two disjoints non-empty open sets ($[x,y]\backslash L$ and $[x,y] \backslash R$), which would contradict the connectedness of $[x,y]$. 
Thus we deduce the existence of $z \in [x,y]$, $r \in [x,x_1]$ and $s \in [y,y_1]$ such that $d(z,r) \leq 2\delta$ and $d(z,s) \leq 2\delta$. 
\end{proof}

A quadrilateral $[x,y,x_1,y_1]$ has thus one of two typical general shapes : the one where a point of $[x,y]$ is close to a point of $[x_1,y_1]$ and the one where all points of $[x,y]$ are far from $[x_1,y_1]$. The following lemma clarifies this alternative in terms of comparing the lengths of the sides of the quadrilateral, in the more specific case where $x_1$ and $y_1$ are projections of $x$ and $y$ on a geodesic.

\begin{Lemma} \label{alternative1}
Let $X$ be a $\delta$-hyperbolic geodesic space and $l$ a bi-infinite geodesic of $X$. Let $x$ and $y$ be two points in $X$. Let $K_x=d(x,l)$, $K_y=d(y,l)$ and $d=d(x,y)$. Consider $x_1$ and $y_1$ two projections of $x$ and $y$ on $l$, that is two points $x_1, y_1 \in l$ satisfying $d(x,x_1)=K_x$ and $d(y,y_1)=K_y$. We denote by $[x_1,y_1]$ the geodesic segment included in $l$ with endpoints $x_1$ and $y_1$.

\begin{itemize}
    \item If $\exists z \in [x,y], d(z,[x_1,y_1]) \leq 2 \delta$, then $d \geq K_x+K_y-4 \delta$
    \item If $\forall z \in [x,y], d(z,[x_1,y_1]) > 2 \delta$, then $d \leq K_x+K_y +4\delta$
\end{itemize}

\end{Lemma}

\begin{proof}
\begin{itemize}
    \item Suppose that there exists $z \in [x,y]$ and $z_1 \in [x_1,y_1]$ such that $d(z,z_1)\leq 2\delta$. By definition, $x_1$ minimizes the distance from $x$ to $l$, and $z_1 \in l$ so $d(x,z_1)\geq K_x$. Then :
    \begin{equation}
    \label{proche-1}
        d(x,z) \geq d(x,z_1)-d(z,z_1)\geq K_x-2\delta 
    \end{equation}
    Similarly, since $y_1$ minimizes the distance from $y$ to $l$, we get :
    \begin{equation}
        \label{proche-2}
        d(z,y) \geq d(y,z_1) -d(z,z_1) \geq K_y - 2 \delta
    \end{equation}
    Thus, since $z \in [x,y]$, we obtain by combining \eqref{proche-1} and \eqref{proche-2} :
    \begin{equation*}
        d=d(x,y)=d(x,z)+d(z,y) \geq K_x+K_y-4\delta.
    \end{equation*}
    
    \item Suppose that for all $z \in [x,y], d(z,[x_1,y_1]) > 2\delta$. Then by Lemma \ref{quadri-hyp}, there exists $z \in [x,y]$, $r \in [x,x_1]$ and $s \in [y,y_1]$ such that $d(z,r) \leq 2 \delta$ and $d(z,s) \leq 2 \delta$. Hence we have on one hand : 
    \begin{align*}
        d(x,z) & \leq d(x,r)+d(r,z) \\
        & = d(x,x_1)-d(r,x_1)+d(r,z) \qquad \text{ because } r \in [x,x_1] \\
        & \leq K_x + 2 \delta \qquad \text{ since } d(r,z) \leq 2 \delta. \\
        \intertext{On the other hand, we obtain in the same way :} 
        d(z,y) & \leq K_y + 2 \delta
    \end{align*}
    and then :
    \begin{equation*}
        d=d(x,z)+d(z,y)\leq K_x +K_y + 4 \delta.
    \end{equation*}
\end{itemize}
\end{proof}

The following lemma and corollary aim to bound, in the same context as in the previous lemma, the distance between $x_1$ and $y_1$.

\begin{Lemma} \label{alternative2}
Let $X$ be a $\delta$-hyperbolic geodesic space, and $l$ a bi-infinite geodesic of $X$. Let $x$ and $y$ be two points in $X$, denote $K_x=d(x,l)$ and $K_y=d(y,l)$. Consider $x_1$ and $y_1$ two projections of $x$ and $y$ on $l$, that is two points $x_1, y_1 \in l$ satisfying $d(x,x_1)=K_x$ and $d(y,y_1)=K_y$. We denote by $[x_1,y_1]$ the geodesic segment included in $l$ with endpoints $x_1$ and $y_1$. We let $d=d(x,y)$ and $d_1=d(x_1,y_1)$. Then we have the following alternative : 
\begin{itemize}
    \item If $\forall z \in [x_1,y_1], d(z,[x,y])>2\delta$, then $d_1 \leq 8 \delta$.
    \item If $\exists z \in [x_1,y_1], d(z,[x,y]) \leq 2 \delta$, then $d_1 \leq d-K_x-K_y+12\delta$.
\end{itemize}
In particular, we always have the inequality $d_1 \leq d +12 \delta$. 
\end{Lemma}

\begin{proof}
\begin{itemize}
    \item Suppose that for all $z \in [x_1,y_1]$ one has $d(z,[x,y])>2\delta$. Denote $[x,x_1]$ (respectively $[y,y_1]$) a geodesic segment with endpoints $x$ and $x_1$ (respectively $y$ and $y_1$), then, by Lemma \ref{quadri-hyp} the exists $z \in [x_1,y_1]$, $s \in [x,x_1]$ and $r \in [y,y_1]$ such that $d(z,s) \leq 2 \delta$ and $d(z,r)\leq 2 \delta$. Therefore 
    \begin{align*}
        d(x_1,z) & \leq d(x_1,s) +d(s,z) \\
         & \leq d(x_1,s)+2\delta \qquad \text{ because } d(s,z) \leq 2\delta \\
         & = d(x,x_1)-d(x,s)+2\delta \qquad \text{ since } s\in [x,x_1] \\
         & \leq d(x,x_1)-d(x,z)+2\delta +2\delta \qquad \text{ because } d(x,z) \leq d(x,s)+d(s,z) \leq d(x,s) + 2 \delta. 
    \end{align*}
    But $x_1$ is a projection of $x$ to $l$ and $z \in l$, so $d(x,x_1) \leq d(x,z)$. Thus : 
    \begin{equation}
    \label{cas-loin-1}
        d(x_1,z) \leq 4 \delta.
    \end{equation}
   Similarly, we show that :
    \begin{equation}
        \label{cas-loin-2}
        d(z,y_1)\leq 4 \delta.
    \end{equation}
    Thus, by summing \eqref{cas-loin-1} and \eqref{cas-loin-2}, we obtain :
    \begin{equation*}
        d_1 \leq 8 \delta.
    \end{equation*}
    
    \item Suppose that there exists $z \in [x,y]$ and $z_1 \in [x_1,y_1]$ such that $d(z,z') \leq 2 \delta$. Consider $[x,z_1]$ a geodesic segment with endpoints $x$ and $z_1$. The triangle $[x,x_1,z_1]$ is $\delta$-thin so, by Lemma \ref{triangle-points-proches} there exists $r\in [x,x_1]$, $s \in [x,z_1]$ and $t \in [x_1,z_1]$ such that $d(t,r) \leq \delta$ and $d(t,s)\leq \delta$. In particular, $d(r,s) \leq 2 \delta$. Thus, we have : 
    \begin{align*}
        d(x_1,t) & \leq d(x_1,r)+d(r,t) \\
        & \leq d(x_1,r)+\delta \qquad \text{ because } d(r,t) \leq \delta \\
        & = d(x,x_1)-d(x,r) +\delta \qquad \text{ because } r \in [x,x_1] \\
        & \leq d(x,x_1)-d(x,t)+\delta +\delta \qquad \text{ since } d(x,t) \leq d(x,r)+d(r,t) \leq d(x,r) +\delta
    \end{align*}
    But $x_1$ is a projection of $x$ on $l$ and $t \in l$ so $d(x,x_1) \leq d(x,t)$. Thus : 
    \begin{equation}
    \label{cas-proche-1}
        d(x_1,t) \leq 2\delta.
    \end{equation}
    Furthermore :
    \begin{align*}
        d(t,z_1) & \leq d(t,s) +d(s,z_1) \\
        & \leq \delta + d(s,z_1) \qquad \text{ because } d(t,s) \leq \delta \\
        & = d(x,z_1)-d(x,s) + \delta \qquad \text{ since } s \in [x,z_1] \\
        & \leq d(x,z)+d(z,z_1)-d(x,t)+\delta + \delta \qquad \text{ because } d(x,t) \leq d(x,s)+d(s,t) \leq d(x,s)+\delta \\
        & \leq d(x,z) + 2 \delta -d(x,t) + 2\delta \qquad \text{ since } d(z,z_1) \leq 2 \delta 
    \end{align*}
    and $d(x,t) \geq d(x,x_1)=K_x$, so
    
    \begin{equation}
    \label{cas-proche-2}
        d(t,z_1) \leq d(x,z) -K_x + 4\delta.
    \end{equation}
    
    Thus, combining \eqref{cas-proche-1} and \eqref{cas-proche-2}, we get: 
    \begin{equation}
        \label{cas-proche-3}
        d(x_1,z_1) = d(x_1,t)+d(t,z_1) \leq d(x,z)-K_x+6\delta
    \end{equation}
    By the same reasoning, considering the triangle $[y,y_1,z_1]$, we get :
    \begin{equation}
        \label{cas-proche-4}
        d(z_1,y_1) \leq d(z,y)-K_y+6\delta
    \end{equation}
    We deduce, using \eqref{cas-proche-3} and \eqref{cas-proche-4} :
    \begin{align*}
        d_1 & = d(x_1,z_1)+d(z_1,y_1) \\
        & \leq d(x,z)-K_x+6\delta + d(z,y)-K_y+6\delta \\
        & = d(x,y)-K_x-K_y+12 \delta \text{ car } z \in [x,y]
    \end{align*}
    Thus :
    \begin{equation*}
        d_1 \leq d -K_x-K_y+12\delta.
    \end{equation*}
\end{itemize}
\end{proof}

Thus, we deduce an upper bound on $d_1$ depending only on the hyperbolic constant  $\delta$ in the case when $d \leq K_x + K_y +6\delta$. (The choice of this particular threshold on $d$ comes from the disjunction made in the proposition \ref{long_ext_banane} at the end of this section.)

\begin{Corollary} \label{d1-delta}
Let $X$ be a $\delta$-hyperbolic geodesic space, and $l$ a bi-infinite geodesic of~$X$. Let $x$ and $y$ be two points in $X$, denote $K_x=d(x,l)$ and $K_y=d(y,l)$. Consider $x_1$ and $y_1$  two projections of $x$ and $y$ on $l$, that is two points $x_1, y_1 \in l$ satisfying $d(x,x_1)=K_x$ and $d(y,y_1)=K_y$. Denote $d=d(x,y)$ and $d_1=d(x_1,y_1)$. Then : 
\begin{equation*}
    \text{  If } d \leq K_x+K_y+6\delta, \text{ then } d_1 \leq 18 \delta.
\end{equation*}

\end{Corollary}

\begin{proof} 
According to Lemma \ref{alternative2}, we have $d_1 \leq \max(8\delta,d-K_x-K_y+12\delta)$. But $d \leq K_x+K_y+6\delta$, so 
\begin{equation*}
    d_1 \leq \max(8\delta,d-K_x-K_y+12\delta) \leq \max(8\delta,18\delta) = 18\delta.
\end{equation*}
\end{proof}

The following lemmas, \ref{ineg-Ldzft}, \ref{ineg-dzft} and \ref{L-petits-bouts}, as well as the propositions \ref{ineg-L-gen} and \ref{long_ext_banane}, aim at bounding from below the length of a path remaining at a distance at least $K$ from a given geodesic. This first lemma, which is a first step, minimizes this length using the distance of a point $z$ on a geodesic connecting the extremal points of the path and a point on the path. \\

\begin{Lemma} \label{ineg-Ldzft}
Let $X$ be a $\delta$-hyperbolic geodesic space, $l$ a bi-infinite geodesic of~$X$ and $K>0$. Consider $f : [a,b] \to X$ a continuous rectifiable path such that : 
\begin{align*}
    \forall t \in [a,b], \hspace{0,5cm} d(f(t),l) \geq K
\end{align*}
Let $x=f(a), y=f(b), d=d(x,y)$ and consider $[x,y]$ a geodesic segment with endpoints $x$ and $y$. Denote $L=\mathrm{length}(f([a,b]))$ and suppose $L > 2 \delta$. \\ 
Then, for all $z \in [x,y]$, there exists $t \in [a,b]$ such that : 
\begin{equation*}
    L \geq (2^{\frac{d(z,f(t))}{\delta}-1}-2)\delta.
\end{equation*}
\end{Lemma}

\begin{proof}
Since the path $f([a,b])$ is rectifiable, there exists a subdivision of $f([a,b])$ in $n$ arcs, each of length $2\delta$ except possibly the last one of length smaller than $2\delta$ (then $n=\lceil \frac{L}{2\delta} \rceil$). Denote by $y_0=x,y_1,\hdots,y_n=y$ the points of this subdivision and consider $Y=[y_0,y_1]\cup[y_1,y_2]\cup \hdots \cup [y_{n-1},y_n]$ the chain of $n$ geodesic segments connecting these points. \\
We can therefore bound the length of the path $f([a,b])$ from below by the length of the first $n-1$ paths of this subdivision: 
\begin{equation} \label{subd}
    L \geq (n-1)2\delta
\end{equation}

The hypothesis $L > 2\delta$ allows us to assert that $n \geq 2$. Then there exists $k \geq 1$ such that $2^{k-1} < n \leq 2^k$. We deduce : 
\begin{equation} \label{subd-k}
    L \geq (2^{k-1}-1)2\delta=(2^k-2)\delta
\end{equation}
We can apply Lemma \ref{gen-delta-fin} to $z \in [y_0,y_n]$ and $Y$ : there exists $p \in Y$ such that $d(z,p)\leq k\delta$. Moreover, there exists $i \in \{0,n-1\}$ such that $p \in [y_i,y_{i+1}]$, and $d(y_i,y_{i+1}) \leq 2\delta$ (because the distance between $y_i$ and $y_{i+1}$ is in particular smaller than $2\delta$) so, there exists $j \in \{0,n\}$ such that $d(p,y_j)\leq \delta$. \\
Let $z \in [x,y]$. Then, by the triangle inequality, we obtain the upper bound $d(z,y_j)\leq d(z,p)+d(p,y_j)\leq k\delta + \delta = (k+1)\delta$ from which we deduce $k \geq d(z,y_j) / \delta -1$. Moreover, since $y_j \in f([a,b])$, we deduce the existence of $t\in [a,b]$ such that $y_j=f(t)$. Thus from the inequality \eqref{subd-k} we obtain the desired inequality: 
\begin{equation*}
        L \geq (2^{\frac{d(z,f(t))}{\delta} -1}-2)\delta.
\end{equation*} 
\end{proof}

Now, we try to give a lower bound on the term $d(z,f(t))$ in the previous lemma \ref{ineg-Ldzft}. We give two lower bounds according to the shape of the quadrilateral.

\begin{Lemma} \label{ineg-dzft}
Let $X$ be a $\delta$-hyperbolic geodesic space, $l$ a bi-infinite geodesic of $X$ and $K>0$. Consider $f : [a,b] \to X$ a continuous rectifiable path such that : 
\begin{align*}
    \forall t \in [a,b], \hspace{0,5cm} d(f(t),l) \geq K.
\end{align*}
Denote $x=f(a), y=f(b), d=d(x,y), K_x=d(x,l), K_y=d(y,l)$ and $L=\mathrm{length}(f([a,b]))$. Let $x_1$ and $y_1$ be two projections of $x$ and $y$ on $l$. \\
Consider $[x,y]$a geodesic segment with endpoints $x$ and $y$ and $[x_1,y_1]$ the geodesic segment included in $l$ with endpoints $x_1$ and $y_1$. \\
We can bound from below the distance from a point $z$ of $[x,y]$ to any point on the path $f([a,b])$ in the following two cases :

\begin{itemize}
  \item Suppose that $\forall z \in [x,y], d(z,[x_1,y_1]) > 2 \delta$. Then : \\
    \begin{equation} \label{ineg-dzft-1}
         \exists z \in [x,y], \forall t \in [a,b], d(z,f(t)) \geq \frac{d}{2} + K -\frac{1}{2}(K_x+K_y)-4 \delta
    \end{equation}
    In particular, if $K_x=K_y=K$, we have :
    \begin{equation} \label{ineg-dzft-2}
        \exists z \in [x,y], \forall t \in [a,b], d(z,f(t)) \geq \frac{d}{2} -4 \delta
    \end{equation}

    \item Suppose that $\exists z \in [x,y], d(z,[x_1,y_1]) \leq 2 \delta$. Then : \\
    \begin{equation} \label{ineg-dzft-3}
        \exists z \in [x,y], \forall t \in [a,b], d(z,f(t)) \geq K - 2 \delta
    \end{equation}
\end{itemize}
\end{Lemma}

\begin{proof}
\begin{itemize}
    \item Suppose that $\forall z \in [x,y], d(z,[x_1,y_1]) > 2 \delta$. Then, considering a quadrilateral with vertices $x,x_1,y_1,y$ (in this order), we deduce by Lemma \ref{quadri-hyp} the existence of $z \in [x,y], r \in [x,x_1]$ and $s \in [y,y_1]$ such that $d(z,r) \leq 2 \delta$ and $d(z,s)\leq 2\delta$. Let $t \in [a,b]$. Then : 
    \begin{align*}
        K & \leq d(f(t),l) \qquad \text{ by hypothesis} \\
          & \leq d(f(t),x_1) \qquad \text{ because } x_1 \in l \\
          & \leq d(f(t),z) +d(z,r)+d(r,x_1) \qquad \text{ by the triangle inequality } \\
          & \leq d(f(t),z) + 2 \delta + d(r,x_1) \qquad \text{ because } d(z,r) \leq 2 \delta \\
          & = d(f(t),z) + 2 \delta + d(x,x_1) - d(x,r) \qquad \text{ because } r \in [x,x_1] \\
          & \leq d(f(t),z) + 2 \delta + K_x - d(x,r) \qquad \text{ since $x_1$ is a projection of $x$ on $l$} \\
          & \leq d(f(t),z) + 2 \delta + K_x - d(x,z) + d(z,r) \qquad \text{ by the triangle inequality } \\
          & \leq d(f(t),z) + 2 \delta + K_x - d(x,z) + 2 \delta \qquad \text{ because } d(z,r) \leq 2 \delta.
    \end{align*}
    We deduce :
    \begin{equation} \label{dxz}
        d(x,z)+K-K_x-4\delta \leq d(f(t),z).
    \end{equation}
    By analogous reasoning changing $x$ to $y$, $x_1$ to $y_1$, $r$ to $s$ and $K_x$ to $K_y$, we also get : 
    \begin{equation} \label{dzy}
        d(z,y)+K-K_y-4\delta \leq d(f(t),z).
    \end{equation}
    Thus, averaging the inequalities \eqref{dxz} and \eqref{dzy}, and since  $d=d(x,y)=d(x,z)+d(z,y)$ we obtain the desired inequality :
    \begin{equation*}
        \frac{d}{2}+K-\frac{1}{2}(K_x+K_y)-4\delta \leq d(f(t),z).
    \end{equation*}
    
    \item Suppose that $\exists z \in [x,y], d(z,[x_1,y_1] \leq 2 \delta$. Then let $z' \in [x_1,y_1]$ such that $d(z,z') \leq 2 \delta$. Let $t\in [a,b]$. We have :
    \begin{align*}
        K & \leq d(f(t),z') \qquad \text{ because } z' \in l \\
          & \leq d(f(t),z)+d(z,z') \qquad \text{ by the triangle inequality} \\
          & \leq d(f(t),z)+2\delta  \qquad \text{ because } d(z,z') \leq 2 \delta 
    \end{align*}
    We have indeed shown that 
    \begin{equation*}
        d(f(t),z) \geq K -2 \delta.
    \end{equation*}
\end{itemize}
\end{proof}

The following proposition is a direct consequence of the two previous lemmas. It gives a lower bound on the length of a path remaining at a distance at least $K$ from a geodesic as a function of the distance between the extremal points of the path, the distances of the extremal points to the geodesic and $K$. 

\begin{Proposition} \label{ineg-L-gen}
Let $X$ be a $\delta$-hyperbolic geodesic space, $l$ a bi-infinite geodesic of $X$ and $K>0$.Consider $f : [a,b] \to X$ a continuous rectifiable path such that : 
\begin{align*}
    \forall t \in [a,b], \hspace{0,5cm} d(f(t),l) \geq K
\end{align*}
Denote $x=f(a), y=f(b), d=d(x,y), K_x=d(x,l), K_y=d(y,l)$ and $L=\mathrm{length}(f([a,b]))$ and suppose that $L > 2 \delta$. Let $x_1$ and $y_1$ be some projections of $x$ and $y$ on $l$. \\
Consider $[x,y]$ a geodesic segment with endpoints $x$ and $y$ and $[x_1,y_1]$ the geodesic segment included in $l$ with endpoints $x_1$ and $y_1$. \\
\begin{itemize}
    \item Suppose that $\forall z \in [x,y], d(z,[x_1,y_1]) > 2 \delta$. Then :
    \begin{equation} \label{ineg-L-gen-1}
        L \geq (2^{\frac{d-K_x-K_y+2K}{2\delta}-5}-2)\delta.
    \end{equation}
    In particular, if $K_x=K_y=K$, then :
    \begin{equation} \label{ineg-L-gen-2}
        L \geq (2^{\frac{d}{2\delta}-5}-2)\delta.
    \end{equation}
    
    \item Suppose that $\exists z \in [x,y], d(z,[x_1,y_1]) \leq 2 \delta$. Then :
    \begin{equation}\label{ineg-L-gen-3}
        L \geq (2^{\frac{K}{\delta}-3}-2)\delta.
    \end{equation}
\end{itemize}
\end{Proposition} 

\begin{proof} 
\begin{itemize}
    \item Suppose that $\forall z \in [x,y], d(z,[x_1,y_1]) > 2 \delta$. \\
    According to the property \eqref{ineg-dzft-1} of Lemma \ref{ineg-dzft} , there exists $z \in [x,y]$ such that for all $t \in [a,b], d(z,f(t)) \geq \frac{d}{2} + K -\frac{1}{2}(K_x+K_y)-4 \delta$. Then, by Lemma \ref{ineg-Ldzft}, we obtain the existence of $t \in [a,b]$ such that
\begin{equation*}
    L \geq (2^{\frac{d(z,f(t))}{\delta}-1}-2)\delta.
\end{equation*}
Indeed we have :
\begin{equation*}
    L \geq (2^{\frac{\frac{d}{2} + K -\frac{1}{2}(K_x+K_y)-4 \delta}{\delta}-1}-2)\delta = (2^{\frac{d-K_x-K_y+2K}{2\delta}-5}-2)\delta.
\end{equation*}

    \item Suppose that $\exists z \in [x,y], d(z,[x_1,y_1]) \leq 2 \delta$. \\
    According to the property \eqref{ineg-dzft-3} of Lemma \ref{ineg-dzft}, there exists $z \in [x,y]$ such that for all $t \in [a,b], d(z,f(t)) \geq K - 2 \delta$. Then, by Lemma \ref{ineg-Ldzft}, we obtain the existence of $t \in [a,b]$ such that
    \begin{equation*}
        L \geq (2^{\frac{d(z,f(t))}{\delta}-1}-2)\delta.
    \end{equation*}
    Indeed we have :
    \begin{equation*}
    L \geq (2^{\frac{K-2\delta}{\delta}-1}-2)\delta=(2^{\frac{K}{\delta}-3}-2)\delta.
    \end{equation*}
\end{itemize}
\end{proof}

The following lemma gives a lower bound on the length of a path remaining at a distance at least $K$ from a geodesic in the specific case where the distance between the endpoints is known as a function of the distances of the endpoints to the geodesic.
 
\begin{Lemma} \label{L-petits-bouts}
Let $X$ be a $\delta$-hyperbolic geodesic, $l$ a bi-infinite geodesic of $X$ and $K>0$. Consider $f : [a,b] \to X$ a continuous rectifiable path such that~: 
\begin{align*}
    \forall t \in [a,b], \hspace{0,5cm} d(f(t),l) \geq K
\end{align*}
Denote $x=f(a), y=f(b), d=d(x,y), K_x=d(x,l), K_y=d(y,l)$ and $L=\mathrm{length}(f([a,b]))$ and suppose that $L > 2 \delta$.\\
Assume that $d=K_x+K_y+6\delta$, then : 
\begin{equation*}
    L \geq (2^{\frac{K}{\delta}-3}-2)\delta
\end{equation*}
\end{Lemma}

\begin{proof}
Let $x_1$ and $y_1$ be two projections of $x$ and $y$ on $l$. \\
Consider $[x,y]$ a geodesic segment with endpoints $x$ and $y$ and $[x_1,y_1]$ the geodesic segment included in $l$ with endpoints $x_1$ and $y_1$. \\
We distinguish two cases :
\begin{itemize}
    \item If $\exists z \in [x,y], d(z,[x_1,y_1]) \leq 2 \delta$, then the inequality \eqref{ineg-L-gen-3} of the proposition \ref{ineg-L-gen} gives $$L \geq (2^{\frac{K}{\delta}-3}-2)\delta.$$
    \item If $\forall z \in [x,y], d(z,[x_1,y_1]) > 2 \delta$, then the inequality  \eqref{ineg-L-gen-1} of the proposition \ref{ineg-L-gen} gives~: 
    \begin{align*}
        L & \geq (2^{\frac{d-K_x-K_y+2K}{2\delta}-5}-2)\delta \\
         & \geq (2^{\frac{6\delta+2K}{2\delta}-5}-2)\delta \text{ because } d=K_x+K_y+6\delta \\
         & \geq (2^{\frac{K}{\delta}-3}-2)\delta.
    \end{align*}
\end{itemize}
\end{proof}

The following proposition is the most important one in this appendix. It gives a lower bound on the length of a path remaining at a distance at least $K$ from a geodesic and whose extremal points are "approximately" at distance $K$ from this geodesic. It distinguishes two regimes : one when the distance to the geodesic is very large compared to the distance between the extremal points and the other in the opposite case. In the first case, the path length grows exponentially with the distance between the extremal points and in the second case, this growth is linear, but with a large constant.

\begin{Proposition} \label{long_ext_banane}
Let $X$ be a $\delta$-hyperbolic geodesic space, $l$ a bi-infinite geodesic of $X$, $K>0$ and $C>0$ two constants. Consider $f : [a,b] \to X$ a continuous rectifiable path such that : 
\begin{align*}
    \forall t \in [a,b], \hspace{0,5cm} d(f(t),l) \geq K \\
    d(f(a),l)\leq K+C \\
    d(f(b),l)\leq K+C 
\end{align*}
Denote $x=f(a), y=f(b), d=d(x,y)$ and $L=\mathrm{length}(f([a,b]))$.

\begin{itemize}
    \item Suppose that $d \leq 2K +6\delta$. Then, after denoting $C'=\max(C,\delta)$ :
    \begin{equation}\label{L-d<K}
        L \geq (2^{\frac{d}{2\delta}-\frac{C'}{\delta}-5}-2)\delta.
    \end{equation}
    \item Suppose that $d > 2K + 6 \delta$. Then there exists an integer $n \geq 2$ such that :
    \begin{equation} \label{L-d>K-n}
    \left\{ 
    \begin{array}{ll} 
        L  \geq (n-1) (2^{\frac{K}{\delta}-3}-2)\delta \\ 
        d  \leq 18n\delta +2K +2C
    \end{array}
    \right. 
    \end{equation}
    
    In particular, we deduce : 
    \begin{equation} \label{L-d>K}
        L \geq \frac{1}{18}(d-2K-2C-18\delta)(2^{\frac{K}{\delta}-3}-2).
    \end{equation}
\end{itemize}

\end{Proposition}

\begin{proof}
Let us first treat separately the case where $L \leq 2\delta$. In this case, we also have $d\leq 2\delta$. Then $d\leq 2K + 6\delta$ on one hand, and on the other hand  $(2^{\frac{d}{2\delta}-6}-2)\delta \leq (2^{-5}-2)\delta < 0 \leq L$, so \eqref{L-d<K} holds. \\
Let us now suppose that $L > 2\delta$. Denote $K_x=d(x,l)$ and $K_y=d(y,l)$. Then by hypothesis, $K_x-K\leq C$ and $K_y-K \leq C$. \\

\begin{itemize}
    \item Suppose that $d \leq 2K+6\delta$. \\
    \begin{itemize}
        \item Suppose that $\forall z \in [x,y], d(z,[x_1,y_1]) > 2 \delta$. Then, according to the inequality \eqref{ineg-L-gen-1} of Lemma \ref{ineg-L-gen}, we have :
        \begin{equation*} 
         L \geq (2^{\frac{d-K_x-K_y+2K}{2\delta}-5}-2)\delta
            \geq (2^{\frac{d-2C}{2\delta}-5}-2)\delta \geq (2^{\frac{d}{2\delta}-\frac{C}{\delta}-5}-2)\delta \geq (2^{\frac{d}{2\delta}-\frac{C'}{\delta}-5}-2)\delta.
        \end{equation*}

        \item Suppose that $\exists z \in [x,y], d(z,[x_1,y_1]) \leq 2 \delta$. Then, according to the inequality \eqref{ineg-L-gen-3} of Lemma \ref{ineg-L-gen}, we have :
    \begin{align*}
        L & \geq (2^{\frac{K}{\delta}-3}-2)\delta \geq (2^{\frac{\frac{d}{2}-3\delta}{\delta}-3}-2)\delta \geq (2^{\frac{d}{2\delta}-6}-2)\delta \geq (2^{\frac{d}{2\delta}-\frac{C'}{\delta}-5}-2)\delta.
    \end{align*}
 
    \end{itemize}
    \item Suppose that $d > 2K +6\delta$. \\
    Let $t_0 \in [a,b]$ and $x_0=f(t_0)$. Let us consider :
    $$g_{t_0}(t)=d(x_0,f(t))-d(x_0,l)-d(f(t),l)-6\delta.$$ The map $g_{t_0}$ is continuous on $[t_0,b]$ and $g_{t_0}(t_0)=d(x_0,x_0)-d(x_0,l)-d(x_0,l)-6 \delta = -2d(x_0,l)-6\delta < 0$. Then if $g_{t_0}(b)\geq 0$, there exists, by the intermediate value theorem, $t_1 \in (t_0,b]$ such that $g_{t_0}(t_1) =0$. Since $g_a(b)=d(x,y)-d(x,l)-d(y,l)-6\delta=d-2K-6\delta \geq 0 $ by hypothesis, we can find a sequence of points $t_0,t_1,\hdots,t_n$ in $[a,b]$ such that, denoting  $x_i=f(t_i), d_i=d(x_i,x_{i+1})$ and $K_i=d(x_i,l)$ we have : 
    \begin{itemize}
        \item $t_0=a,t_n=b$ 
        \item $\forall i \in \{0,n-2\}, g_{t_i}(t_{i+1})=0$, namely $d_i=K_i+K_{i+1}+6\delta$
        \item $d_{n-1}<K_{n-1}+K_n+6 \delta$.
    \end{itemize}
    Let $L_i=\mathrm{length}(f([t_i,t_{i+1}]))$, we can use lemma \ref{L-petits-bouts} between $x_i$ and $x_{i+1}$ to show that $L_i \geq (2^{\frac{K}{\delta}-3}-2)\delta$ as soon as $i \in \{0,n-2\}$. This justifies in particular that the number $n$ of points cutting the path $f$ as above is indeed finite. Moreover, the hypothesis $d > 2K+6\delta$ implies that $n \geq 2$. \\
    Thus we can lower bound the length of the path $f$ on $[a,b]$ by : 
    \begin{equation*}
        L=\sum_{i=0}^{n-1} L_i \geq \sum_{i=0}^{n-2} L_i \geq (n-1) (2^{\frac{K}{\delta}-3}-2)\delta
   \end{equation*}
    which indeed gives us the first inequality of \eqref{L-d>K-n}. \\
    
    Now consider for all $i \in \{0,\hdots,n\}$, $p_i$ a projection of $x_i$ on $l$ (we choose in particular $p_0=x_1$ and $p_n=y_1$). Denote $d_{1,i}=d(p_i,p_{i+1})$. Since $d_i \leq K_i+K_{i+1}+6\delta$, for all $i \in \{0,\hdots,n-1\}$, the corollary \ref{d1-delta} states that $d_{1,i} \leq 18\delta$. Therefore :
    \begin{equation} \label{d1}
        d_1=d(x_1,y_1)=d(p_0,p_n)\leq \sum_{i=0}^{n-1} d(p_i,p_{i+1}) = \sum_{i=0}^{n-1} d_{1,i} \leq 18n\delta.
    \end{equation}
    Thus :
    \begin{align*}
        d & =d(x,y) \text{ by definition } \\
          & \leq d(x,x_1)+d(x_1,y_1)+d(y_1,y) \text{ by the triangle inequality } \\
          & = d_1+2K+2C \\
          & \leq 18n\delta +2K +2C \text{ by the inequality \eqref{d1}}
    \end{align*}
   which indeed gives us the second inequality of \eqref{L-d>K-n}. \\
   The inequality \eqref{L-d>K} is immediate using the two inequalities of \eqref{L-d>K-n}.
\end{itemize}
\end{proof}

\vspace{2cm}
\printbibliography

@article{bowditch_markoff_1998,
	title = {Markoff triples and quasifuchsian groups},
	volume = {77},
	issn = {00246115},
	url = {http://doi.wiley.com/10.1112/S0024611598000604},
	doi = {10.1112/S0024611598000604},
	language = {en},
	number = {3},
	urldate = {2021-03-11},
	journal = {Proceedings of the London Mathematical Society},
	author = {Bowditch, B.H.},
	month = nov,
	year = {1998},
	pages = {697--736},
}

@article{series_primitive_2019,
	title = {Primitive stability and {Bowditch}'s {BQ}-condition are equivalent},
	url = {http://arxiv.org/abs/1901.01396},
	language = {en},
	urldate = {2021-03-11},
	journal = {arXiv:1901.01396 [math]},
	author = {Series, Caroline},
	month = jan,
	year = {2019},
	note = {arXiv: 1901.01396},
}

@article{minsky_dynamics_2013,
	title = {On dynamics of {Out}({F}\_n) on {PSL}(2,{C}) characters},
	volume = {193},
	url = {http://arxiv.org/abs/0906.3491},
	language = {en},
	urldate = {2021-03-24},
	journal = {Israel J. Math.},
	author = {Minsky, Yair N.},
	year = {2013},
	note = {arXiv: 0906.3491},
	pages = {47--70},
}

@article{lee_bowditchs_2019,
	title = {Bowditch's {Q}-conditions and {Minsky}'s primitive stability},
	volume = {373},
	issn = {0002-9947, 1088-6850},
	url = {http://arxiv.org/abs/1812.04237},
	doi = {10.1090/tran/7953},
	language = {en},
	number = {2},
	urldate = {2022-11-24},
	journal = {Transactions of the American Mathematical Society},
	author = {Lee, Jaejeong and Xu, Binbin},
	month = oct,
	year = {2019},
	note = {arXiv:1812.04237 [math]},
	pages = {1265--1305},
}

@article{series_geometry_1985,
	title = {The geometry of markoff numbers},
	volume = {7},
	issn = {0343-6993},
	url = {http://link.springer.com/10.1007/BF03025802},
	doi = {10.1007/BF03025802},
	language = {en},
	number = {3},
	urldate = {2022-11-24},
	journal = {The Mathematical Intelligencer},
	author = {Series, Caroline},
	month = sep,
	year = {1985},
	pages = {20--29},
}

@article{gilman_enumerating_2011,
	title = {Enumerating palindromes and primitives in rank two free groups},
	volume = {332},
	issn = {00218693},
	url = {https://linkinghub.elsevier.com/retrieve/pii/S0021869311000871},
	doi = {10.1016/j.jalgebra.2011.02.010},
	language = {en},
	number = {1},
	urldate = {2022-11-24},
	journal = {Journal of Algebra},
	author = {Gilman, Jane and Keen, Linda},
	month = apr,
	year = {2011},
	pages = {1--13},
}

@article{tan_generalized_2008,
	title = {Generalized {Markoff} maps and {McShane}'s identity},
	volume = {217},
	issn = {00018708},
	url = {https://linkinghub.elsevier.com/retrieve/pii/S0001870807002630},
	doi = {10.1016/j.aim.2007.09.004},
	language = {en},
	number = {2},
	urldate = {2022-11-24},
	journal = {Advances in Mathematics},
	author = {Tan, Ser Peow and Wong, Yan Loi and Zhang, Ying},
	month = jan,
	year = {2008},
	pages = {761--813},
}

@book{coornaert_geometrie_1990,
	address = {Berlin, Heidelberg},
	series = {Lecture {Notes} in {Mathematics}},
	title = {Géométrie et théorie des groupes},
	volume = {1441},
	isbn = {978-3-540-52977-4 978-3-540-46294-1},
	url = {http://link.springer.com/10.1007/BFb0084913},
	urldate = {2022-11-25},
	publisher = {Springer},
	author = {Coornaert, Michel and Delzant, Thomas and Papadopoulos, Athanase},
	year = {1990},
	doi = {10.1007/BFb0084913},
}

@article{series_primitive_2020,
	title = {Primitive stability and the {Bowditch} conditions revisited},
	url = {https://arxiv.org/abs/2006.10403},
	language = {en},
	urldate = {2022-11-25},
	journal = {arXiv:2006.10403 [math]},
	author = {Series, Caroline},
	month = jun,
	year = {2020},
	note = {Number: arXiv:1901.01396
arXiv:1901.01396 [math]},
}

@book{farb_primer_2012,
	address = {Princeton, NJ},
	series = {Princeton mathematical series},
	title = {A primer on mapping class groups},
	isbn = {978-0-691-14794-9},
	language = {en},
	publisher = {Princeton University Press},
	author = {Farb, Benson and Margalit, Dan},
	year = {2012},
}

@article{goldman_topological_1988,
	title = {Topological components of spaces of representations},
	volume = {93},
	issn = {0020-9910, 1432-1297},
	url = {http://link.springer.com/10.1007/BF01410200},
	doi = {10.1007/BF01410200},
	language = {en},
	number = {3},
	urldate = {2023-04-30},
	journal = {Inventiones Mathematicae},
	author = {Goldman, William M.},
	month = oct,
	year = {1988},
	pages = {557--607},
}

@article{canary_dynamics_2015,
	title = {Dynamics on character varieties: a survey},
	volume = {Vol. II},
	language = {en},
	journal = {Handbook of Group Actions},
	author = {Canary, Richard D},
	year = {2015},
	pages = {175--200},
}

@article{palesi_dynamique_2014,
	title = {Dynamique de l’action du groupe modulaire et triplets de {Markov}},
	volume = {31},
	issn = {2118-9242},
	url = {https://proceedings.centre-mersenne.org/articles/10.5802/tsg.298/},
	doi = {10.5802/tsg.298},
	language = {fr},
	urldate = {2023-05-01},
	journal = {Séminaire de théorie spectrale et géométrie},
	author = {Palesi, Frédéric},
	year = {2014},
	pages = {137--161},
}

@article{maloni_type-preserving_2021,
	title = {On type-preserving representations of thrice punctured projective plane group},
	volume = {119},
	issn = {0022-040X},
	url = {https://projecteuclid.org/journals/journal-of-differential-geometry/volume-119/issue-3/On-type-preserving-representations-of-thrice-punctured-projective-plane-group/10.4310/jdg/1635368618.full},
	doi = {10.4310/jdg/1635368618},
	number = {3},
	urldate = {2023-05-01},
	journal = {Journal of Differential Geometry},
	author = {Maloni, Sara and Palesi, Frédéric and Yang, Tian},
	month = nov,
	year = {2021},
	note = {Publisher: Lehigh University},
	pages = {421--457},
}

@article{maloni_character_2020,
	title = {On the character variety of the three–holed projective plane},
	volume = {24},
	issn = {1088-4173},
	url = {https://www.ams.org/ecgd/2020-24-04/S1088-4173-2020-00349-5/},
	doi = {10.1090/ecgd/349},
	language = {en},
	number = {4},
	urldate = {2023-05-01},
	journal = {Conformal Geometry and Dynamics of the American Mathematical Society},
	author = {Maloni, Sara and Palesi, Frédéric},
	month = mar,
	year = {2020},
	pages = {68--108},
}

@article{maloni_character_2015,
	title = {On the character variety of the four-holed sphere},
	volume = {9},
	issn = {1661-7207},
	url = {https://ems.press/doi/10.4171/ggd/326},
	doi = {10.4171/GGD/326},
	language = {en},
	number = {3},
	urldate = {2023-05-01},
	journal = {Groups, Geometry, and Dynamics},
	author = {Maloni, Sara and Palesi, Frédéric and Tan, Ser Peow},
	year = {2015},
	pages = {737--782},
}

@article{delzant_displacing_2011,
	title = {Displacing representations and orbit maps},
	volume = {Univ. Chicago Press},
	language = {en},
	journal = {Geometry, rigidity and group actions},
	author = {Delzant, Thomas and Guichard, Olivier and Labourie, Francois and Mozes, Shahar},
	year = {2011},
	pages = {494--514},
}

@book{knapp_lie_1996,
	address = {Boston, MA},
	title = {Lie {Groups} {Beyond} an {Introduction}},
	isbn = {978-1-4757-2455-4 978-1-4757-2453-0},
	url = {http://link.springer.com/10.1007/978-1-4757-2453-0},
	language = {en},
	urldate = {2023-06-02},
	publisher = {Birkhäuser Boston},
	author = {Knapp, Anthony W.},
	year = {1996},
	doi = {10.1007/978-1-4757-2453-0},
}

@article{nielsen_uber_1918,
	title = {Über die {Isomorphismen} unendlicher {Gruppen} ohne {Relation}},
	volume = {79},
	issn = {0025-5831, 1432-1807},
	url = {http://link.springer.com/10.1007/BF01458209},
	doi = {10.1007/BF01458209},
	language = {de},
	number = {3},
	urldate = {2023-06-02},
	journal = {Mathematische Annalen},
	author = {Nielsen, J.},
	month = sep,
	year = {1918},
	pages = {269--272},
}

@article{nielsen_isomorphismengruppe_1924,
	title = {Die {Isomorphismengruppe} der freien {Gruppen}},
	volume = {91},
	issn = {0025-5831, 1432-1807},
	url = {http://link.springer.com/10.1007/BF01556078},
	doi = {10.1007/BF01556078},
	language = {de},
	number = {3-4},
	urldate = {2023-06-02},
	journal = {Mathematische Annalen},
	author = {Nielsen, Jakob},
	month = sep,
	year = {1924},
	pages = {169--209},
}

\end{document}